\documentstyle{amltd2004}
\begin{document}
\annalsline{155}{2002}
\received{December 29, 1999}
\startingpage{709}
\def\bye{\end{document}}
 \font\tenrm=cmr10
\input amssym.def
\input amssym.tex
%--------------- Author macros ---------------
%for Bbb in amstex
\catcode`\@=11
\font\twelvemsb=msbm10 scaled 1100
\font\tenmsb=msbm10
%\font\ninemsb=msbm7 scaled 1100%msbm9
\font\ninemsb=msbm10 scaled 800
\newfam\msbfam
\textfont\msbfam=\twelvemsb  \scriptfont\msbfam=\ninemsb
  \scriptscriptfont\msbfam=\ninemsb
\def\msb@{\hexnumber@\msbfam}
\def\Bbb{\relax\ifmmode\let\next\Bbb@\else
 \def\next{\errmessage{Use \string\Bbb\space only in math
mode}}\fi\next}
\def\Bbb@#1{{\Bbb@@{#1}}}
\def\Bbb@@#1{\fam\msbfam#1}
\catcode`\@=12

 \catcode`\@=11
\font\twelveeuf=eufm10 scaled 1100
\font\teneuf=eufm10
\font\nineeuf=eufm7 scaled 1100%eufm9
\newfam\euffam
\textfont\euffam=\twelveeuf  \scriptfont\euffam=\teneuf
  \scriptscriptfont\euffam=\nineeuf
\def\euf@{\hexnumber@\euffam}
\def\frak{\relax\ifmmode\let\next\frak@\else
 \def\next{\errmessage{Use \string\frak\space only in math
mode}}\fi\next}
\def\frak@#1{{\frak@@{#1}}}
\def\frak@@#1{\fam\euffam#1}
\catcode`\@=12

 \def\bcw{\mathbin{\bigcirc\mkern-15mu\wedge}}

%THIS FILE CONTAINS A SPECIAL DIACRITIC MACRO TO
%ACCOMMODATE THE SPECIAL DIACRITICS IN THIS FILE
%WITHIN EQNARRAY ENVIRONMENT (MATH ENVIRONMENT)
% IT FOLLOWS:

\def\r#1{{\mathop{#1}\limits^\circ}}
\def\ring{\r}

%Note I designed the above macro to accommodate the o wedge

\def\nint{\mathbin{\int\mkern-18mu\diagup \;} }

\newcommand{\namelistlabel}[1] {\mbox{#1}\hfil}
\newenvironment{namelist}[1]{%
\begin{list}{}
{\let\makelabel\namelistlabel
\settowidth{\labelwidth}{#1}
\setlength{\leftmargin}{1.1\labelwidth}}
}{%
\end{list}}
%-------------- Author entries --------------------

\title{An equation of Monge-Amp\`ere type in\\  conformal geometry, and
four-manifolds\\ of positive Ricci curvature}
\shorttitle{Four-manifolds of positive Ricci curvature}  % Shortened version for headline title

 \acknowledgements{The research of the first author was supported in
part by NSF Grant DMS-9706864 and a Guggenheim Foundation Fellowship. The research of the second author was supported in
part by NSF Grant DMS-9801046 and an Alfred P. Sloan Foundation Research Fellowship. The research of the third author was
supported in part by NSF Grant DMS-9706507.}
  \twoauthors{Sun-Yung A. Chang, Matthew J. Gursky,}{Paul C. Yang}
 \institutions{Princeton University, Princeton, NJ\\
{\eightpoint {\it E-mail address\/}:   chang@math.princeton.edu}\\
\vglue6pt
University of Notre Dame, Notre Dame, IN\\
{\eightpoint {\it E-mail address\/}: mgursky@nd.edu}\\
\vglue6pt
 Princeton University, Princeton, NJ\\
{\eightpoint {\it E-mail address\/}: yang@math.princeton.edu}}

 %-------------- Article Text--------------------

\centerline{\bf Abstract}
 \vglue12pt
We formulate natural conformally invariant conditions on a 4-manifold for the existence of a metric whose Schouten tensor satisfies a 
quadratic inequality. This inequality implies that the eigenvalues of the Ricci tensor are positively pinched.

\vglue12pt
\intro

 Conformal geometry in two dimensions is distinguished by its relationship
to complex analysis.  In higher dimensions the landscape becomes more
complicated, and in the absence of some special structure (e.g., K\"ahler)
even an extensive knowledge of the theory of Riemann
surfaces is no longer a reliable guide.  

Our setting in this paper is four dimensions, and one of our goals is to
propose a point of view which emphasizes certain parallels between
conformal geometry in two and four dimensions.  To illustrate this, let
us begin by recalling the Gauss-Bonnet formula for compact surfaces:
\begin{equation}
2 \pi \chi (M^2) \, = \, \int_{M^2} \, K  \, d\mu  
\end{equation}
\noindent
where $ d\mu $ is the area element, and $K$ is the Gauss curvature of the
surface. In four dimensions the Chern-Gauss-Bonnet integrand is a quadratic 
polynomial in the curvature; nevertheless there is a strong analogy with (0.1).

To see this, let $(M^4, g )$ be a compact Riemannian
four-manifold, and let $W$, ${\rm Ric}$, and $R$ denote respectively the
Weyl curvature tensor, Ricci tensor, and scalar curvature of $g$.
To express the Chern-Gauss-Bonnet formula it will be helpful to introduce the
elementary symmetric functions $\sigma_k : \Bbb R^n \rightarrow \Bbb R $, 
$ 1\leq k \leq n$. Given a section $A$ of the bundle of symmetric
 two-tensors, we can use the metric to raise an index and view $A$ as a tensor 
of type (1,1), or equivalently as a section of ${\rm End}(TM^4)$. Under this
identification, $\sigma_k(A)$ means $\sigma_k$ applied to the eigenvalues of 
$A$.

In particular, let $A = {\rm Ric} - \frac{1}{6} R
g$ denote the Schouten tensor.  With the notation described above the 
Chern-Gauss-Bonnet formula can be written 
\begin{equation}
8 \pi^2 \chi (M^4) \, = \frac{1}{4} \int | W |^2 \, dv \, + \, \int \sigma_2 ( A ) dv ,
\end{equation} 
where $dv$ denotes the volume form.  When we compare (0.1) and (0.2), a certain
parallel emerges between the Gauss curvature of a surface and the quantity
$\sigma_2(A)$ of a four-manifold, despite the presence of the Weyl curvature 
term (0.2). Indeed, this term actually strengthens the analogy: recall that 
the Weyl tensor measures whether the four-manifold is locally conformally flat (LCF). But every surface is LCF, so the
obstruction is vacuous and the corresponding  term is absent in (0.1) (or if one prefers, it is zero).

A further parallel between $\int K \, d \mu$ and $\int \sigma_2 (A) dv$
is that both are conformally invariant.  This is obvious for the Gauss
curvature; for $\sigma_2 (A)$ it follows from (0.2) and the conformal
invariance of $\int |W|^2 dv$.

This is actually a special case of a more general phenomenon.  Let $(
M^{2k}, g )$ be a compact, LCF Riemannian manifold of dimension $n = 2k$. If we define $A = {\rm Ric} - \frac{1}{2 ( n - 1)} Rg = {\rm Ric} -
\frac{1}{2(2k-1)} Rg$, then the integral 
$$
\int \sigma_k ( A ) \, dv
$$
\noindent
is conformally invariant (see [V-1]).  Moreover,  

$$
\chi (M^{2k} ) \, = \, c_k  \, \int \sigma_k ( A ) \, dv .
$$

Returning to four dimensions, we have a further parallel between the
integrals in (0.1) and (0.2): if $\int_{M^2} K dv > \, 0$ then $M^2$ has
genus zero; on the other hand, if $( M^4 , g )$ has positive scalar
curvature and $\int \sigma_2 (A) dv > 0$, then the first Betti number
$b_1 (M^4) = 0$ (see [G-1]).

The foregoing observation provided the motivation for the main result of
the present paper.  To understand how, recall the classical vanishing
theorem of Bochner: a compact Riemannian manifold of positive Ricci
curvature has $b_1 = 0$.  The assumption on the Ricci curvature is
natural in light of the famous Weitzenbock formula for harmonic
one-forms $\omega$:
$$
\frac{1}{2} \, \Delta \, |\omega|^2 = \,| \nabla  \omega|^2 \, + {\rm Ric} (\omega, 
\omega).
$$
\noindent
It is a little surprising, therefore, that the positivity of conformal
invariants like $\int \sigma_2 (A) dv$ and the Yamabe invariant \pagebreak should
also imply the vanishing of~$b_1$.  So one is lead to conjecture:
Suppose $(M^4, g_0)$ is a compact four-manifold with $\int \sigma_2
(A_0) dv_0 > 0$ and Yamabe invariant $Y ( g_0 ) > 0$.  Is there a
conformal metric $g = e^{2w} g_0$ with strictly positive Ricci
curvature?

%Suppose we strengthened our assumption from $\int \sigma_2 (A_0) dv_0 >
%0$ to $\sigma_2 (A_0) > 0$ in a {\it pointwise} sense.  In this case,
%from the positivity of the Yamabe invariant it is easy to see that the
%scalar curvature $R_0$ must be everywhere positive (cf. Lemma 1.1).
%Moreover, the Ricci curvature of $g_0$ satisfies the following
%inequality (see Lemma 1.2):

To understand our approach to this conjecture, it will be illuminating to
point out a relationship between the Ricci tensor and $\sigma_2(A)$. 
Suppose $g = e^{2w} 
g_0$ with $\sigma_2 (A_g ) > 0$. From the positivity of the Yamabe invariant of 
$g_0$, it is not difficult to conclude that the scalar curvature $R$ of $g$ 
must 
also be positive (cf. Lemma 1.1). Moreover, the Ricci curvature of $g$ satisfies
the inequality  
$$
{\rm Ric} \,\geqslant \, \frac{3 \sigma_2 (A)}{R} g. $$
\noindent
In particular, ${\rm Ric}$ is positive.\footnote{In dimensions greater than four, 
$\sigma_2(A) >0$ implies that the scalar curvature has a sign, but not the 
Ricci curvature. In three dimensions $\sigma_2(A)>0$ actually implies a 
sign on the {\it sectional} curvature; see [GV].}

%eturning to our theme at the beginning, and keeping in
%ind the parallels between $\int K d \mu$ and $\int \sigma_2 (A) dv$,
%he connection between this problem and the classical uniformization
%heorem should be clear (except, of course, that our goal is not to
%roduce a conformal metric with constant $\sigma_2 (A)$, just one with
%onstant sign.  More about this later).  It would also seem that this
%roblem is more tractable than working with the Ricci tensor directly -
%fter all, we are trying to prescribe a scalar function (i.e., 
%\sigma_2 (A)$) through the choice of a scalar function (the conformal factor
%w$).  But this apparent reduction in dimension is actually deceiving,
%s we now attempt to explain.
%
The preceding tells us that we can solve our conjecture in the affirmative 
by constructing a conformal metric $g= e^{2w} g_0$ with $\sigma_2(A_g)>0$, 
assuming only that
$\int \sigma_2(A_0)dv_0>0$ and $Y(g_0)>0$.

We remark that the positivity
of $\sigma_2 (A)$ is a much stronger condition than positive Ricci curvature.
 Indeed, if we define 
$S = \, - {\rm Ric} + \frac{1}{2} R g $ (in general relativity this is the
gravitational tensor) then $S$ also satisfies the inequality
$$
S \,\geqslant \, \frac{3 \sigma_2 (A)}{R} \, g 
$$
\noindent
(see Lemma 1.2).  Thus, $\sigma_2 (A) > 0$ imposes a pinching condition
on the Ricci curvature; it implies that each eigenvalue of ${\rm Ric}$ is
positive, but less than the sum of the other three. Moreover, if
$(M^4 , g )$ is oriented with $\int \sigma_2 (A) dv > 0$, then the
Euler characteristic $\chi (M^4)$ and signature $\tau (M^4)$ must
satisfy the inequality
\begin{equation}
\chi (M^4) \, > \frac{3}{2} \, | \tau ( M^4 ) |
\end{equation}
\noindent
(see [G-2]).  However, there are many examples of four-manifolds
with positive Ricci curvature which violate (0.3) (see [ShYa]).  This is
discussed in more detail Section 8.

The main analytic difficulty arising in the study of $\sigma_2 (A)$ is 
the following. If we write $g = e^{2w} g_0$, then the tensor $A$
of $g$ is related to the tensor $A_0$ of $g_0$ by the identity
\begin{equation}
A \, = \, A_0 - 2 \nabla^2_0 w \, + \, 2 dw \, \otimes \, dw \, - |
dw|^2 \, g_0 ,
\end{equation}
\noindent
where $\nabla^2_0$ denotes   the Hessian with respect to $g_0$.  In
light of (0.4), the equation under consideration is
\begin{equation}
\sigma_2 (A_0 - 2 \nabla^2_0 w \, + \, 2 dw \, \otimes \, dw \, - |
dw|^2 \, g_0 ) \, = \, f > 0 .
\end{equation}
\noindent
This is an example of a fully nonlinear equation of Monge-Amp\`ere type.
Introducing local coordinates, one can view (0.5) as an equation of
the form
\begin{equation}
F[\partial _i\partial _j w, \partial _kw,w, x]=f
\end{equation}
where
$$
F: \Bbb R ^{n\times n} \times \Bbb R ^n \times \Bbb R \times \Bbb R ^n
\rightarrow
\Bbb R,
$$
given by $F(r_{ij}, v_k, s, x) = f$.
Then (0.6) is {\it elliptic} at a solution $w$ if the matrix
$\frac{\partial F}{\partial {r_{ij}}}$
is positive definite. In the case of (0.5), this will hold provided $f >0$
(see Proposition 1.5).
The classical
techniques for analyzing such equations usually begin by assuming that one has some kind of approximate solution
(sub-, super-, viscosity, etc.) which is elliptic.  As we shall
see, equation (0.5) is elliptic at $w$ if and only if the tensor $S = -
{\rm Ric} + \frac{1}{2} R g$ is positive definite for the metric $g = e^{2w}
g_0$.  That is, we are confronted with the difficulty of solving
(0.5) without knowing {\it a priori} that our linearized operator is
elliptic.

Despite these difficulties, we are able to resolve the conjecture in the
affirmative: 

\nonumproclaim{Theorem A}  Let $(M^4 , g_0 )$ be a compact
four\/{\rm -}\/manifold satisfying {\rm (i)} $\int \sigma_2 (A_0) d v_0 > 0$ and {\rm
(ii)} $Y ( g_0 ) > 0${\rm .}  Then there is a conformal metric $g = e^{2w} g_0$ with
$\sigma_2 (A_g )> 0${\rm .}\endproclaim

\nonumproclaim{{C}orollary B}  Under the assumptions of
Theorem {\rm A,} there is a conformal metric $g = e^{2w} g_0$ with 
\begin{itemize}
\item[{\rm (i)}] ${\rm Ric} > 0${\rm ,}
\item[{\rm (ii)}] $S \, = \, - \, {\rm Ric} \, + \, \frac{1}{2} R g >
0${\rm .}
\end{itemize}
\endproclaim

It is natural to ask for which manifolds the assumptions of the theorem are likely to
hold.  This question is addressed in Section 8, where we use Freedman's [F]
work to give a list of the simply connected candidates (up to
homeomorphism).  In addition, we construct some explicit examples.

Assuming $M^4$ is orientable, by combining the Chern-Gauss-Bonnet
formula
$$
8 \pi^2 \, \chi ( M^4) \, =\,\frac{1}{4}
\int | W |^2 \, dv \, + \, \,
\int \sigma_2 (A) \, dv
$$
\noindent
with the signature formula
$$
12 \pi^2 \tau (M^4) \, =\, \frac{1}{4}
\int ( | W^+ |^2 \, - \, | W^- |^2 ) \, dv 
$$
\noindent
we obtain
$$
2 \pi^2 ( 2 \chi (M^4) \, + \, 3 \tau (M^4) ) \, = \,\frac{1}{4}
\int | W^+ |^2 \, dv \, + \, \frac{1}{2} \, \int \sigma_2 (A) \,
dv .
$$
\noindent
Therefore, an equivalent formulation of Theorem A is
 
\nonumproclaim{Corollary C}  If $(M^4 , g_0)$ is an
oriented compact four\/{\rm -}\/manifold satisfying $Y ( g_0 ) > 0$ and  
$$
\frac{1}{4}\int | W^+_0 |^2 \, dv_0 \, < \, 2 \pi^2 ( 2 \chi (M^4) \, + \, 3
\tau (M^4)),
$$
then there is a conformal metric $g = e^{2w} g_0$ with $\sigma_2 (A_g) >
0$ {\rm (}\/hence with ${\rm Ric} > 0$, $S > 0$\/{\rm ).}
\endproclaim

The
problem of conformally deforming a metric with $\sigma_2 (A) > 0$ to one
with $\sigma_2 (A) \, \equiv$ constant is addressed - but not resolved -
in [V-2], where degree-theoretic arguments are used.  What is lacking
are $L^\infty$-estimates for solutions of (0.5).  In a subsequent paper
we present an alternative approach, including\break {\it a priori}
$L^\infty$-bounds for solutions of (0.5) on manifolds that are not
conformally equivalent to the round four-sphere [CGY-2].

We conclude the introduction with some remarks about the structure of
the proof and the organization of the paper.

To overcome the lack of ellipticity for the linearized problem, we will
regularize our equation by a geometrically natural fourth order term.
The regularized equation actually arises in spectral theory, in the
context of the zeta functional determinant of a conformally invariant
operator (see Section~2).  More precisely, our regularized equation is 
\begin{equation}
\sigma_2 (A) \, = \, \frac{\delta}{4} \, \Delta R \, - \, 2 \gamma_1 \,
| \eta |^2
\end{equation}
\noindent
where $\Delta R$ is the Laplacian of the scalar curvature, $\delta > 0$
is small, $\gamma_1 < 0 $ is a fixed constant, and $\eta$ is a
nowhere-vanishing section of the bundle of symmetric two-tensors.  For
each sufficiently small $\delta > 0$, we are able to show that (0.7)
admits a smooth solution with positive scalar curvature (see Section 4).

The next (and most involved) step is obtaining {\it a priori} estimates for
solutions of (0.7) that are independent of $\delta$.  This is
accomplished in Sections 3,5, and 6 - at least up to a point.  There
seem to be technical obstructions preventing us from establishing
anything beyond $C^{1,\alpha}$-estimates.  But these estimates are adequate to
prove that the regularizing term $\delta \Delta R$ in (0.7) is
approaching zero in (roughly) an $L^2$-sense as $\delta \rightarrow 0$.

The final step of the proof is an application of heat equation
techniques.  Using the Yamabe flow, we show that solutions to (0.7) can
be perturbed to give metrics with $\sigma_2 (A) > 0$, once $\delta$ is
sufficiently small. This is explained in Section 7. 

\demo{Acknowledgements} The authors wish to express their appreciation to
the referees for their time and helpful comments. This paper was completed while the second author was a Visiting
Fellow and the third author was a Visiting Professor at Princeton University. They
would like to thank the Department of Mathematics for its support and hospitality.
\enddemo

\section{Background}

In this section we establish our notation, sketch some background
material, and prove some technical lemmas.

\demo{{\rm a.}   The curvature of four-manifolds}
To begin, let $(M^4 , g )$ be a compact four-manifold.  The curvature
tensor will be denoted $Rm$, and usually viewed as a $(0,4)$-tensor.  We
let $W$, ${\rm Ric}$, and $R$ denote respectively the Weyl curvature, Ricci
curvature, and scalar curvature of $g$.  There are various ways to
decompose the curvature tensor under the action of the orthogonal
group, depending on the basis one chooses.  If $E = {\rm Ric} - \frac{1}{4}
Rg$ denotes the trace-free Ricci tensor, then
\begin{equation} \,
Rm = W + \frac {1}{2} E \bcw g + \frac{1}{24} R\,g \bcw g
\end{equation}
\noindent
where $\bcw$ is the Kulkarni-Nomizu product (see [Be, 1.110]).
Alternatively, if $A = {\rm Ric} - \frac{1}{6} Rg$, then we
have the somewhat simpler decomposition
\begin{equation}
Rm \, = \, W \, + \, \frac{1}{2} \, A \bcw g .
\end{equation}
\noindent
In conformal geometry there are certain computational advantages to
working with (1.2) instead of (1.1).
\medskip
If $\chi(M^4)$ denotes the Euler characteristic of $M^4$, then by the
Chern-Gauss-Bonnet formula 
\begin{equation}
\chi (M^4) \, = \, \int \,{\rm  P f}f (Rm)
\end{equation}
\noindent
where ${\rm Pf} f (Rm)$ denotes the Pfaffian of the curvature (now viewed as a
matrix of Lie algebra-valued two-forms).  Using the basis in (1.2),
we can also   express (1.3) as
\begin{equation}
\chi (M^4) \, = \, \frac{1}{ 8 \pi^2} \, \int \left( \frac{1}{4}| W |^2 \, + \, \sigma_2 ( A ) \right) dv.
\end{equation}   
\noindent
If $M^4$ is oriented, let $* : \Omega^p ( M^4) \rightarrow \Omega^{4-p}
( M^4)$ denote the Hodge operator.  Then we have the splitting $\Omega^2
( M^4) = \Omega^2_+ ( M^4) \oplus \Omega^2_- ( M^4)$ into the
sub-bundles of self-dual and anti-self-dual two-forms.  This splitting
induces a decomposition of the Weyl curvature into $W^{\pm} :
\Omega^2_{\pm} ( M^4) \rightarrow \Omega^2_\pm (M^4)$, viewed as as
bundle endomorphism. Combining the signature formula
$$
12 \pi^2 \tau ( M^4) \, = \,\frac{1}{4} \int\left( | W^+ |^2 \, - \, | W^{-} |^2\right)
$$
\noindent
with (1.4) we obtain 
\begin{equation}
2 \pi^2 ( 2 \chi ( M^4) \, + \, 3 \tau ( M^4) ) \, = \,\frac{1}{4}
\int | W^+ |^2 \, + \, \frac{1}{2} \, \sigma_2 (A).
\end{equation}
It is clear from (1.4) and (1.5) that the positivity of $\sigma_2 (A)$
implies global topological information. But it also implies local
geometric information, as the following lemmas show.

\proclaimtitle{see [V-1, Lemma 23]}
\proclaim{Lemma}   $R^2 \geq 24 \sigma_2 ( A )$
with equality if and only if $E = 0${\rm .}  In particular{\rm ,} if $\sigma_2 ( A)
> 0$ on $M^4$ then either $R > 0$ or $R < 0$ on $M^4${\rm .}
\endproclaim

\demo{Proof}  This is immediate, since 
\vglue8pt
\hfill${\displaystyle
\sigma_2 ( A ) \, = \, - \, \frac{1}{2} \, | E |^2 \, + \, \frac{1}{24}
\, R^2 \, \leq \, \frac{1}{24} \, R^2 .
}$
\enddemo

\proclaim{Lemma}
Let $P \in M^4$ and $X \in T_P M^4$ be a tangent vector at $P${\rm .}  If the
scalar curvature $R$ of $g$ is positive at $P${\rm ,} then
\begin{eqnarray}
S ( X, X) & = & - \, {\rm Ric} (X , X ) \, + \, \frac{R}{2} \, g ( X , X )
 \\
& \geq & \frac{3}{R} \, \sigma_2 (A) \, g ( X , X )  
 \nonumber \\
{\rm Ric} ( X , X ) & \geq & \frac{3}{R} \, \sigma_2 ( A ) \, g ( X , X ) .
\end{eqnarray} 
\endproclaim

\demo{Proof}  To simplify notation we often denote
$g ( X , X ) = | X |^2 = \langle X , X \rangle$.  
In terms of the trace-free Ricci
tensor, 
\begin{equation}
S \, = \, - E \, + \, \frac{1}{4} \, Rg,
\end{equation}
\noindent
so that
$$
S ( X , X ) \, = \, - E ( X , X ) \, + \, \frac{1}{4} \, R g ( X , X ).
$$
\noindent
Since $E$ is trace-free, we have the sharp inequality $|E(X,X)| \leq
\frac{\sqrt{3}}{2} \, | E | \, | X |^2$ (see [SW, p.\ 234]).  Thus
\begin{eqnarray*}
S ( X , X ) & \geq & - \, \frac{\sqrt{3}}{2} \, | E | \, |X|^2 \, + \,
\frac{1}{4} \, R | X |^2 \\
& = & -2 \left( | E | {\sqrt\frac{3}{2R}} \; \right) \,
\left ( {\sqrt\frac{R}{8}} \; \right) \, | X |^2 \, + \, \frac{1}{4} \, R |
X |^2
\\
& \geq & - \left(
| E | \, {\sqrt\frac{3}{2R}} \; \right)^2 \, | X |^2 \, - \left(
{\sqrt\frac{R}{8}} \; \right)^2 \, | X |^2 \, + \, \frac{1}{4} \, R | X|^2
\\
& = & \left( - \frac{3}{2} \, \frac{|E|^2}{R} \, + \, \frac{1}{8} \, R
\right) | X |^2 \\
& = & \frac{3}{R} \, \sigma_2 \, ( A ) \, | X |^2 .
\end{eqnarray*}
The proof of (1.7) is essentially the same.  We begin with
\begin{equation}
{\rm Ric} \, = \, E \, + \, \frac{1}{4} \, Rg . 
\end{equation}
\noindent
Then
$$
{\rm Ric} ( X , X ) \, \geq \, - \frac{\sqrt{3}}{2} \, | E | \, | X |^2 \, +
\, \frac{1}{4} \, R | X |^2,
$$
\noindent
and we can argue as before. 
\enddemo

{\it Remark}.  It is worthwhile comparing (1.8)
and (1.9).  Recall that every symmetric two-tensor can be decomposed
into a trace-free part and a pure trace part.  The identities (1.8) and
(1.9) show that $S$ and ${\rm Ric}$ have the same pure trace component under
this decomposition, but their trace-free components differ by a sign. 
\vglue4pt

Arguing exactly as in the proof of Lemma 1.2 we have

\proclaim{Lemma}
Let $P \in M^4$ and $X \in T_P M^4${\rm .}  If $R < 0$ at $P$ then
\begin{eqnarray*}
S ( X , X ) & \leq & \frac{3}{R} \, \sigma_2 (A) \, g ( X , X ), \\
{\rm Ric} ( X , X ) & \leq & \frac{3}{R} \, \sigma_2 (A) \, g ( X , X ) .
\end{eqnarray*}
\endproclaim

Combining the preceding lemmas we conclude:

\proclaim{{C}orollary} If $\sigma_2 ( A ) > 0$ on $M^4$ then either $S > 0$ and
${\rm Ric} > 0$ on $M^4${\rm ,} or $S < 0$ and ${\rm Ric} < 0$ on $M^4${\rm ,} depending on the
sign of the scalar curvature {\rm (}\/which is necessarily constant by Lemma
{\rm 1.1). }
\endproclaim

 b.  {\it Conformal changes of metric}.
Now denote our four-manifold by $(M^4, g_0)$.  We will usually write
conformal metrics in the form $g = e^{2w} g_0$.  Also, metric-dependent
quantities which have $0$ as a subscript or superscript are understood
to be with respect to $g_0$, while those without are with respect to
$g$.  For example, $\nabla^2_0 \varphi$ denotes the Hessian of
$\varphi$ with respect to $g_0$ and $\Delta_0 \varphi = t r_{g_0} 
\nabla^2_0 \varphi$ the Laplacian; while $\nabla^2 \varphi$ and
$\Delta \varphi = {\rm tr}_g \nabla^2 \varphi$ denote the Hessian
and Laplacian with respect to $g$.

Of basic importance are the transformation laws for the various
components of the curvature tensor under a conformal change of metric:
\begin{eqnarray}
R & =&  e^{-2w} (R_0 - 6 \Delta_0 w - 6 | \nabla_0 w|^2 ), \\ 
{\rm Ric} & =& {\rm Ric}_0 \, - 2 \nabla^2_0 w - \Delta_0 wg_0 \, + \, 2 d w
\otimes dw \, - 2 | \nabla_0 w|^2 g_0 , \\
A & = & A_0 \, - 2 \nabla^2_0 w \, + \, 
2dw \otimes dw \, -| \nabla_0 w|^2 g_0 , \\
S & = & S_0 + 2 \nabla^2_0 w \, 
- 2 \Delta_0w g_0 \, 
- 2dw \otimes dw \, - | \nabla_0 w|^2 g_0 .
\end{eqnarray}
\noindent
It will often be useful to rewrite the above identities so that the
covariant derivatives are taken with respect to $g$ instead of $g_0$.
In this case, 
\begin{eqnarray}
R & = & R_0 e^{-2w} \, - 6 \Delta w \, + \, 6 | \nabla w|^2, \\
{\rm Ric} & = & {\rm Ric}_0 \, - 2 \nabla^2 w \, - \Delta w g \, - 2 dw
\otimes dw \, + \, 2 | \nabla w|^2 g, \\
A & = & A_0 \, - 2 \nabla^2 w \, - 2dw \otimes dw \, + \, | \nabla w|^2
g, \\
S & = & S_0 \, + \, 2 \nabla^2 w \, - 2 \Delta w g \, 
+ 2dw \otimes dw \, + \, | \nabla w|^2 g .
\end{eqnarray}
The Bach tensor plays a prominent role in our analysis.  It is defined
by (see [De])
$$
B_{ij} \, = \, 
\nabla^k 
\nabla^\ell 
W_{kij\ell} \, + \,
\frac{1}{2} \, 
R^{k \ell} \, W_{k i j \ell} .  $$
\noindent
Using the Bianchi identities, we can   rewrite this  as 
\begin{eqnarray}
B_{ij} &=& - \frac{1}{2} \Delta E_{ij} \, +  \frac{1}{6} \,
\nabla_i \nabla_j R \, - 
\frac{1}{24} \, \Delta R g_{ij} - E^{k \ell} 
W_{ikj\ell}  
\\
&&+\ E_i^k E_{jk} \, - \frac{1}{4} \, |E|^2 g_{ij} \, + \, \frac{1}{6} \,
R E_{ij} \nonumber 
\end{eqnarray}
where $\Delta E_{ij} = g^{k \ell} \nabla_k \nabla_\ell
E_{ij}$.  Although it has several interesting properties, for our
purposes the most important feature of the Bach tensor is its conformal
invariance: if $g = e^{2w} g_0$, then
\begin{equation}
B \, =\, e^{-2w} B_0 .
\end{equation}

\demo{{\rm c}.  Equations of Monge-Amp\`ere type}
Since our eventual goal is to produce conformal metrics with $\sigma_2
(A ) > 0$, it will be helpful to provide some background for the
analytic aspects of the problem.  If we fix a background metric $g_0$,
then by (1.12) we are attempting to solve the equation
\begin{equation}
\sigma_2 (A_0 \, - 2 \nabla^2_0 w \, + \,
2 dw \otimes dw \, - | \nabla_0 w |^2 g_0 ) \, = \, f
\end{equation}
\noindent
for some $f > 0$.  This is an example of a fully nonlinear equation of
Monge-Amp\`ere type (see [CNS-1], [CNS-2], [CKNS]).  Many of the relevant properties
of (1.20) are summarized by the following result: 

\proclaim{Proposition} The equation {\rm (1.20)} is elliptic at a solution $w$ if $f >
0${\rm .}  The linearized operator
$$
L[\varphi]= \frac{\partial F}{\partial r_{ij}} (\nabla _0^2 \varphi)_{ij}
$$
{\rm (}\/see the introduction\/{\rm )} is given by
\begin{equation}
L [ \varphi ] \, = \,
-2 S^{ij} \, \nabla_i^0 \, \nabla^0_j \, \varphi ,
\end{equation}
\noindent
where $S^{ij} \, = \, e^{-4w} (g_0)^{ik} \, (g_0)^{j\ell} \, S_{k \ell}${\rm ,} and
\begin{eqnarray*}
S_{k\ell}&= & S^0_{k \ell} \, + \,
2 \, \nabla^0_k \, \nabla^0_\ell w \, - 2 ( \Delta_0 w
)(g_0)_{k \ell}
\\
&&- \, 2 \nabla^0_k w \, \nabla^0_\ell w \, - |
\nabla_0 w|^2 ( g_0)_{k \ell} 
\end{eqnarray*}
is given by {\rm (1.13).}  If the scalar curvature $R$ of $g = e^{2w} g_0$ is
positive{\rm ,} then the ellipticity constants of $L$ satisfy
\begin{equation}
\frac{1}{2} Rf | \xi |^2 \, \geq \,
S_{ij} \xi^i \xi^j \, \geq \, \frac{3}{R} f \, | \xi |^2.
\end{equation}
\endproclaim 

A proof of Proposition 1.5 can be found in [V-1].  We only remark that
the estimates (1.22) follow from Lemma 1.2.
%%%%%%%%%%%%%%%%%%%%%%%%%%%%%%%%%%%%%%%%%%%%%%
 
\vglue-8pt
\section{The functional determinant}
\vglue-4pt
 
Let $(M^4, g_0)$ be a compact four-manifold.  A metrically defined
differential operator $L$ is said to be conformally covariant of bidegree
$(a,b)$ if under the conformal change of metric $g = e^{2w} g_0$,
\begin{equation}
L_g (\varphi) \, = \, e^{-bw} L_0 (e^{aw} \varphi ) .
\end{equation}
\noindent
In [BO] an explicit formula for $F [ w ] = \log ( \det L_g / \det L_0 )$
is computed, which may be expressed as
\begin{equation}
F [ w ] \, = \, \gamma_1
{\rm I} [ w] \, + \, \gamma_2 \, {\rm II} \, [ w ] \, + \,
\gamma_3 \, {\rm III} \, [ w ]
\end{equation}
\noindent
where $\gamma_i = \gamma_i (L)$ are constants and
\begin{eqnarray}
\qquad {\rm I}  \, [ w ] &\hskip-6pt = \hskip-6pt &\int 4 | W_0 |^2 w d v_0 \, - \left(
\int | W_0|^2 dv_0 \right) \, \log \, \nint \ e^{4w} d v_0 ,
\\
{\rm II} \, [w] &\hskip-6pt = \hskip-6pt & \int w P_0 w \, d v_0 \, + \, \int 4 Q_0 w \,
dv_0 \, -
\left(
\int Q_0 \, d v_0
\right) \,
\log \, \nint \, e^{4w} \, d v_0 , \nonumber\\
{\rm III}  \, [w] &\hskip-6pt = \hskip-6pt & 12 \left(
Y [ w ] \, - \frac{1}{3} \, \int \Delta_0 R_0 w \, d v_0 \right), \nonumber \\ 
Y [ w ] &\hskip-6pt = \hskip-6pt & \int \left(
\Delta_0 w \, + \, | \nabla_0 w |^2 \right)^2 \, dv_0 \, -
\frac{1}{3} \, \int R_0 | \nabla_0 w |^2 \, d v_0 . \nonumber
\end{eqnarray}
%ATTACHMENT begins here
Here $P$ denotes the {\it Paneitz} operator [P]:
$$
P = (\Delta )^2 + d^* \left(\frac {2}{3} R g -2 {\rm Ric}\right) d ,
$$
where d is the exterior derivative, $d^*$ is the adjoint of $d$, and $Q$
is the fourth order
curvature invariant:
\begin{eqnarray*}
Q&= &\frac {1}{12} \left(- \Delta R + \frac {1}{4} R^2 - 3 |E|^2\right).
\\
\noalign{\noindent Thus}
Q&=& \frac {1}{2} \sigma_2 (A)\,+ \frac {1}{12} (- \Delta R).
\end{eqnarray*}
Before we discuss the existence theory some remarks are in order,
explaining the significance of these formulas.  First, if we consider
the functional II alone, then critical points satisfy
$$
  P_0 w \, +  \, 2 Q_0 = \, 2 \left(
\int Q_0 \, dv_0 \right) \,
e^{4w} .
$$
\noindent
In general, if $g = e^{2 w} g_0$ is a conformal change of metric, then
the quantity $Q$ transforms according to the formula
$$
 P_0 w \, + \, 2 Q_0\, = \, 2 Qe^{4w}
$$
\noindent
where $Q = Q (e^{2w} g_0 )$.  We therefore conclude that critical points
of II are precisely those metrics which satisfy $Q \equiv$ constant. 

To understand III, it is helpful to rewrite it.  Let $R$ and $dv$ denote
the scalar curvature and volume form of the metric $g = e^{2w} g_0$;
then 
$$
{\rm III} \, [w] \, = \, \frac{1}{3} \,
\left[
\int R^2 \, d v \, - \int R^2_0 \, d v_0 \right]. $$
\noindent
From this expression it is easy to see that critical points of III
satisfy $\Delta R \equiv$ constant.  Since $M^4$ is compact, this
implies that $R$ is constant.  Thus III is the quadratic version of the
Yamabe functional.

In part, the interest of the functional determinant resides in the fact
that it is a natural Lagrangian arising in spectral theory whose Euler 
equation combines these geometrically natural ``sub-functionals.''

In order to state the relevant existence result of [CY-1] we need to
 define further the conformal invariant
\begin{equation}
\kappa_d \, = \,
\gamma_1 \, \int | W_0|^2 \, dv_0 \, + \, \gamma_2 \, \int Q_0 \, dv_0 .
\end{equation}

\proclaimtitle{[CY-1, Th.~1.1]}
\proclaim{Theorem}
  Let $( M, g_0)$ be a compact
four\/{\rm -}\/manifold{\rm .}  If $\gamma_2 , \gamma_3 > 0$ and $\kappa_d < 8 \gamma_2
\pi^2${\rm ,} then {\rm inf} $F (w)$ is attained by some function $w \in W^{2,2}$
and the metric $g = e^{2w} g_0$ satisfies
\begin{equation}
\gamma_1 | W |^2 \, + \, \gamma_2 Q \, - \gamma_3 \Delta R \, = \,
\kappa_d \, {\rm vol} ( g )^{-1} .
\end{equation}
Furthermore{\rm ,} $g$ is smooth {\rm ([CGY-1]).}
\endproclaim

{\it Remark}.
We warn the reader that the notation of [CY-1] is different: the signs of
$\gamma _i$ are reversed. What is crucial is that $\gamma _2$ and $\gamma _3$
have the same sign. If they have opposite signs then the existence theory
of [CY-1] is not applicable. There are examples arising in applications in
other contexts in which  $\gamma _2$ and $\gamma _3$ have opposite signs;
see [Br].
\vglue4pt

It will suit our purposes to modify slightly  the functional studied in\break
[CY-1] and [CY-2].  To describe our variant, we begin by pointing out that the
functional I in (2.2) does not involve any derivatives of the
conformal factor~$w$.  In addition, I introduces the term $\gamma_1 |
W|^2$ into the Euler equation~(2.5) --- or more precisely the term
$\gamma_1 e^{-4 w} | W_0|^2$.  In fact, any geometric quantity which
transforms in the same manner as the Weyl curvature under a conformal
change of metric behaves similarly.  Let us illustrate this with a
specific example.

Let $S_2 (M^4) = \Gamma ( {\rm Sym}( T^* M^4 \otimes T^* M^4 ))$ denote sections of
the bundle of symmetric $( 0 , 2 )$-tensors on $M^4$.  
Then $T^* M^4 \otimes T^* M^4$ 
inherits a bundle metric in the usual way from $TM^4$.
Moreover, if 
$g = e^{2w} g_0$ 
and $\eta \in S_2 (M^4)$ then $| \eta |^2_g
= e^{- 4w} | \eta |^2_0$.  In particular, this example enjoys the same
conformal scaling properties as the norm of the Weyl curvature:
\begin{eqnarray*}
|W|^2 & = & | W |^2_g \\
      & = & e^{-4w} | W_0 |^2_{g_0} \\
      & = & e^{-4w} | W_0 |^2 .
\end{eqnarray*}
\noindent
Analogous to (2.2) we can introduce the functionals
\begin{eqnarray}
I [w] & =&\, \int 4 | \eta |^2_0 \, w \, dv_0 - 
\, \left( 
\int | \eta |^2_0 dv_0 \right) \log \nint \, e^{4w} dv_0,
\\
F [w]&= & \gamma_1 {\rm I} [ w ] \, + \, \gamma_2 {\rm II} [ w ] \, + \,
\gamma_3 \, {\rm III} [ w ] ,
\end{eqnarray}
where II and III are defined as in (2.3).  We then have the
corresponding existence result:
 
\proclaim{{C}orollary}
Let $(M^4 , g_0)$ be a compact four\/{\rm -}\/manifold and $\eta \in S_2 (M^4)${\rm .}
If $\gamma_2, \gamma_3 > 0$ and $\kappa_d \equiv \gamma_1 \int | \eta
|^2_0 dv_0 + \gamma_2 \int Q_0 \, dv_0 < 8 \gamma_2 \pi^2${\rm ,} then {\rm inf} $F$ is
attained by some function $w \in W^{2,2}${\rm .}  The metric $g = e^{2w} g_0$
satisfies the Euler equation
\begin{equation}
\gamma_1 | \eta |^2 + \gamma_2 Q \, - \, \gamma_3 \Delta R \, = \,
\kappa_d vol(g)^{-1}
\end{equation}  
and is moreover smooth{\rm .}
\endproclaim 

The proof of Corollary 2.2 is identical in its  details to the proof of
Theorem 2.1, and will therefore be omitted. 

What is the purpose of modifying our functional in this way?  It is
difficult to give a complete answer to this question in advance of the
description of our regularized problem.  Eventually, though, we will
construct a conformal metric with the property that $\sigma_2 (A)$ is
bounded below by a positive constant times $| \eta |^2$, plus error terms.  Now it is easy
to choose a section of $S_2 (M^4)$ 
which is nowhere vanishing on $M^4$, and this
means that $\sigma_2 (A)$ is positive.  But if the minimum of $\sigma_2
(A)$ depended instead on $| W|^2$, then we could not in general rule out
the possibility that $|W|^2$ (and hence $\sigma_2 (A)$) vanishes
somewhere.
\pagebreak
 
To this end, let us fix once and for all a section
$\eta \in S_2 (M^4)$ 
which is nowhere vanishing (this is always possible: just let
$\eta$ be an arbitrary Riemannian metric on $M^4$).  Let $\delta \in ( 0
, 1]$, and set 
\begin{eqnarray}
\gamma_1 & = & - \int Q_0 d v_0 \,  \Bigg/ \, \int | \eta |^2_0 \, d v_0
\\
& = & - \frac{1}{2} \, \int \sigma_2 (A_0) \, d v_0 \, \Bigg/ \, \int |
\eta |^2_0 \, d v_0 , \nonumber \\ 
\gamma_2 & = & 1 , \nonumber \\
\gamma_3 & = & \frac{1}{24} \, ( 3 \delta - 2 ) .\nonumber   
\end{eqnarray} 

\demo{{R}emark}  
When $\delta < \frac{2}{3}$, according to (2.9), $\gamma _3 < 0$.
Thus we are considering values of $\gamma _3$ for which the existence result of Corollary 2.2
is inapplicable. This will be explained below.

Notice that if $\int \sigma_2 ( A_0) dv_0 > 0$ then $\gamma_1 < 0$.
With this choice of $( \gamma_1  , \gamma_2 , \gamma_3 )$,
\begin{equation}
\kappa_d \, = \, \gamma_1 \, \int | \eta |^2_0 \, d v_0 \, + \,
\int Q_0 dv_0 \, = \, 0 .
\end{equation}
\noindent
To write down the corresponding functional, let us introduce the
quantity
\begin{eqnarray}
U^\delta_0 & = & U^\delta ( g_0) \\
& = & \gamma_1 | \eta |^2_0 \, + \, Q_0 \, - \frac{1}{24} \, (3 \delta
-2 ) \, \Delta_0 R_0 . \nonumber
\end{eqnarray}
\noindent
Then according to (2.3), (2.6), (2.7), and (2.9),
\begin{eqnarray}
F [ w ] \, = \, F_\delta [w] & = & \gamma_1 {\rm I} [w] \, + \, 
{\rm II} [w] \, + \,
\frac{1}{24} \, (3 \delta -2 ) \, {\rm III} [w] \\
& = & \int 4 U^\delta_0 \, w \, d v_0 \, + \, \int w P_0 w \, d
v_0 \nonumber \\
&& + \, \frac{1}{2} ( 3 \delta -2 ) \, Y [ w ] . \nonumber
\end{eqnarray}  
\noindent
Note that $F_\delta$ is scale-invariant; i.e., $F_\delta [w + c] \, = \,
F_\delta[w]$ for any constant $c$.  By (2.7) and (2.8) the corresponding Euler
equation for $F_\delta$ is 
$$
{\gamma_1 | \eta|^2 \, + \, Q \, - \, \frac{1}{24} \, (3\delta - 2) \,
\Delta R \, = \, 0  \hspace{.5in} (*)_\delta} $$
\noindent
which can be rewritten as either
\begin{eqnarray*}
{\delta \Delta R} &= & 8 \gamma_1 | \eta|^2 \, - \, 2 | E|^2 \, +
\,  \frac{1}{6} R^2  \hspace{.5in} (*)_\delta \\
\noalign{\vskip-4pt}
\noalign{\noindent 
or }
\noalign{\vskip-4pt}
{\sigma_2 (A)}& = & \frac{\delta}{4} \, \Delta R \, - \, 2 \gamma_1
| \eta|^2 .  \hspace{.95in} (*)_\delta
\end{eqnarray*}

The latter way of writing $(*)_\delta$ reveals the motivation for
introducing the functional $F_\delta$.  For if $\delta = 0$, then
$(*)_0$ becomes $\sigma_2 (A) = \, - 2 \gamma_1 | \eta |^2$.  Now recall
that $\int \sigma_2 (A_0) dv_0 > 0$ implies that $\gamma_1 < 0$, so
in this case we conclude that $\sigma_2 (A) > 0$.  This observation
suggests the following strategy: to construct a conformal metric with
$\sigma_2 (A) > 0$, it suffices to show that $F_\delta$ admits a
critical point when $\delta = 0$.  This approach, however, presents some
serious technical difficulties.  
In some sense $F_\delta$ actually degenerates
as $\delta \rightarrow 0$.  
One can see that this is the case by writing down
just the highest order terms in (2.12):
\begin{eqnarray*}
F_\delta [w] \, = \,
\int \frac{3}{2} \delta ( \Delta_0 w)^2 & + & 
 ( 3 \delta -2 ) \, \Delta_0 w | \nabla_0 w|^2 \, + \, \frac{1}{2} ( 3
\delta -2 ) | \nabla_0 w|^4 \\
& + & {\rm (lower \ order \ terms)}.
\end{eqnarray*}
\noindent
When $\delta = 0$ the leading term is absent.  This behavior is
reflected in the Euler equation for $F_\delta$: when $\delta \ne 0$ then
$(*)_\delta$ is fourth order in the metric, but only second order when
$\delta = 0$.

Instead of studying $F_0$ directly we instead rely on a limiting
argument.  That is, we begin by showing that for any sufficiently small
$\delta > 0$, $(*)_\delta$ admits a smooth solution with positive scalar
curvature.  Even when $\delta > 0$, though, things are hardly routine:
recall that once $\delta < \frac{2}{3}$ then $\gamma_3 < 0$ while
$\gamma_2 > 0$, so that the existence theory of [CY-1] does not apply.  The
next (and most involved) step is to obtain {\it a priori} estimates for
solutions of $(*)_\delta$ that are independent of $\delta$.  For
technical reasons that we will explain at the appropriate time, the
optimal estimates we can derive give $W^{2,s}$-bounds on solutions with $s<5$.
This is sufficient to apply heat equation techniques and obtain a smooth conformal metric with
$\sigma_2 (A) > 0$.

To lay the groundwork for our study of $(*)_\delta$, let us begin by
fixing a $\delta_0 \in ( 0 , 1 )$ and defining 
$$
S   =  
 \{ \delta \in [ \delta_0 , 1 ] | (*)_\delta   \hbox{ \rm admits  a
 smooth solution with   positive  scalar  curvature} \}.
$$
\noindent
In Section 2 we will use the continuity method to show that $S = [
\delta_0, 1]$.  Since $\delta_0$ is arbitrary, we will conclude that
$(*)_\delta$ always admits a smooth solution of positive scalar
curvature for any $\delta \in(0 , 1]$.  We end this section with a
preliminary result which uses the existence theory of [CY-1] for the
functional determinant in order to show that $S$ is nonempty.

\proclaim{Proposition}
If $\int \sigma_2 (A_0) dv_0 > 0$ and $Y (g_0) > 0${\rm ,} then $1 \in S${\rm .}
\endproclaim 

{\it Proof}. When $\delta = 1$, $\gamma_3 =
\frac{1}{24}$.  It follows from Corollary 2.2 that there is a smooth
extremal metric $g = e^{2w} g_0$ satisfying $(*)_1$.  In particular, 
$$
\Delta R \, = \, 8 \gamma_1 | \eta |^2 \, - 2 | E |^2 \, + \,
\frac{1}{6} R^2.
$$
\noindent
Also, $\int \sigma_2 (A_0 ) dv_0 > 0$ implies that $\gamma_1 < 0$.  Thus

\centerline{${\displaystyle
\Delta R \, \leq \, \frac{1}{6} R^2}
$}

\noindent
on $M^4$.  It follows from [G-1, Lemma 1.2] that $R > 0$ on $M^4$. 
\hfill\qed

%%%%%%%%%%%%%%%%%%%%%%%%%%%%%%%%%%%%%%%%%%%%%%%%%%%%%%
\section{The regularized equation --- {\it a priori} estimates}

 In this section, we will derive some {\it a priori} estimates for smooth
solutions of the regularized equation $(*)_\delta$.

Let $F_\delta$ denote the functional as defined in (2.12).  That is,
$F_\delta$ is the functional (2.7) with coefficients $\gamma_1, \gamma_2,
\gamma_3$ chosen as in (2.9).  The main result in this section is:
\proclaim{Theorem}
Suppose $g = e^{2w} g_0$ is a smooth solution of $(*)_\delta$ with
positive scalar curvature{\rm ,} normalized so that $\int w d v_0 = 0${\rm .}  Then
there exist constants $C_0 , C_1 , C_p$ etc{\rm .,} all depending only on
$g_0${\rm ,} so that
\begin{equation}
w \, \geq \, C_0,
\end{equation}
\begin{equation}
\int [ \delta ( \Delta_0w)^2 \, + \, | \nabla_0 w |^4 ] \, d v_0 \,
\leq \, C_1,\; and
\; 
 \int ( - \Delta_0 w ) \, | \nabla_0 w|^2 \, dv_0 \, \leq \, C_1.\quad
\end{equation}
\noindent
Moreover{\rm ,} for  any  real  number $\alpha${\rm ,}
\begin{equation}
\int e^{\alpha w} \, d v_0 \, \leq \, C_\alpha .
\end{equation}
Finally{\rm ,} for any  positive  integer  $p${\rm ,} and for $0 \, < \, \delta \, \leq \,
\frac{1}{3}${\rm ,}
\begin{equation}
\int | \nabla_0 w |^4 \, | w |^p \, dv_0 \, \leq \, C_p .
\end{equation}
\endproclaim

We begin the proof of Theorem 3.1 with an identity.

\proclaim{Lemma}  
Suppose $g = e^{2w}g_0$  is a solution of $(*)_\delta${\rm .}
Then for any $\varphi \in W^{2,2} (M^4)${\rm , }
%\endproclaim
\begin{eqnarray}
&& \int \frac{3}{2} \,
\delta \Delta_0 w \, \Delta_0 \varphi  +  
\frac{1}{2} 
( 3 \delta -2 ) 
\left[
\Delta_0 \varphi | \nabla_0 w |^2 
  +  2 \Delta_0 w
\langle \nabla_0 \varphi, \nabla_0 w \rangle_0   \right.\\
&&\hskip2in\left.
+   2 | \nabla_0 w|^2 \langle \nabla_0 \varphi, \nabla_0 w
\rangle_0 \right]\nonumber \\
&& \quad\qquad = \
\int - 2 U^\delta_0 \varphi \, +\,
2 {\rm Ric}_0 ( \nabla_0 \varphi, \nabla_0 w) \, + \,
\frac{1}{2} ( \delta -2 ) R_0 \langle \nabla_0 \varphi ,
\nabla_0 w \rangle . \nonumber  
\end{eqnarray}
\endproclaim 

{\it Remark}.  Although we implicitly assume in
the proof that $w$ is smooth, it follows from a standard limiting
argument that (3.5) is valid if $w \in W^{2,2} (M^4)$.   Indeed, we
shall take (3.5) as our definition of a (weak) $W^{2,2}$-solution of
$(*)_\delta$.

\demo{Proof}  From a straightforward
computation (cf.\  [CY-1, (1.8)], or [BO]), $w$ satisfies
\begin{equation}
0 \, = \,
2 U^\delta_0 \, + \,
P_0 w \, + \,
\frac{1}{2}\, ( 3 \delta - 2 ) \,
\{ b_1 (w) \, + \, b_2 (w) \, + \, b_3 (w) \}
\end{equation}
\noindent
where  
\begin{eqnarray}
b_1 ( w ) & = &
\Delta^2_0 w \, + \,
\frac{1}{3} \,
R_0 \Delta_0 w \, + \, \frac{1}{3} 
\, 
\langle \nabla_0 R_0, \nabla_0 w \rangle_0 ,  \\
b_2 ( w ) & = & \Delta_0 | \nabla_0 w|^2 \, - 
2 ( \Delta_0 w )^2 \, - 2 \langle \nabla_0 w , \nabla_0
( \Delta_0 w ) \rangle_0 , \nonumber \\
b_3 ( w ) & = & - 2 | \nabla_0 w|^2 \Delta_0 w \, - 2 \langle
\nabla_0 w , \nabla_0 | \nabla_0 w |^2 \rangle_0 .\nonumber
\end{eqnarray} 
\noindent
Therefore,
\begin{equation}
0 \, = \, \int 2 U^\delta_0 \varphi \, + \, 
\int \varphi P_0 w \, + \,
\frac{1}{2} \,
( 3 \delta - 2 ) \left\{
\int \varphi b_1 ( w ) \, + \,
\int \varphi b_2 ( w ) \, + \, \int \varphi b_3 ( w ) \right\} .
\end{equation}  
\noindent
Proceeding term by term, we begin with the definition of the Paneitz
operator to get 
\begin{eqnarray}\qquad
\int \varphi P_0 w & = & 
\int \varphi \left[
\Delta^2_0 w \, + \,
d^* \left(
\frac{2}{3} \, R_0 g_0 \, - 
2 {\rm Ric}_0 \right) d w \right]  \\
& = &
\int \Delta_0 \varphi \Delta_0 w \, + \,
\frac{2}{3} \, R_0 \, 
\langle \nabla_0 \varphi , \nabla_0 w \rangle_0 \,
- 2 \, {\rm Ric}_0 ( \nabla_0 \varphi, \nabla_0 w ) .\nonumber 
\end{eqnarray} 
\noindent
Using the definitions in (3.7) and integrating by parts we also have
\begin{eqnarray}\qquad
\int \varphi b_1 ( w ) & = & 
\int \varphi \Delta^2_0 w \, + \,
\frac{1}{3} \, \varphi R_0 \Delta_0 w \, + \,
\frac{1}{3} \, \varphi \langle \nabla_0 R , \nabla_0 w
\rangle_0  \\
& = & \int \Delta_0 \varphi \Delta_0 w \, - \frac{1}{3} \,
\varphi \langle \nabla_0 R, \nabla_0 w \rangle_0 \,
- \frac{1}{3} \, R_0 \,
\langle \nabla_0 \varphi , \nabla_0 w \rangle_0 \nonumber
\\
&& + \, \frac{1}{3} \, \varphi \, 
\langle \nabla_0 R_0 , \nabla_0 w \rangle_0 \nonumber \\
& = & \int \Delta_0 \varphi \Delta_0 w \, - 
\frac{1}{3} \, R_0 \langle \nabla_0 \varphi , \nabla_0 w
\rangle_0 ; \nonumber
\end{eqnarray} 
\begin{eqnarray} 
\noalign{\vskip-8pt} &&\\
\int \varphi b_2 ( w ) & = &
\int \varphi \Delta_0 | \nabla_0 w |^2 \, - 
2 \varphi ( \Delta_0 w )^2 \, -
2 \varphi \langle \nabla_0 w, \nabla_0 ( \Delta_0 w)
\rangle_0\nonumber   \\
& = &
\int \Delta_0 \varphi | \nabla_0 w |^2 \, - 2 \varphi (
\Delta_0 w )^2 \, + \,
2 \varphi ( \Delta_0 w )^2 \, + \,
2 \Delta_0 w \langle \nabla_0 \varphi, \nabla_0 w 
\rangle_0  
\nonumber \\
& = &
\int \Delta_0 \varphi | \nabla_0 w |^2 \, + \,
2 \Delta_0 w \langle \nabla_0 \varphi, \nabla_0 w
\rangle_0 ; \nonumber
\end{eqnarray}
\begin{eqnarray}
\noalign{\vskip-8pt}
&&\\
\int \varphi b_3 ( w ) & = &
\int - 2 \varphi | \nabla_0 w |^2 \Delta_0 w \, - 
2 \varphi \langle \nabla_0 w , \nabla_0 | \nabla_0 w|^2 \rangle_0
 \nonumber\\
& = & 
\int - 2 \varphi | \nabla_0 w |^2 \Delta_0 w \, + \,
2 \varphi \Delta_0 w | \Delta_0 w|^2 \, + \,
2 | \Delta_0 w |^2 \langle \nabla_0 \varphi , \nabla_0 w \rangle_0
\nonumber \\
& = & 
\int 2 | \nabla_0 w|^2 \, \langle \nabla_0 \varphi, \nabla_0 w
\rangle_0 . \nonumber
\end{eqnarray}
Substituting (3.9)--(3.12) into (3.8) we arrive at (3.5).
\enddemo

\demo{Proof of Theorem {\rm 3.1}} 
In the following, we let $C$ denote various constants whose
value may change but depend only on $g_0$.

\vglue4pt {\it Proof of} (3.1). If $R$ denotes the
scalar curvature of $g$, then $ R > 0$ by our assumption, and
by (1.10)  
\begin{equation}
\Delta_0 w \, + \, | \nabla_0 w |^2 \, + \, \frac{1}{6} \, R
e^{2w} \, = \, \frac{1}{6} R_0 .
\end{equation}
\noindent
Hence
\begin{equation}
\Delta_0 w \, + \, | \nabla_0 w |^2 \, \leq \, \frac{1}{6} \,
R_0 ,
\end{equation}
\noindent
and in particular
\begin{equation}
\Delta_0 w \, \leq \, \frac{1}{6} R_0.
\end{equation}
\noindent
Now by Green's formula, 
\begin{equation}
- \, w ( x ) \, + \, \bar{w} \, = \,
\int G ( x , y ) \, \Delta_0 w ( y ) \, dv_0 (y),
\end{equation}
\noindent
where $G (x,y)$ is the Green's function for $(M , g_0 )$, and
$\bar{w} = \int w dv_0 = 0$.  Since $M$ is a compact manifold,
we may add a constant to $G$ and assume it is positive. Then (3.1) 
follows from (3.15) and (3.16).

\vglue4pt {\it Proof of} (3.2). By integrating
(3.13) over $M^4$, we have
\begin{equation}
\int | \nabla_0 w |^2 \, \leq \, C.
\end{equation}
\noindent
Since $\int w = 0$, by the Poincar\'e inequality we conclude 
\begin{equation}
\int w^2 \, \leq \, C .
\end{equation}
Now, taking $\varphi = w$ in (3.5) we have 
\begin{eqnarray*}
&&\int \frac{3}{2} \,
\delta ( \Delta_0 w )^2   + 
\frac{3}{2} \, ( 3 \delta - 2 ) \, \Delta_0 w | \nabla_0 w|^2 \, + \,
( 3 \delta - 2 ) | \nabla_0 w |^4 \\
\\ 
&&\qquad\qquad =  
\int - 2 U_0^\delta w \, + \,
2 {\rm Ric}_0 ( \nabla_0 w, \nabla_0 w ) \, + \,
\frac{1}{2} ( \delta - 2 ) R_0 | \nabla_0 w|^2 .
\end{eqnarray*}
\noindent
Using (3.17) and (3.18) we conclude
\begin{eqnarray}\qquad\quad
\int \frac{3}{2} \, \delta ( \Delta_0 w )^2 \, + \, 
\frac{3}{2} \, ( 3 \delta - 2 ) \, \Delta_0 w | \nabla_0 w |^2 \, + \,
( 3 \delta- 2 ) \, | \nabla_0 w |^4 \, \leq \, C .
\end{eqnarray}
\noindent
There are now two cases to consider.  First, suppose that $\delta \in
\left[ \frac{2}{3} , 1 \right]$, i.e., that $3 \delta - 2 \in [ 0 , 1
]$.  It then follows from the inequality 
$$
\frac{3}{2} \, ( 3 \delta - 2 ) x y \, \geq \, - \, \frac{9}{16} \,
(3 \delta - 2 ) x^2 \, - \,
(3 \delta - 2 ) y^2
$$
\noindent
that 
$$
\int \frac{3}{16} \, ( 6 - \delta ) \, ( \Delta_0 w )^2 \, \leq \,
C $$
\noindent
which implies
\begin{equation}
\int ( \Delta_0 w )^2 \, \leq \, C .
\end{equation}
On the other hand, suppose $\delta \in ( 0 , \frac{2}{3} )$.  Then $3
\delta - 2 \in ( - 2 , 0 )$, and by (3.13)
$$
| \nabla_0 w|^2 \, \Delta_0 w \, + \, | \nabla_0 w|^4 \, \leq \,
\frac{1}{6} \, R_0 | \nabla_0 w |^2
$$
\noindent
whence
$$
\int ( 3 \delta - 2) \, \Delta_0 w | \nabla_0 w |^2
\, \geq \,
\int - \, (3 \delta - 2 ) \, | \nabla_0 w|^4 \, - C .
$$
\noindent
Substituting this into (3.19) gives 
\begin{equation}
\int \frac{3}{2} \, \delta ( \Delta_0 w )^2 \, -
\, \frac{1}{2} \, ( 3\delta - 2 ) \, | \nabla_0 w |^4 \, \leq \, C, 
\end{equation}   
\noindent
and
$$
\int \frac{3}{2} \, \delta ( \Delta_0 w)^2 \, - \,
\frac{1}{2} \, ( 3 \delta - 2) \, \Delta_0w | \nabla_0 w|^2 \,
\leq C . $$
\noindent
Finally, using (3.14) again we observe that
\begin{eqnarray*}
\int | \nabla_0 w|^4 & \leq & 
\int - \Delta_0 w | \nabla_0 w|^2 \, + \, \frac{1}{6} \, R_0 |
\nabla_0 w |^2 \\
& \leq & \left( \int ( \Delta_0 w )^2 \right)^{\frac{1}{2}} \,
\left( \int | \nabla_0w|^4 \right)^{\frac{1}{2}} \, + \, C,
\end{eqnarray*}
\noindent
which implies
\begin{equation}
\int | \nabla_0 w |^4 \, \leq \,
\int ( \Delta_0 w )^2 \, + \, C.
\end{equation}
To complete the proof of (3.2), notice that when $\delta \in \left[
\frac{2}{3}, 1 \right]$, then (3.2) follows from (3.20) and (3.22).
If $\delta \in [\frac 13, \frac 23) $, then (3.2) follows from
(3.22) and the first half of (3.21). If $ \delta \in (0, \frac 13]$,
then (3.2) follows from the first and second half of (3.21).
\enddemo

{\it Proof of {\rm (3.3)}}.
This is a direct consequence of (3.2) and the sharp Sobolev embedding
theorem of Moser [M] - Trudinger [T]: if $w \, \varepsilon W^{1,n}_0 ( \Omega )$,
then $w$ is in the Orliz class $e^{L^{ \frac{n}{n-1}}} ( \Omega )$ 
for any bounded domain
$\Omega$ in $\Bbb{R}^n$. In particular, $e^{\alpha w}$ is integrable for
any $\alpha$. Trudinger's result was later generalized ([F],
[BCY]) to functions $w \varepsilon W^{1,n} (M^n )$
with $\int w dv_0 = 0$ on $(M^n , g_0)$, with $(M^n, g_0) $ any compact
Riemannian manifold.

\vglue4pt {\it Proof of} (3.4).  We will prove the
statement inductively on $p$.  We first observe that by (3.1), we may
replace $w$ with $w+C_0$ and so assume that $w \geq 0$.  
We now substitute $\varphi =
w^p$ in (3.5) and call the expression $(3.5)_{p}$.

Integrating by parts on the left-hand side of $(3.5)_{p}$, we get for
$\delta < \frac{2}{3}$,  
\begin{eqnarray} \noalign{\vskip-4pt}
&&\\[-6pt]
{\rm L.H.S. \ of \ (3.5)}_p & = & \frac{3 \delta}{2} \, p
\, 
\int ( \Delta_0 w ) \, 
\left[ \Delta_0 w \, w^{p-1} \, + \, (p - 1) \,
| \nabla_0 w|^2 \, w^{p-2} 
\right] \nonumber \\
& & +\, \frac{1}{2} \,
( 3 \delta -2 )p \,
\int \left[ 
3 \Delta_0 w | \nabla_0 w |^2 \, w^{p-1} \,   \right.\nonumber \\ 
\nonumber \\
&&\,  + \, 2 \left. | \nabla_0 w|^4 \, w^{p-1} 
\, + \, (p-1) \, | \nabla_0 w |^4 \, w^{p-2} \right] \nonumber \\
& = & \frac{3 \delta}{2} \, p \, \int ( \Delta_0 w )^2 \, w^{p-1} \, + \,
\frac{3 \delta}{2} \, p ( p-1) \, \int ( \Delta_0 w ) \, | \nabla_0 w|^2 \, w^{p-2} \nonumber \\
& &+ \, \frac{1}{2} \, ( 2 - 3 \delta ) \, p \, 
\left[
2 \, \int ( - \Delta_0 w \, - \, | \nabla_0 w|^2 ) \, | \nabla_0 w|^2
\, w^{p-1} \right. \nonumber \\
&&
+ \,  \left. \int ( - \Delta_0 w ) \, | \nabla_0 w|^2 \, w^{p-1} 
\right] \nonumber \\
& &-  \, \frac{1}{2} \, ( p-1) \, ( 2 - 3 \delta ) \, \int |\nabla_0 w |^4 \, w^{p-2}\nonumber \\  
& \geq & I_p \, + \,
II_p \, - \,
\frac{1}{6} \,
( 2-3 \delta) \, p \, \int R_0 \, | \nabla_0 w |^2 \, w^{p-1}, \ \ \
{\rm (by \ (3.14)) } \nonumber 
\end{eqnarray} 
\noindent
where 
\begin{equation}
I_p \, = \,
\frac{3 \delta}{2} \, p \, \int ( \Delta_0 w )^2 \,
w^{p-1} \, + \,
\frac{1}{2} ( 2 - 3 \delta) p \, \int ( - \Delta_0 w ) \, | \nabla_0 w |^2 w^{p-1} ,
\end{equation}
\begin{equation}
II_p \, = \, \frac{3 \delta}{2} \, p (p-1) \, \int ( \Delta_0 w ) \, | \nabla_0 w |^2 \, w^{p-2} \, - \,
\frac{1}{
2} \,
(p - 1 ) \, ( 2 - 3 \delta ) \, \int | \nabla_0 w |^4 \, w^{p-2} .
\end{equation}
On the right-hand side of $(3.5)_p$ we have
\begin{eqnarray}
{\rm R.H.S. \ of \ (3.5)_p} & = & 
\int - 2 U^\delta_0 \, w^p \, + \,
2 p \, \int {\rm Ric}_0 \, | \nabla_0 w |^2 \, w^{p-1}  \\
&& + \, \frac{1}{2} \,
( \delta -2 ) p \, \int R_0 \, | \nabla_0 w |^2 \, w^{p-1} \nonumber
\\
& \lesssim & \int w^p \, + \,
p \left( \int | \nabla_0 w |^4 \right)^{1/2} 
\, \left( \int w^{2 ( p - 1)} \right)^{1/2}
\nonumber \\
& \lesssim & C_p  \ \ \ \ \ {\rm (by \ (3.2)). }\nonumber
\end{eqnarray} 
Combining (3.23) and (3.26), we conclude that
\begin{equation}
I_p \, + \, II_p \, \lesssim \, C_p .
\end{equation}
We now claim that 
\begin{equation}
- \, II_p \, \lesssim \,
C_p \, I_{p-1} \, + \,
C_p \, \leq \, C_p , \hspace{.5in} {\rm for} \ p \, \geq \, 2.
\end{equation}
To see (3.28), we first observe that from (3.14), we have 
\begin{eqnarray}
\int | \nabla_0 w |^4 \, w^{p-1} & \leq & \int \left(
\frac{1}{6} \, R_0 - \Delta_0 w \right) \, | \nabla_0 w |^2 \, w^{p-1}
  \\
& \leq & \int \left( - \, \Delta_0 w \right) \, | \nabla_0 w |^2 \,
w^{p-1} \, + \, C_p .\nonumber
\end{eqnarray}
\noindent
Thus for $p \geq 2$, $\delta < 2/3$,
\begin{eqnarray*}
- \, II_p & \lesssim & C_p \, \int ( - \Delta_0 w ) \, | \nabla_0 w
|^2 \, w^{p-2} \, + \, C_p 
\nonumber \\
& \lesssim & C_p \, I_{p-1} \, + \, C_p .
\end{eqnarray*}
\noindent 
When $p = 2$, $I_1$ is bounded via (3.2), thus $- II_2$ is bounded
and hence $I_2$ is bounded via (3.27).  Thus it is clear we can
establish (3.28) and (3.27) inductively for all $p \geq 2$. 
Also, note that the constant $C_p$ will depend on $(\frac{2}{3} - \delta)^{-1}$;
thus we assume $\delta \leq \frac 13 $ to eliminate this dependence.
\hfill\qed

\section{The regularized equation: Existence and regularity}

 In this section, we will show that for all sufficiently small
$\delta > 0$, $(*)_\delta$ admits a smooth solution with
positive scalar curvature.  To accomplish this, we will apply the continuity 
method.  Fix
$\delta_0 \in (0,1)$, and define 
$$
S   =  
 \{ \delta \in [ \delta_0 , 1 ] | (*)_\delta \ \hbox{\rm admits a
  smooth solution    with   positive   scalar  curvature} \}.
$$
Following the usual practice, we will show that $S = [ \delta_0 , 1]$ by
arguing that $S$ is both open and closed.  
Since we already saw that $ 1 \, \in S$ by Proposition 2.1, 
the desired result will follow.
\proclaim{Proposition} If $\int \sigma_2 (A_0) dv_0 > 0$ then $S$ is open{\rm .}
\endproclaim  

{\it Proof}. The proof of this fact relies (as
usual) on a perturbation result.  Consequently, we will need to study
the linearized problem.

\proclaim{Lemma}
Let ${\cal{L}}_\delta$ denote the 
linearization of $(*)_\delta$ at a
solution $g$ of positive scalar curvature{\rm .}  Then for any $\varphi \in
W^{2,2}${\rm ,}
\begin{equation}
\langle \varphi , {\cal{L}}_\delta \varphi \rangle_{L^2} \, \geq \,
\int \frac{3}{13} \, \delta^2 ( \Delta \varphi )^2 \, + \,
\frac{\delta}{16} \, R | \nabla \varphi|^2. 
\end{equation}
In particular{\rm ,} ${\rm Ker} \, {\cal{L}}_\delta = {\Bbb{R}}${\rm .}
\endproclaim

{\it Remark}. The kernel of ${\cal{L}}_\delta$ is due
to the scale-invariance of $F_\delta$.

\demo{Proof}  By a straightforward computation
(Theorem 2.1 in [CY]), 
\begin{equation}
\langle \varphi , {\cal{L}}_\delta \varphi \rangle_{L^2} \, = \,
\int 3 \delta ( \Delta \varphi )^2 \, - \,
4 E ( \nabla \varphi , \nabla \varphi ) \, + \,
( 1 - \delta ) \, R | \nabla \varphi|^2 .\hskip.5in
\end{equation} 
We start with the sharp inequality of [SW, p.\ 234]:
$$
\int - 4 E ( \nabla \varphi , \nabla \varphi ) \, \geq \,
\int - 4 \left( \frac{\sqrt{3}}{2} \right) \, | E | \, |
\nabla \varphi |^2 .
$$
By the inequality $2 x y$  $ \leq \, \varepsilon x^2 + \varepsilon^{-1}
y^2$, which holds for any $\varepsilon > 0$, it follows that 
\begin{equation}
\int - 4 E ( \nabla \varphi , \nabla \varphi ) \, \geq
\,
\int - 2 \varepsilon \left ( \frac{\sqrt{3}}{2} \right)^2 \,
\frac{|E|^2}{R} \, | \nabla \varphi |^2 \, - 2 \varepsilon^{-1} R |
\nabla \varphi|^2 . \hskip.35in
\end{equation}
\noindent
Since $g$ satisfies $(*)_\delta$,  
$$
- |E|^2 \, = \,
\frac{\delta}{2} \, \Delta R \, - \,
\frac{1}{12} \, R^2 \, - \, 4 \gamma_1 | \eta |^2 .
$$
\noindent
Also, $\int \sigma_2 (A_0 ) dv_0 > 0$ implies $\gamma_1 < 0$, so  that $ -
|E|^2 \geq \frac{\delta}{2} \, \Delta R \, - \frac{1}{12} R^2$.
Substituting this into (4.3) gives
\begin{equation}
\int - 4 E ( \nabla \varphi , \nabla \varphi ) \, \geq
\,
\int \frac{3 \varepsilon \delta}{4 } \,
\frac{\Delta R}{R} \, | \nabla \varphi|^2 \,
- \, \left( \frac{2}{\varepsilon} \, + \, \frac{\varepsilon}{8} \right) 
\, R | \nabla \varphi|^2.\hskip.35in
\end{equation}
\noindent
Integrating by parts in the first term on the right-hand side of (4.4)
we get
\begin{eqnarray*}
\int \frac{\Delta R}{R} \, | \nabla \varphi |^2 & = & 
\int - \,
\nabla R \, \nabla ( R^{-1} ) \, | \nabla
\varphi|^2 \, - \, \frac{\nabla R}{R} \, \nabla | \nabla \varphi|^2 \\ 
\\
& = & \int \frac{| \nabla R|^2}{R^2} \, 
| \nabla \varphi |^2 \, - \,
2 \, 
\nabla^2 \varphi 
\left( 
\nabla \varphi, \,
\frac{\nabla R}{R}
\right) .
\end{eqnarray*}
\noindent
From the inequality 
$ 2 | \nabla^2 \varphi ( 
\nabla \varphi ,
\frac{\nabla R}{R} ) 
| \, \leq \, 
| \nabla^2 \varphi |^2 \, + \, 
\frac{| \nabla R|^2}{R^2} \, 
| \nabla \varphi|^2 $,
\noindent
this becomes
$$
\int \frac{\Delta R}{R} \, | \nabla \varphi |^2 \, \geq \,
\int - 
| \nabla^2 \varphi |^2.
$$
\noindent
By the integrated Bochner formula, 
$$ \int - \, | \nabla^2 \varphi|^2 \, = \,
\int - ( \Delta \varphi )^2 \, + \,
E ( \nabla \varphi , \nabla \varphi ) \, + \,
\frac{1}{4} \, R | \nabla \varphi|^2 .
$$
\noindent
Therefore,
\begin{eqnarray*}
\int \frac{\Delta R}{R} \, | \nabla \varphi |^2 & \geq & 
\int - \,  
( \Delta \varphi )^2 \, + \, 
E ( \nabla \varphi , \nabla \varphi ) \\
\\ 
& &+ \, \frac{1}{4}  R | \nabla \varphi |^2 .
\end{eqnarray*}
\noindent
Substituting this into (4.4) gives
\begin{eqnarray*}
\int - 4 E ( \nabla \varphi , \nabla \varphi ) & \geq &
\int - \, 
\frac{3 \varepsilon \delta}{4} \,
( \Delta \varphi )^2 \, + \, 
\frac{3 \varepsilon \delta}{4} \, E 
( \nabla \varphi , \nabla \varphi ) \\
\\
&& + \, 
\left(
\frac{ 3 \varepsilon \delta}{16} \, - \, 
\frac{2}{\varepsilon} \, - \, 
\frac{\varepsilon}{8}
\right)  \, R | \nabla \varphi |^2 .
\end{eqnarray*}
\noindent
Now take $\varepsilon = \frac{4 (
4-\delta)}{4-3 \delta}$, which yields
\begin{eqnarray*}
\int - 4 E ( \nabla \varphi , \nabla \varphi ) & \geq &
\int \frac{- 12 \delta ( 4 - \delta )}{(16 - 3 \delta^2)} \, (
\Delta \varphi )^2 \\
\\
&& + \, 
\frac{(3 \delta^3 - 44 \delta^2 + 112 \delta - 64)}
{(4 - \delta ) (16 - 3 \delta^2 )} \,
R | \nabla \varphi |^2 .
\end{eqnarray*}
\noindent
Substituting this inequality into (4.2), we get 
$$
\langle \varphi, {\cal{L}}_\delta \varphi \rangle_{L^2} \, \geq \,
\int \frac{3 \delta^2 ( 4 - 3 \delta )}{(16 - 3 \delta^2)} \,
( \Delta \varphi )^2 \, + \,
\frac{\delta ( - 3 \delta^3 + 18 \delta^2 - 40 \delta + 32)}
{(4 - \delta ) (16-3 \delta^2)} \,
R | \nabla \varphi |^2 .
$$
\noindent
Since $\delta \in [ 0 , 1 ]$, we can estimate the above expressions to
arrive at 
\vglue8pt
\hfill ${\displaystyle
\langle \varphi, {\cal{L}}_\delta \varphi \rangle_{L^2} \, \geq \,
\int \frac{3}{13} \, \delta^2 ( \Delta \varphi )^2 \, + \,
\frac{\delta}{16} \, R | \nabla \varphi |^2 .
}$ \enddemo

{\it Remark}.  Lemma 4.2 is a generalization of
[G-2, Th.~A], which considered the case where $\delta = \frac{2}{3}$.
This corresponds to an eigenvalue estimate for the Paneitz operator.  It
is remarkable that, despite the coefficient $\delta$ in the leading term
of (4.2), one can still show that ${\cal{L}}_\delta$ is invertible
(modulo constants) for {\it all} $\delta > 0$.  
\vglue4pt

Define the differential operator
$$
G[w]=e^{4w} \left( \sigma_2(A)-\frac{\delta}{4}\Delta R+2\gamma_1|\eta|^2\right).
$$
If $G[w]=0$, then $g=e^{2w}g_0$ satisfies $(\ast)_\delta$.  From conformal invariance we see that $\int G[w]dv_0=0$.  Thus,
$G:W^{2,2}_0\to L^2_0$, where the subscript $0$ denotes functions with mean value zero.  If we linearize $G$ at a solution of
$(\ast)_\delta$, it follows from Lemma 4.2 that the linearization is invertible.

Now suppose that $\delta_1 \in S$, and that $g_1 = e^{2 w_1} g_0$ is a
smooth solution of $(*)_{\delta_1}$ with positive scalar curvature.
It
follows from [ADN, Th.~13.1] that there is a unique (up to scaling)
smooth solution of $(*)_\delta$ for all $\delta$ sufficiently close to
$\delta_1$. Moreover, since the scalar curvature of $g_1$ is positive,
by taking solutions in a small enough $C^{2 , \alpha}$-neighborhood of
$g_1$ we may conclude that the solutions of $(*)_\delta$ will also have
positive scalar curvature, for $\delta$ close enough to $\delta_1$.  It
follows that $S$ is open, and the proof of Proposition 4.1 is complete.
\hfill\qed

\proclaim{Proposition}
$S$ is closed{\rm .}
\endproclaim 

{\it Proof}.
The proof of Proposition 4.3 consists of two parts.  First, an {\it a priori}
estimate for solutions of $(*)_\delta$ with positive scalar curvature.
A consequence of this estimate will be the following:  if $\{ \delta_k
\}$ is a sequence in $S$, and $\delta_k \rightarrow \bar{\delta}$, then
$(*)_{\bar{\delta}}$ admits a weak $W^{2,2} (M^4)$-solution.  The second
part of the proof is a local estimate which, when combined with the
regularity theory for extremals of the functional determinant developed
in [CGY], will allow us to conclude that this weak solution of
$(*)_{\bar{\delta}}$ is actually smooth with positive scalar curvature.
It then follows that $S$ is closed.

Now, let $\{ \delta_k \}$ be a sequence in $S$, and suppose $\delta_k
\rightarrow \bar{\delta}$. For each $k$, let $g_k = e^{2 w_k} g_0$ be the
corresponding solution of $(*)_{\delta_k}$, normalized so that $ \int w_k \, = \, 0$.
\noindent
Since $\delta_k \in S$, the scalar curvature $R_k$ of $g_k$ is positive.
Therefore, by (3.2) we have the estimate 

\centerline{${\displaystyle
\int \delta_k ( \Delta_0 w_k )^2 \, + \, | \nabla_0 w_k|^4 \, \leq
\, C_0 .
}$}
\vglue4pt\noindent 
From this we conclude that a subsequence of $\{ w_k \}$ (also denoted
$\{ w_k \}$) converges (i) weakly in $W^{2,2} (M^4)$, (ii) strongly in
$W^{1,2} (M^4)$, (iii) almost everywhere to $w \in W^{2,2} (M^4)$.
Moreover, it is clear that $g = e^{2w} g_0$ satisfies
$(*)_{\bar{\delta}}$ weakly, in the sense that (3.5) holds with $\delta
= \bar{\delta}$ for every $\varphi \in W^{2,2}(M^4)$.

To see that  $w$ is smooth, we need a growth estimate on the integral of
$(\Delta_0 w)^2$ on a small ball of radius $r$.  Therefore, fix $P \in
M^4$ and let $\rho > 0$ be small enough so that the geodesic ball $B (
\rho)$ of radius $\rho$ (measures in the $g_0$ metric) centered at $P$
admits normal coordinates $\{ x^i\}$.  Then in $B ( \rho )$ we also
have 
the Euclidean metric and associated Laplacian, gradient, and volume
form:  
\begin{eqnarray*}
\noalign{\vskip-12pt}
ds^2 & = & \sum^{4}_{i = 1} \, dx^i \, \otimes \, dx^i , \\
\bar{\Delta} & = & \sum^{4}_{i = 1} \, \left(
\frac{\partial}{\partial x^i} \right)^2 , \\
\bar{\nabla} & = & \frac{\partial}{\partial x^i} , \\
dx & = & dx^1 \wedge \ldots \wedge \, d x^4 .
\end{eqnarray*}
\noindent
For $r > 0$ sufficiently small, say $r < r_0$, let $\bar{B} (r)$ denote
the Euclidean ball of radius $r$:
$$
\bar{B} (r) \, = \,
\left\{ Q \in M^4 \, | \, \sum^{4}_{i = 1} \, ( x^i (Q))^2 < r \right\} .
$$ 
\noindent
Now fix $r \in ( 0 , r_0 )$ and define $h$ to be the biharmonic
extension of $w$ on $\bar{B} (r):$
$$
\left\{
\begin{array}{llll}
\bar{\Delta}^2 h \, = \, 0 & {\rm in} & \bar{B}(r) , \\ 
\frac{\partial h}{\partial n} \, = \, \frac{\partial w}{\partial n} & {\rm
on} & \partial \bar{B} (r) , \\ 
h \, = \, w & {\rm on} & \partial \bar{B} (r) ,
\end{array}
\right.
$$ 
\noindent
where $\frac{\partial}{\partial n}$ denotes the outward normal
derivative on $\partial \bar{B} (r)$ in the Euclidean metric.  Define
$\varphi \in W^{2,2} (M^4)$ by
$$
\varphi \, = \, \left\{
\begin{array}{lll}
w - h & {\rm in} & \bar{B} (r) , \\ 
0 & {\rm outside \ of } & \bar{B} (r) .
\end{array}
\right.
$$
\noindent
By (3.5),
\begin{eqnarray*}
&&\int_{\bar{B} (r)} \,
\frac{3}{2} \, 
\bar{\delta} \,
\Delta_0 w \,
\Delta_0 ( w - h ) \\ 
&&\qquad + \,
\frac{1}{2} \, ( 3 \bar{\delta}-2 ) 
\Bigg[
\Delta_0 ( w - h ) \, | \nabla_0 w|^2 \, + \,
2 \Delta_0 w \langle \nabla_0 ( w - h ) , \nabla_0 w \rangle_0 \\ 
&& \hskip1.5in + \, 2 | \nabla_0 w|^2 \, \langle \nabla_0 ( w - h ) ,
\nabla_0 w \rangle_0 \Bigg] \\ 
&&\quad= \, \int_{\bar{B} ( r )} \, - 2 U^\delta_0 \,
( w - h ) \, + \, 2 {\rm Ric}_0 ( \nabla_0 (w - h ) , \nabla_0 w ) \\[5pt]
&&\qquad + \,
\frac{1}{2} \, ( \bar{\delta} - 2 ) R_0 \langle  \nabla_0 ( w - h ) ,
\nabla_0 w \rangle_0 
\end{eqnarray*}
 implies
\begin{eqnarray}
&&\int_{\bar{B}(r)} \, \frac{3}{2} \bar{\delta} \,
(\Delta_0 w )^2 \, +  \, \frac{3}{2} \,
( 3 \bar{\delta} - 2 ) \,
\Delta_0 w | \nabla_0 w|^2 \, + \, ( 3 \bar{\delta} - 2 ) | \nabla_0
w |^4  \\[5pt]
&&\quad = \, \int_{\bar{B}(r)} \, \frac{3}{2} \,
\bar{\delta} \, \Delta_0 w \Delta_0 h \, + \,
\frac{1}{2} \, ( 3 \bar{\delta} - 2 ) \, \Delta_0 h | \nabla_0 w|^2
\nonumber\\[5pt]
&&\qquad + \, ( 3 \bar{\delta}-2 ) \, \Delta_0 w \langle \nabla_0 h ,
\nabla_0 w \rangle_0 \nonumber \\[5pt]
&&\qquad + \, ( 3 \bar{\delta} - 2 ) | \nabla_0 w |^2 \langle \nabla_0 h ,
\nabla_0 w \rangle_0 \, - 2 U^\delta_0 ( w - h ) \nonumber \\[5pt]
&&\qquad +\,  2 \, {\rm Ric}_0 ( \nabla_0 ( w - h ) , \nabla_0 w ) \, + \,
\frac{1}{2} \, ( \bar{\delta} - 2 ) R_0 \langle \nabla_0 ( w - h ) ,
\nabla_0 w \rangle_0 .\nonumber
\end{eqnarray} 
Let us denote the expression on the left-hand side (respectively, 
right-hand side) of (4.5) by LHS (resp., RHS).

\proclaim{Lemma}
{\rm (i)}  There is a constant $C_1 = C_1 ( \bar{\delta} , g_0)$ such
that
\begin{equation}
{\rm LHS} \, \geq \, C_1 \left[
\int_{\bar{B} (r)} \,
(\Delta_0 w)^2 \, + \,
| \nabla_0 w |^4 \right] \, - \,
C_1 r^2 .
\end{equation}
\noindent
{\rm (ii)} There is a constant $C_2 = C_2 ( \bar{\delta}, g_0 )$ such
that 
\begin{equation}
{\rm RHS}\, \leq \, C_2 \left[
\int_{\bar{B}(r)} ( \Delta_0 h )^2 \, + \,
| \nabla_0 h |^4 \right] \, + \,
C_2 r^2 .
\end{equation}
\noindent
{\rm (iii)}   There is a constant $C_3 = C_3 ( \bar{\delta}, g_0)$
such that
\begin{equation}
\int_{\bar{B}(r)} \,
\left[
( \bar{\Delta} w )^2 \, + \,
| \bar{\nabla} w |^4 \right]
dx \, \leq \, C_3 \, \int_{\bar{B}(r)} \, \left[
( \bar{\Delta}h)^2 \, + \,
| \bar{\nabla} h|^4 \right] dx \, + \, C_3 r^2. \hskip.2in
\end{equation}
\endproclaim

\demo{Proof} (i) This inequality essentially
follows from the arguments in the proof of (3.2) in Theorem 3.1.  
We begin with a
claim:

\vglue4pt {\it Claim}. $\Delta_0 w \, + \,
| \nabla_0 w|^2 \, \leq \,
\frac{1}{6} \, R_0$ almost everywhere on $M^4$.  

\medskip
To prove the claim, recall that for each $k$,
\begin{equation}
\Delta_0 w_k \, + \,
| \nabla_0 w_k |^2 \, + \,
\frac{1}{6} \, R_k e^{2w_k} \, = \,
\frac{1}{6} \, R_0.
\end{equation}
\noindent
Since $R_k > 0$, if we multiply both sides of (4.9) by any smooth $\psi \geq
0$ and then integrate over $M^4$, we obtain
$$
\int \psi \Delta_0 w_k \, + \,
\psi \, | \nabla_0 w_k | \, \leq \, \int \frac{1}{6} \,
R_0 \psi .
$$
\noindent
Since $w_k  \, \rightharpoonup \, w$ weakly in $W^{2,2} ( M^4)$, this implies
that 
$$
\int \psi \Delta_0 w \, + \,
\psi \, | \Delta_0 w |^2 \, \leq \, \int \frac{1}{6} \, R_0 \psi
$$
\noindent
and the claim follows.
\vglue4pt

Now, to verify (4.6), we consider two different cases.  First, suppose
$\bar{\delta} \in \left[
\frac{2}{3} , 1 \right]$.  Then $3 \bar{\delta} - 2 \geq 0$, so as
before we have
\begin{eqnarray}
\qquad\qquad {\rm LHS} &\hskip-6pt =\hskip-6pt & 
\int_{\bar{B}(r)} \, \frac{3}{2} \,
\bar{\delta} ( \Delta_0 w )^2 \, + \,
\frac{3}{2} \,
( 3 \bar{\delta} - 2 ) \,
\Delta_0 w | \Delta_0 w|^2 \, + \,
(3 \bar{\delta} - 2 ) | \nabla_0 w|^4 \\[5pt]
& \hskip-6pt\geq \hskip-6pt& \int_{\bar{B} (r)} \, \frac{3}{2} \,
\bar{\delta} ( \Delta_0 w)^2 \, - \,
\frac{9}{16} \,
( 3 \bar{\delta} - 2 ) ( \Delta_0 w)^2\nonumber \\[5pt]
&\hskip-6pt\hskip-6pt& -\,
(3 \bar{\delta}-2 ) | \nabla_0 w |^4 \, + \, (3 \bar{\delta} - 2 ) |
\nabla_0 w |^4 \nonumber \\[5pt]
&\hskip-6pt = \hskip-6pt& \int_{\bar{B}(r)} \, \frac{3}{16} \,
( 6 - \bar{\delta} ) ( \Delta_0 w )^2 \nonumber \\[5pt]
& \hskip-6pt\geq \hskip-6pt& \int_{\bar{B}(r)} \, \frac{15}{16} \,
( \Delta_0 w )^2 .\nonumber 
\end{eqnarray} 
By the claim above,
\begin{eqnarray}\hskip.5in
\int_{\bar{B} (r)} \, | \nabla_0 w|^2 \, \Delta_0 w   +  
| \nabla_0 w |^4 & \leq & \int_{\bar{B}(r)}  \frac{1}{6}  
R_0 | \nabla_0 w |^2 \\[5pt]
& \leq & \left( 
\int_{\bar{B}(r)} | \nabla_0 w |^4 \right)^{\frac{1}{2}} 
\left( \int_{\bar{B} (r)} 
\left( \frac{R_0}{6} \right)^2 
\right)^{\frac{1}{2}} .\nonumber 
\end{eqnarray}
\noindent
Notice that $w \in W^{2,2} (M^4) \subset W^{1,4} (M^4)$ implies that
each integral in (4.11) is well-defined.  Moreover, 
\begin{eqnarray*}
\left(
\int_{\bar{B}} \, | \nabla_0 w |^4 \right)^{\frac{1}{2}} \,
\left(
\int_{\bar{B}} \, 
\left( \frac{R_0}{6} \right)^2 \right)^{\frac{1}{2}} & \leq & C \,
\left(
\int_{\bar{B} (r)} \right)^{\frac{1}{2}} \\
\\
& \leq & C r^2 ,
\end{eqnarray*}
\noindent
and we conclude
\begin{equation}
\int_{\bar{B} (r)} \,
| \nabla_0 w|^2 \Delta_0 w \, + \,
| \nabla_0 w |^4 \, \leq \, C r^2 .
\end{equation}
\noindent
From (4.12) we also have
\begin{eqnarray*}
\int_{\bar{B}(r)} \, | \nabla_0 w|^4 & \leq &
\int_{\bar{B}(r)} \, | \nabla_0 w|^2 \, ( - \Delta_0 w) \, + \, Cr^2
\\
& \leq &
\frac{1}{2} \,
\int_{\bar{B}(r)} \, | \nabla_0 w|^4 \, + \,
\frac{1}{2} \,
\int_{\bar{B}(r)} \,
( \Delta_0 w)^2 \, + \, Cr^2
\end{eqnarray*}
implies
\begin{equation}
\int_{\bar{B}(r)} \,
| \nabla_0 w|^4 \, \leq \,
\int_{\bar{B}(r)} \,
( \Delta_0 w)^2 \, + \,
Cr^2 .
\end{equation} 
\noindent
Combining (4.10) and (4.13) we see that (4.6) holds when 
$\bar{\delta} \in 
\left[
\frac{2}{3} , 1 
\right]$.  If 
$\bar{\delta} \in \left(
0 , \frac{2}{3} \right)$, then 
$3 \bar{\delta} - 2 < 0$ so from (4.12)
we get 
$$
\int_{\bar{B}(r)} \,
\frac{3}{2} \,
\left( 
3 \bar{\delta} - 2 
\right) \,
\Delta_0 w | \nabla_0 w |^2 \, \geq \,
\int_{\bar{B}(r)} \, - 
\frac{3}{2} \,
( 3 \bar{\delta} - 2 ) \,
| \nabla_0 w |^4 \, - C r^2
$$
\noindent
so that
$$
{\rm LHS} \, \geq \,
\int_{\bar{B}(r)} \,
\frac{3}{2} \, 
\bar{\delta} ( \Delta_0 w )^2 \, -
\frac{1}{2} \,
( 3 \bar{\delta} - 2 ) \,
| \nabla_0 w |^4 \, - C r^2 .
$$
\noindent
This completes the proof of (i). \vglue4pt

To prove (ii), we need to consider each term on the RHS separately.  This
will be considerably simplified if we begin with the following crude
estimate:
\begin{eqnarray}
{\rm RHS} & \lesssim & \int_{\bar{B}(r)} \, | \Delta_0 w | \, | \Delta_0
h | \, + \,
\int_{\bar{B}(r)} \,
| \nabla_0 w |^2 | \Delta_0 h| \\[5pt]
&& +\,  \int_{\bar{B}(r)} \,
| \Delta_0 w | \, | \nabla_0 w | \, | \nabla_0 h | \, + \,
\int_{\bar{B}(r)} \, | \nabla_0 w|^3 \, | \nabla_0 h | \nonumber \\[5pt]
&& + \, \int_{\bar{B}} \,
| w - h | \, + \, \int_{\bar{B}} \, | \nabla_0 w | \, | \nabla_0 ( w - h
) | \nonumber 
\end{eqnarray} 
where $\lesssim$ means that the inequality holds up to a multiplicative
constant which depends on $\bar{\delta}$ and $g_0$.  Using the
arithmetic-geometric mean inequality we can estimate the first four
integrals in (4.14) as follows: 
\begin{eqnarray}
\int_{\bar{B}(r)} \, 
| \Delta_0 w | \, | \Delta_0 h | & \leq &
\frac{1}{2} \, \varepsilon \,
\int_{\bar{B}(r)} \, 
( \Delta_0 w )^2 \, + \, \frac{1}{2} \,
\varepsilon^{-1} \, \int ( \Delta_0 h )^2 ,  \\[5pt]
\int_{\bar{B}(r)} \,
| \nabla_0 w |^2 \, | \Delta_0 h | & \leq &
\frac{1}{2} \, \varepsilon \,
\int_{\bar{B}(r)} \, | \nabla_0 w |^4 \, + \, 
\frac{1}{2} \, \varepsilon^{-1} \,
\int_{\bar{B}(r)} \, ( \Delta_0 h)^2 , \nonumber \\[5pt]
\int_{\bar{B}(r)} \,
| \Delta_0 w | \, | \nabla_0 w | \, 
| \nabla_0 h | & \leq &
\frac{1}{2} \,
\varepsilon 
\left[
\int_{\bar{B}} \, ( \Delta_0 w )^2 \, + \, | \nabla_0 w|^4 
\right]
\nonumber\\
&& + \, \frac{1}{8} \, \varepsilon^{-3} \, \int_{\bar{B}(r)} \, |
\nabla_0 h|^4 ,  \nonumber \\[5pt]
\int_{\bar{B}(r)} \,
| \nabla_0 w |^3 \, | \nabla_0 h| & \leq & 
\varepsilon \, \int_{\bar{B} (r)} \,
| \nabla_0 w |^4 \, + \,
\frac{1}{8} \,
\varepsilon^{-3} \,
\int_{\bar{B}(r)} \, | \nabla_0 h |^4 .\nonumber
\end{eqnarray} 
\pagebreak

\noindent
To estimate the last two terms in (4.14), observe that on $M^4$
$$
\parallel h \parallel_{2,2} \, \leq \, C ( \parallel w \parallel_{2,2} )
$$
\noindent
where $C$ is independent of $r$ (see [CGY-1, p.\ 237]).  Therefore,
\begin{equation}
\int_{\bar{B}(r)} \, | w - h | 
\, \leq \, 
\left( \int_{\bar{B}(r)} \, (w - h )^2 \right)^{\frac{1}{2}} \,
\left( \int_{\bar{B}(r)} \right)^{\frac{1}{2}} \, \leq \, Cr^2, 
\end{equation}
\begin{eqnarray}
&&\int_{\bar{B}(r)} 
| \nabla_0 w |  
| \nabla_0 (w - h ) | \\
&&\qquad\quad \leq  \,
\left(
\int_{\bar{B}(r)}   | \nabla_0 w |^4 \right)^{\frac{1}{4}}  
\left( \int_{\bar{B}(r)} \, | \nabla_0 ( w - h ) |^4
\right)^{\frac{1}{4}} 
\left( \int_{\bar{B}(r)} \right)^{\frac{1}{2}}\nonumber  \\
&&\qquad\quad \leq \, C r^2 .\nonumber
\end{eqnarray}
\noindent
By substituting (4.15)--(4.17) into (4.14) and choosing $\varepsilon > 0$
sufficiently small we get (4.7).  

\medskip
Finally, to prove (iii), we combine (4.6) and (4.7) to get
\begin{equation}
\int_{\bar{B}(r)} \,
\left[
( \Delta_0 w )^2 \, + \,
| \nabla_0 w|^4 \right] \, \lesssim \,
\int_{\bar{B}(r)} \,
\left[
( \Delta_0 h)^2 \, + \,
| \nabla_0 h |^4 \right] \, + \, r^2 .\hskip.35in
\end{equation}
\noindent
Then (4.8) follows from (4.18) by appealing to [CGY-1, (3.1)], which
compares the Euclidean volume form, Laplacian, and gradients appearing
in (4.8) with their Riemannian counterparts in (4.18).  The details
will be omitted. \hfill\qed\vglue4pt

Inequality (4.8) is precisely the conclusion of [CGY-1, Lemma 3.4].  We
can therefore apply the subsequent arguments of [CGY-1] to conclude that
$w \in C^\infty (M^4)$.  Moreover, it follows from the claim in Lemma
4.4 that the scalar curvature $R$ of $g = e^{2w} g_0$ is nonnegative.
Since $g$ satisfies $(*)_{\bar{\delta}}$, the scalar curvature satisfies
\begin{eqnarray*}
\bar{\delta} \Delta R & = &
8 \gamma_1 | \eta|^2 \, - 2|E|^2 \, + \, \frac{1}{6} \, R^2 \\
& \leq & \frac{1}{6} \, R^2 .
\end{eqnarray*}
\noindent
Thus, by the minimum principle, $R > 0$ on $M^4$.  This completes the
proof of Proposition 4.3. \hfill\qed
\vglue4pt

Combining Propositions 4.1 and 4.3 we conclude:
\proclaim{{C}orollary}
If $\int \sigma_2 (A_0 ) dv_0 > 0$ and $Y ( g_0) > 0$ then for each
$\delta > 0${\rm ,} $(*)_\delta$ admits a smooth solution with positive scalar
curvature{\rm .}
\endproclaim
\pagebreak

\section{{\it A priori} $W^{2,3}$ estimates}

%\end{document}

\hskip .2in Our goal in this section is to establish the following {\it a priori}
estimate:

\proclaim{Theorem}
Let $g = e^{2w} g_0$ be a solution of $(*)_\delta$
with positive scalar curvature{\rm ,} normalized so that $\int
w dv_0 = 0${\rm ,} and assume
\begin{equation}
\int \sigma_2 (A_0 ) d v_0 \, = \, \int \sigma_2 (A) dv \, > \, 0 .
\end{equation}
Then there are constants $C = C ( g_0 )$ and $ 0 < \delta_0 < 1$ such that
\begin{equation}
\int | \nabla^2_0 w|^3 d v_0 \, + \, \int | \nabla_0 w |^{12} d v_0 \, \leq \, C
\end{equation}
for $0 < \delta < \delta_0${\rm .} In particular{\rm ,} for any $\alpha
\in (0, \frac{2}{3})$ there is a constant $C_{\alpha} =
C(\alpha, g_0)$ such that
\begin{equation}
|| w ||_{C^{\alpha}} \leq C_{\alpha}.
\end{equation}
\endproclaim 

The proof of (5.2) is quite involved, and will be obtained through a
series of lemmas and propositions.  The basic estimate we will need is:

\proclaim{Proposition}
Under the same hypotheses of Theorem {\rm 5.1,} there is a constant $C = C (
g_0 )$ such that 
\begin{equation}
\int \left( \frac{R}{6} \right)^3 dv \, \leq \,
( 1 + C \delta ) \, \int | \nabla w|^6 dv \, + \,
C \, \int R^2 dv \, + \, C 
\end{equation}
\noindent
for $\delta$ sufficiently small.
\endproclaim

\demo{Proof} We  begin  with some
notational conventions.  First, all integrals (unless otherwise
specified) are with respect to the volume form $dv$ of $g$. With this
understood, we will suppress it from now on.  Second, our calculations
will sometimes be facilitated by introducing local coordinates.  These
coordinates are assumed to be normal at some point, and certain
identities involving covariant derivatives are understood to hold only
at that point.  For this reason we will not bother to distinguish
between raised and lowered indices;  all indices will be subscripts.
For example, the gradient of a function will be denoted by $\nabla_j
\varphi$, and by our conventions ${\rm Ric} ( \nabla \varphi , \nabla \varphi
) =  R_{ij} \nabla_i \varphi \nabla_j \varphi$ in local coordinates.
Finally, $C$ denotes a constant which depends at most on $g_0$.

Our arguments are somewhat imitative of the {\it a priori} $C^2$-estimates for
Monge-Amp\`ere equations as described in [CNS-1], [CNS-2], [Ev], and [K].  However, the
estimates in these references are pointwise in nature and involve the
maximum principle.  Since our regularized equation is fourth order, such
techniques cannot work for us.  Instead, we rely on integral estimates
which (as we will see) present their own difficulties.  Since our
calculations are quite involved, it may be helpful to begin with an
overview of the argument.

We begin with a simple identity.  Let $f \in C^\infty (M^4)$.  Then by
the divergence theorem,
\begin{eqnarray*}
\noalign{\vskip-8pt}
0 & = & \int \nabla_i ( S_{ij} \, \nabla_j f ) \\
& = & \int \nabla_i S_{ij} \, \nabla_j f \, + \, S_{ij} \, \nabla_i
\nabla_j f .
\end{eqnarray*}
\noindent
Also, the contracted second Bianchi identity implies that $S$ is
divergence-free: $\nabla_i S_{ij} = \nabla_i ( - R_{ij} + \frac{1}{2} \,
R g_{ij} ) = - \nabla_i R_{ij} + \frac{1}{2} \, \nabla_j R = 0$.
Therefore, 
\begin{equation}
0 \, = \,
\int S_{ij} \, \nabla_i \, \nabla_j f , \hbox{ for  any }
  f \in C^\infty ( M^4) .
\end{equation} 
\noindent
We will apply (5.5) to two different choices of $f$, resulting in two
different inequalities.  First, we let $f = R$.  Then differentiating
$(*)_\delta$ twice and using (5.5) we obtain the
inequality (see (5.25))
\begin{equation}
0 \, = \,
\int S_{ij} \, \nabla_i \, \nabla_j R \, \geq \,
\int 6 \, {\rm tr}\,E^3  + \frac{1}{12} \, R^3 \nonumber
 + ( \rm lower \ order \ terms ),\hskip.4in
\end{equation}
\noindent
where ${\rm tr}\, E^3 = E_{ij} E_{ik} E_{jk}$.  Next, we let $f = 12 | \nabla
w|^2$, resulting in the inequality
\begin{eqnarray}
0 & = & \int S_{ij} \nabla_i \nabla_j ( 12 | \nabla w |^2 ) 
\\
& \geq & \int -  \, 6 {\rm tr}\, E^3 \, + \,
\frac{1}{12} \, R^3 \, - \,
6 \, R | \nabla w|^4 \, + \, \hbox{(l.o.t.)}.\nonumber 
\end{eqnarray}
\noindent
Adding (5.6) to (5.7), we see that the term ${\rm tr}\, E^3$ cancels to give
$$
\int \left(
\frac{R}{6} \right)^3 \, \leq \,
\int \left( \frac{R}{6} \right) \, | \nabla w|^4 \, + \,
\hbox{( l.o.t.)},
$$
\noindent
and (5.4) can be shown to follow from this inequality.  A further argument
is needed to derive (5.2), but this brief overview should  
provide the reader with a rough guide to the estimates which follow.
To begin, let
\begin{equation}
I \, = \,
\int S_{ij} \, \nabla_i \, \nabla_j R .
\end{equation}
\noindent
By (5.5), $I = 0$.  Also, we have:
\pagebreak

\proclaim{Proposition}
\begin{equation}
I \, \geq \, 
\int 3 \left(
| \nabla E|^2 \, - \, \frac{1}{12} \,
| \nabla R|^2 \right) \, + \,
6 \, {\rm tr}\, E^3 \, +  \, \frac{1}{12} \, R^3 \, - CR^2 \, - \, C .\qquad
\end{equation}
\endproclaim

{\it Proof}. Inequality (5.9) is a consequence 
of the following fundamental identity:
 
\proclaim{Lemma}
Let $( M^4 , g)$ be any Riemannian $4$\/{\rm -}\/manifold{\rm .} Then
\begin{eqnarray}
\qquad S_{ij} \nabla_i \nabla_j R & = & 3 \Delta \sigma_2 (A) \, + \,
3 \left( | \nabla E|^2 \, - \, \frac{1}{12} \,
| \nabla R|^2 \right)  \\
& & +\, 6 \, {\rm tr}\, E^3 \, + \, R|E|^2 \, - \,
6 \, W_{ijk\ell} \, E_{ik} E_{j \ell} \, - \, 6 \, E_{ij} \, B_{ij}\nonumber
\end{eqnarray}
\noindent
where $B_{ij}$ denotes the Bach tensor{\rm .}
\endproclaim

\demo{Proof} By (1.18),
\begin{eqnarray*}
\Delta E_{ij} & = & \frac{1}{3} \, \nabla_i \nabla_j R \, - \,
\frac{1}{12} \, \Delta \, R g_{ij} \, + \,
2 E_{ik} E_{jk} \, - \,
\frac{1}{2} \, | E |^2 g_{ij} \\
\\
& & + \, \frac{1}{3} \, R E_{ij} \, - \, 2W_{ikj\ell} \, E_{k \ell} \, - 2
B_{ij} . 
\end{eqnarray*}
\noindent
Thus, 
\begin{eqnarray*}
\frac{1}{2} \, \Delta \, |E|^2 & = & | \nabla E|^2 \, + \,
E_{ij} \, \Delta E_{ij} \\
\\
& = & | \nabla E|^2 \, + \, \frac{1}{3} \, E_{ij} \, \nabla_i \nabla_j R
\, + \,
2 \, {\rm tr}\, E^3 \, + \,
\frac{1}{3} R | E|^2 \\
\\
& & - \, 2 W_{ikj\ell} \, E_{ij} \, E_{k \ell} \, - \, 2B_{ij} \, E_{ij} ,
\end{eqnarray*}
\begin{eqnarray}
\qquad\quad\Delta \sigma_2 (A) & = & \Delta \left( - \frac{1}{2} \, | E|^2 \, + \,
\frac{1}{24} \, R^2 \right)  \\
& = & - | \nabla E|^2 \, + \, \frac{1}{12} \, | \nabla R|^2 \, - \,
\frac{1}{3} \,
E_{ij} \, \nabla_i \nabla_j R \, + \, \frac{1}{12} \, R \Delta R
\nonumber \\
& & - \, 2 \, {\rm tr}\, E^3 \, - \, \frac{1}{3} \, R | E|^2 \, + \, 2 \, W_{ikj\ell} \,
E_{ij} \, E_{k \ell} \, + \, 2B_{ij} \, E_{ik} .\nonumber
\end{eqnarray}
Note that $E_{ij} = \, - S_{ij} + \, \frac{1}{4} \, R g_{ij}$, so that
(5.11) can be rewritten
\begin{eqnarray*}
\Delta \sigma_2 (A) & = & - \left(
| \nabla E|^2 \, - \,
\frac{1}{12} \, | \nabla R |^2 
\right) \, + \, \frac{1}{3} \, S_{ij} \,
\nabla_i \nabla_j R \\
&& - \, 2 \, {\rm tr}\, E^3 \, - \,
\frac{1}{3} \,
R | E|^2 \, + \, 2 \, W_{ikj\ell} \, E_{ij} \, E_{k \ell} \, + \, 2 B_{ij} \,
E_{ij} , 
\end{eqnarray*}
\noindent
and (5.10) follows.
\hfill\qed

Continuing the proof of Proposition 5.3, we integrate
(5.10) over $M^4$: 
\begin{eqnarray}
I & = &
\int S_{ij} \nabla_i \nabla_j R \\
& = & \int 3 \left(
| \nabla E|^2\, - \frac{1}{12} \,
| \nabla R|^2 \right) \, + \,
6 \, {\rm tr}\, E^3 \, + \, R | E |^2 \nonumber \\
&& - \, 6 \, W_{ikj \ell} \, E_{ij} \, E_{k \ell} \, - \,
6 \, B_{ij} \, E_{ij} . \nonumber
\end{eqnarray} 
\noindent
Now, 
$$
| W_{ikj \ell} E_{ij} \, E_{k \ell} | \, \lesssim \, | W | \, | E|^2 .
$$
\noindent
By the transformation law for   Weyl curvature, $|W| = e^{-2w} | W_0
|_0$.  From (3.1), we conclude that 
$|W| \, \leq \, C$.  Thus
$$
| W_{i k j \ell} \, E_{ij} \, E_{k \ell} | \, \leq \,
C | E |^2 .
$$
\noindent
Similarly,
$$
| B_{ij} \, E_{ij} | \, \leq \, | E | \, | B| ,
$$
\noindent
and $|B| = e^{-4w} | B_0 |_0 \, \leq \, C$.  Thus, the last two terms in
(5.12) can be estimated by 
\begin{eqnarray}
&&\int - 6 \, W_{ikj \ell} \, E_{ij} \, E_{k \ell} \, - \, 6 \, B_{ij} \, E_{ij}
 \\
&&\qquad \qquad\gtrsim  \int - \, | E |^2 \, - \, |E|  \nonumber \\
&&\qquad \qquad \gtrsim  - \, \int ( |E|^2 \, + \, C ) .    \nonumber
\end{eqnarray}
\noindent
Note that the estimate (3.1) implies that the volume of $g$ has 
uniform upper and lower bounds: $\int dv \sim  1 $.  It is therefore irrelevant 
whether we place the constant in (5.13) inside the integrand or outside.

By (5.1), 
$$
0 < \int \sigma_2 (A) \, = \, \int - \,
\frac{1}{2} \, | E |^2 \, + \, \frac{1}{24} \, R^2 ,
$$
\noindent
so that $\int - | E |^2 \, \geq \,
- \, \frac{1}{12} \, \int R^2$.  Substituting this inequality into
(5.13) and combining with (5.12) we obtain
\begin{eqnarray}
I & \geq & \int \, 3 \left(
| \nabla E |^2 \, - \frac{1}{12} \, | \nabla R |^2 \right) \, + \,
6 \, {\rm tr}\, E^3 \, + \, R |E |^2 \\
& & - \, CR^2 \, - \, C. \nonumber
\end{eqnarray} 
\noindent
The conclusion of Proposition 5.3 will then follow from
(5.14) along
with part (i) of the following lemma:

\proclaim{Lemma}
\, {\rm(i)} \ $\int R | E|^2 \, \geq \, \int \frac{1}{12} \, R^3 \, - \, CR^2 \, -\, C${\rm , } 
 \vglue4pt
{\rm (ii)} $\int R | E |^2 \, \leq \, \int \frac{1}{2} \, \delta | \nabla R|^2 \, + \, \frac{1}{12} \, R^3${\rm .}
\endproclaim

\demo{Proof} From $(*)_\delta$,
\begin{eqnarray*}
  \int \delta R \, \Delta R &  = & \int R \left[
8 \gamma_1 | \eta|^2 \, - \,
2 | E |^2 \, + \, 
\frac{1}{6} \, R^2 
\right] \\
  \Longrightarrow \,  \int R | E |^2  & = &  \int \frac{1}{2} \, \delta | \nabla R |^2 \, + \, 4 \gamma_1 R | \eta |^2 \,
+
\, \frac{1}{12} \, R^3 .
\end{eqnarray*}
\noindent
Since $\gamma_1 < 0$ by (5.1), (ii) is immediate.  To see (i), observe
that $| \eta |^2 = e^{- 4 w}| \eta |^2_0 \, \leq \, C$, so that
$\int 4 \gamma_1 R | \eta |^2 \, \gtrsim \, \int - R^2 \, - 1$. \enddemo

In the next two lemmas we undertake some precise estimates of the terms
in (5.9).
\proclaim{Lemma}
For any $p \geq 0${\rm ,}
\begin{eqnarray}
&&\hskip-24pt\int 3 \, R^p \left(
| \nabla E|^2 \, - \,
\frac{1}{12} \, | \nabla R \,|^2 \right)\\
&&\qquad \geq 
\int \frac{3}{2} \, \delta R^{p-1} \, ( \Delta R )^2 + \, \frac{3}{2} \, 
\delta p \, R^{p-2} \, | \nabla R|^2 \, \Delta R\nonumber \\
&&\qquad\quad
+ \, 12 \gamma_1 \, R^{ p-1} \langle \nabla R, \nabla | \eta |^2 \rangle   - \, 12 \gamma_1 \, R^{p-2} | \eta |^2 \, | \nabla R|^2
. \nonumber
\end{eqnarray}
\endproclaim

\demo{Proof}
Differentiating $(*)_\delta$, we obtain 
$$
0 \, = \, \delta \nabla ( \Delta R) \, - \,
8  \gamma_1
  \nabla | \eta |^2 \, + \,
4 |E| \nabla |E| \, - \, \frac{1}{3} \, R \nabla R . $$
\noindent
Now take the inner product of both sides with $R^{p-1} \nabla R$ and
integrate:
\begin{eqnarray}
0 & = & \int \delta R^{p-1} \langle \nabla R, \nabla ( \Delta R )
\rangle \, - \, 8 \gamma_1  R^{p-1} \langle \nabla R, \nabla | \eta |^2
\rangle  \\
  && + \, 4R^{p-1} | E | \langle \nabla R, \nabla |E| \rangle \, - \,
\frac{1}{3} \, R^p | \nabla R |^2 .\nonumber
\end{eqnarray}
Using the arithmetic-geometric mean inequality (henceforth referred to
as the AGM inequality),  we obtain
$$ 
\displaystyle{\int} 4R^{p-1} | E | \langle \nabla R, \nabla | E |
\rangle \, \leq \,
\displaystyle{\int} 2 R^p | \nabla | E | |^2 \, + \,
2 R^{p-2} | E |^2 | \nabla R|^2 .
$$
\noindent
By Kato's inequality $| \nabla | E | |^2 \leq | \nabla E|^2$,  
\begin{equation}
\int 4R^{p-1} | E | \langle \nabla R, \nabla | E | \rangle \, \leq
\, \int 2 R^p | \nabla E|^2 \, + \,
2R^{p-2} | E |^2 | \nabla R|^2 .\qquad
\end{equation}
\noindent
We now use $(*)_\delta$ to substitute $|E|^2$ into the last term of
(5.17):
\begin{eqnarray}
&& \int 2 R^{p-2} | E |^2 \, | \nabla R|^2\\
& &\qquad\quad=\,
\int 2 R^{p-2} \left\{
- \frac{\delta}{2} \, \Delta R \, + \,
4 \gamma_1 | \eta|^2 \, + \, \frac{1}{12} \, R^2 \right\} | \nabla R|^2
 \nonumber\\
&&\qquad\quad=\,\int - \delta R^{p-2} | \nabla R|^2 \, \Delta R \, + \, 8
\gamma_1   R^{p-2} | \eta |^2 | \nabla R |^2 
+ \, \frac{1}{6} \, R^p | \nabla R |^2 . \nonumber
\end{eqnarray}
By (5.17) and (5.18),
\begin{eqnarray}\qquad\quad
\int 4 \, R^{p-1} | E | \langle \nabla R, \nabla | E | \rangle  & \leq&
\, 
\int 2 R^{p} | \nabla E|^2 \, + \,
\frac{1}{6} R^p | \nabla R |^2  \\
&&-\, \delta R^{p-2} | \nabla R|^2 \Delta R \, + \, 8 \gamma_1 R^{p-2} |
\eta |^2 | \nabla R|^2 .\nonumber
\end{eqnarray}
\noindent
Returning to (5.16), we integrate by parts in the first term:
\begin{equation}
\int \delta R^{p-1} \, \langle \nabla R, \nabla ( \Delta R ) \rangle
\, = \,
\int - \delta R^{p-1} ( \Delta R)^2 \, - \,
\delta ( p - 1 ) R^{p-2} | \nabla R |^2 \, \Delta R .
\end{equation} 
\noindent
Substituting (5.19) and (5.20) into (5.16), multiplying by
$\frac{3}{2}$, then rearranging terms, we get (5.15). \hfill\qed

\proclaim{{C}orollary}
\begin{equation}
\int 3 
\left(
| \nabla E|^2 - \, \frac{1}{12} \,
| \nabla R|^2 \right) \, \geq \,
\int \frac{3}{2} \, \delta \, \frac{( \Delta R)^2}{R} \, - \, C .
\end{equation}
\endproclaim

\demo{Proof}  If we take $p = 0 $ in
(5.15),
then 
\begin{eqnarray}\quad
\int 3 \left(
| \nabla E |^2 -  
\frac{1}{12} \, | \nabla R |^2 \right) &\hskip-7pt \geq\hskip-7pt &
\int \frac{3}{2}   \delta 
\frac{(\Delta R )^2}{R}  +  
12 \gamma_1 R^{-1} \langle \nabla R , \nabla | \eta |^2 \rangle)
 \\
&\hskip-7pt\hskip-7pt& - \, 12 \gamma_1 R^{-2} | \eta |^2   | \nabla R |^2 .\nonumber
\end{eqnarray}
\noindent
Since $\gamma_1 < 0$, by the AGM inequality,
\begin{eqnarray}
&&12 \gamma_1 R^{-1} \langle \nabla R , \nabla
| \eta |^2 \rangle \, - \,
12 \gamma_1 R^{-2} | \eta |^2 \, | \nabla R|^2   \\
&&\qquad\quad  =   24 \gamma_1 R^{-1} \, | \eta | \, \langle \nabla R , \nabla |
\eta | \rangle \, - \, 12 \gamma_1 R^{-2} | \eta |^2 | \nabla R |^2
\nonumber \\
&&\qquad\quad  \geq   12 \gamma_1 | \nabla | \eta | \, |^2 .\nonumber
\end{eqnarray} 
Now, 
\begin{eqnarray}
  \int | \nabla | \eta | \, |^2 \, dv & = &
\int e^{2w} | \nabla_0 ( e^{-2w} | \eta |_0 ) |^2 \, dv_0 
\\
& \lesssim & \int e^{-2w} \, dv_0 \, + \, \int e^{-2w} | \nabla_0
w|^2 \, dv_0 \nonumber \\
& \lesssim & 1 \, + \,
\int | \nabla_0 w|^2 \, dv_0 \, \leq \, C ,  \nonumber
\end{eqnarray}
where the last inequality follows from (3.1). Then (5.21) follows
from (5.22)--(5.24).
\enddemo

\proclaim{{C}orollary}
\begin{equation}
I \, \geq \,
\int \frac{3}{2} \, \delta \,
\frac{(\Delta R)^2}{R} \, + \,
6 {\rm tr}\, E^3 \, + \,
\frac{1}{12} \, R^3 \, - \,
CR^2 \ - \, C.  
\end{equation}
\endproclaim

\demo{Proof} This follows from (5.9)
and (5.21).
\enddemo

The next result amounts to a technical lemma, which will be useful in
our subsequent estimates.

\proclaim{Lemma}
If $\delta < \frac{1}{2}${\rm ,} then
\begin{equation}
\int \delta | \nabla R|^2\, \lesssim \, \int \delta R^3 \, + \, R^2 \, + \, 1 .
\end{equation}
\endproclaim

\demo{Proof} Since $E$ is trace-free, we have the sharp
inequality
$$
6 {\rm tr}\, E^3 \, \geq \,  - \, \frac{6}{\sqrt{3}} \, | E |^3
$$  
\noindent
Thus,
\begin{eqnarray*}
6 {\rm tr}\, E^3 \, + \, R |E|^2  & \geq & - \,
\frac{6}{\sqrt{3}} \, | E |^3 \, + \, R |E|^2 \\
\\
& = & | E |^2 \left(
- \, 2 \sqrt{3} | E | \, + \, R \right).
\end{eqnarray*}
\noindent
Using the AGM inequality 
$$
- 2 \sqrt{3} \, | E | \, = \, - \, 2 ( \sqrt{6} \, | E | R^{-1/2} ) \,
\left( \frac{1}{\sqrt{2}} R^{1/2} \right) \, \geq \,
- \, 6 | E |^2 R^{-1} \, - \frac{1}{2} R , $$
\noindent
we get
\begin{eqnarray}
\int 6 {\rm tr}\, E^3 \, + \, R|E|^2 & \geq & \int \frac{|E|^2}{R} \, \left( - \, 6 |E|^2 \, + \, \frac{1}{2} R^2 \right)
 \\
& = & \int \frac{|E|^2}{R} \, ( 12 \sigma_2 (A) ) \nonumber \\
& = & \int \frac{|E|^2}{R} \,
( 3 \delta \Delta R \, - \, 24 \gamma_1 | \eta|^2 ) \nonumber \\
& = & \int 3 \delta \, \frac{\Delta R}{R} \, | E |^2 \, - \,
24 \gamma_1  \, \frac{|E|^2}{R} \, | \eta |^2 \nonumber \\
& \geq & \int 3 \delta \, \frac{\Delta R}{R} \, | E |^2 ,\nonumber
\end{eqnarray}  
the last line following from the fact that $\gamma_1 < 0$.  Using the
AGM inequality once again,  we obtain
$$
\int 3 \delta \, \frac{\Delta R}{R} \, | E |^2 \, \geq \,
\int - \, \frac{1}{2} \, \delta \,
\frac{(\Delta R)^2}{R} \, - \, \frac{9}{2} \, \delta \,
\frac{|E|^4}{R} ,
$$
\noindent
and substitution of $(*)_\delta$ into the last term above gives 
\begin{eqnarray*}
\int 3 \delta \, \frac{\Delta R}{R} \, | E |^2 & \geq & 
\int - \, \frac{1}{2} \, \delta\,
\frac{(\Delta R)^2}{R} \, - \,
\frac{9}{2} \, \delta \,
\frac{|E|^2}{R} \, \left(
- \, \frac{\delta}{2} \, \Delta R \, + \,
4 \gamma_1 | \eta |^2 \, + \, \frac{1}{12} R^2 \right) \\ 
& = & \int - \,
\frac{1}{2} \, \delta \,
\frac{(\Delta R)^2}{R} \, + \,
\frac{9}{4} \, \delta^2 \, 
\frac{\Delta R}{R} \, | E |^2 \, - 18 \delta \gamma_1 \, \frac{|E|^2}{R} \, |
\eta |^2 \\
&& - \, \frac{3}{8} \, \delta R | E |^2 .
\end{eqnarray*}
Combining terms, we have
\begin{eqnarray}
&&\int 3 \delta \left( 1 - \, \frac{3}{4} \, \delta \right) \,
\frac{\Delta R}{R} \,
| E |^2 \, \geq \, \int - \, \frac{1}{2} \, \delta \, 
\frac{(\Delta R)^2}{R} \, - \,
\frac{3}{8} \, \delta R | E |^2
\\
&\Longrightarrow&
\int 3 \delta \, \frac{\Delta R}{R} \, | E |^2 \, \geq \,
\int - \, \frac{\delta}{2} \, 
\left( 1 - \, \frac{3}{4} \, \delta \right)^{-1} \,
\frac{(\Delta R)^2}{R} \, - \,
\frac{3}{8} \,
\delta ( 1 - \,
\frac{3}{4} \, \delta )^{-1} \, R | E |^2 . \nonumber
\end{eqnarray} 
\noindent 
Substitution of (5.28) into (5.27) gives 
\begin{eqnarray}
\int 6 {\rm tr}\, E^3 \, + \, R | E |^2 & \geq & \int - \,
\frac{\delta}{2} \, \left(
1 - \, \frac{3}{4} \, \delta \right)^{-1} \, 
\frac{( \Delta R )^2}{R}  \\
&& -   \, \frac{3}{8} \, \delta 
\left( 1 - \, \frac{3}{4} \, \delta \right)^{-1} 
\, R | E |^2 .\nonumber
\end{eqnarray} 
\noindent
We now substitute (5.29) into (5.25) and obtain 
\begin{eqnarray*}
0   =   I & \geq& \int \delta \, 
\frac{(2 - 3 \delta )}{(4 - 3 \delta )} \,
\frac{(\Delta R)^2}{R}  \\
&& - \,
\frac{3}{2} \,
\delta \, ( 4 - 3 \delta )^{-1} \, R | E |^2- \, C R^2 \, - \, C ,
\end{eqnarray*}
when $\delta < \frac{1}{2}$.  Note that we can roughly estimate the
terms above by 
$$
0 \, \geq \, \int \frac{\delta}{8} \,
\frac{(\Delta R)^2}{R} \, - \,
\delta R | E |^2 \, - \,
CR^2 \, - \, C . $$
\noindent
By Lemma 5.5 (ii) this implies
\begin{equation}
0 \, \geq \, \int \frac{\delta}{8} \,
\frac{(\Delta R)^2}{R} \, - \,
\frac{\delta^2}{2} \,
| \nabla R |^2 \, - \,
\frac{\delta}{12} \,
R^3 \, - \,
CR^2 \, - \, C .
\end{equation}
\noindent
Finally, notice that
\begin{eqnarray*}
\int 2 \delta \, | \nabla R |^2 & = & \int - 2 \delta R \Delta R
\\
\\
& \leq & \int \frac{\delta}{8} \,
\frac{( \Delta R )^2}{R} \, + \, 8 \delta R^3 \\
\Longrightarrow
\int \frac{\delta}{8} \,
\frac{(\Delta R)^2}{R} & \geq & \int 2 \delta | \nabla R|^2 \, - \,
8 \delta R^3 .
\end{eqnarray*}
\noindent
Substituting this into (5.30) we conclude
\vglue6pt
\hfill ${\displaystyle
0 \, \geq \,
\int 2 \delta \left( 1 - \frac{\delta}{4} \right) \, | \nabla R |^2
\, - \,
C \delta R^3 \, - \,
CR^2 \, - \, C . }$ \enddemo

Now define $V = \frac{1}{2} | \nabla w|^2$.  By (5.5), 
\begin{equation}
{\rm II} \equiv \: \int S_{ij} \nabla_i\nabla_j V \, = \, 0 .
\end{equation} 
\proclaim{Lemma}
\begin{eqnarray}
  S_{ij} \, \nabla_i \nabla_j V &\hskip-7pt =\hskip-7pt&S_{ij}\, \nabla_i \nabla_k w
\nabla_j\nabla_k w \, - \,
\frac{1}{2} \,
\nabla_k w \nabla_k A_{ij} S_{ij} \\
&\hskip-7pt\hskip-7pt& + \,  \frac{1}{2} \, \nabla_k w \, 
\nabla_k A^0_{ij} \, S_{ij} \, - \,
S_{ij} \, \nabla_i | \nabla w|^2 \, \nabla_j w \nonumber\\
&\hskip-7pt\hskip-7pt& + \, \frac{1}{2} \,
R \langle \nabla w , \nabla | \nabla w|^2 \rangle   +\,  R_{ikjm} \, S_{ij} \, \nabla_m w \nabla_k w ,\nonumber 
\end{eqnarray} 
\noindent
where $R_{ikjm}$ denotes the components of the curvature tensor of $g${\rm .}
\endproclaim

\demo{Proof} Clearly,
$\nabla_j V = \nabla_j \left( \frac{1}{2} | \nabla w|^2 \right) = 
\nabla_j \nabla_k w \nabla_k w$.  Thus, $\nabla_i \nabla_j V = \nabla_i
\nabla_k w \nabla_j \nabla_k w + \nabla_i \nabla_j \nabla_k w \nabla_k
w$.  Since the Hessian is symmetric, $\nabla_i \nabla_j \nabla_k w =
\nabla_i \nabla_k\nabla_j w$.  Commuting derivatives,  we find
\begin{eqnarray*}
\nabla_i \nabla_k \nabla_j w & = & 
\nabla_k \nabla_i \nabla_j w \, + \, R_{ikjm} \, \nabla_m w. 
\end{eqnarray*}
\noindent
Therefore,
\begin{equation}
\nabla_i \nabla_j V \, = \,
\nabla_i \nabla_k w \nabla_j\nabla_k w \, + \,
\nabla_k \nabla_i \nabla_jw \nabla_k w \, + \,
R_{ikjm} \, \nabla_m w \nabla_k w .\qquad
\end{equation}    
\noindent
Note that by (1.16),
\begin{equation}
\nabla_i \nabla_j w \, = \,
- \, \frac{1}{2} \,
A_{ij} \, + \,
\frac{1}{2} \,
A^0_{ij} \, - \,
\nabla_i w \nabla_j w \, + \,
\frac{1}{2} \, | \nabla w |^2 g_{ij} .
\end{equation}
\noindent
Hence,
\begin{eqnarray*}
\nabla_k \nabla_i \nabla_j w & = &
- \, \frac{1}{2} \, \nabla_k A_{ij} \, + \,
\frac{1}{2} \,
\nabla_k \, A^0_{ij} \, - \,
\nabla_k \nabla_i w \nabla_j w \\
& &- \, \nabla_i w \nabla_k \nabla_j w \, + \,
\frac{1}{2} \, \nabla_k | \nabla w|^2 \, g_{ij},
\end{eqnarray*}
\noindent
which we substitute into (5.33) to get
\begin{eqnarray}\qquad
\nabla_i \nabla_j V & = & \nabla_i \nabla_k w \, \nabla_j
\nabla_k w \, - \,
\frac{1}{2} \, \nabla_k w \, \nabla_k \, A_{ij} \, + \,
\frac{1}{2} \, \nabla_k w \, \nabla_k \, A^0_{ij}   \\
&& - \, \nabla_i \nabla_k w \nabla_j w \nabla_k w \, - \,
\nabla_j \nabla_k w \nabla_i w \nabla_{k} w\nonumber\\
&&+ \,
\frac{1}{2} \, \nabla_k w \nabla_k | \nabla w |^2 g_{ij}  +   R_{ikjm} \, \nabla_m w \, \nabla_k w .\nonumber
\end{eqnarray} 
Pairing both sides of (5.35) with $S_{ij}$ and using the identity
$S_{ij} g_{ij} = R$, we have (5.32). 
\enddemo

Since we will need to rewrite several of the terms in (5.32), let us
denote
\begin{equation}
S_{ij} \nabla_i \nabla_j V \, = \,
{\rm II}_1 \, + \cdots + \, {\rm II}_6 .
\end{equation}
\pagebreak

\proclaim{Lemma}
\begin{eqnarray}
{\rm II}_1  & \equiv & S_{ij} \, \nabla_i \nabla_k w \, \nabla_j
\nabla_k w \\
& = & - \,
\frac{1}{4} \, {\rm tr}\, E^3 \, + \, \frac{1}{48} \, R|E|^2 \, + \, 
\frac{1}{576} \, R^3 \nonumber \\
&& + \, S_{ij} \, A_{jk} \, \nabla_i w \, \nabla_k w \, - \, | \nabla w|^2
\, \sigma_2 (A) \, + \, \frac{1}{4} \, R | \nabla w |^4 \nonumber \\
& & -\, \frac{1}{2} \,
S_{ij} \, A_{ik} \, A^0_{jk} \, + \,
\frac{1}{4} \,
S_{ij} \, A^0_{ik} \, A^0_{jk} \, - \,
S_{ij} \, A^0_{ik} \, \nabla_j w \nabla_k w \nonumber \\
& &+ \, \frac{1}{2} \, S_{ij} \, A^0_{ij} | \nabla w |^2 .\nonumber 
\end{eqnarray}
\endproclaim

\demo{Proof} By (5.34),
\begin{eqnarray}
\quad{\rm II}_1 & = & S_{ij} \nabla_i \nabla_k w \nabla_j \nabla_k w 
\\
& = & S_{ij} \left(
- \, \frac{1}{2} \, A_{ik} \, +\, \frac{1}{2} \, A^0_{ik} \, - \,
\nabla_i w \nabla_k w \, + \, \frac{1}{2}\, | \nabla w|^2 g_{ik} \right)
\nonumber\\
&&\cdot 
\left( - \, \frac{1}{2} \, A_{j k}  + \, \frac{1}{2} \, A^0_{jk} \, - \, \nabla_j w \nabla_k w \, + \,
\frac{1}{2} \, | \nabla w|^2 \, g_{jk} \right) \nonumber \\
& = & \frac{1}{4} \, S_{ij} \, A_{ik} \, A_{jk} \, - \,
\frac{1}{2} \, S_{ij} \, A_{ik} \, A^0_{jk} \, + \, \frac{1}{4} \,
S_{ij} \, A_{ik}^0 \, A_{jk}^0 \nonumber \\
& &+ \, S_{ij} A_{jk} \nabla_i w \nabla_k w \, - \,
S_{ij} A^0_{ik} \, \nabla_j w \nabla_k w \, - \,
\frac{1}{2} \, S_{ij} \, A_{ij} \, | \nabla w |^2 \nonumber \\
&& + \, \frac{1}{2} \, S_{ij} \, A^0_{ij} | \nabla w |^2 \, + \, \frac{1}{4}
\, R | \nabla w|^4 .\nonumber
\end{eqnarray}
Then (5.37) follows from (5.38) and the following two identities:
\begin{eqnarray}
\frac{1}{4} \, S_{ij} \, A_{ik} \, A_{jk} & =& - \, 
\frac{1}{4} \, {\rm tr}\, E^3 \, + \,
\frac{1}{48} \, R | E |^2 \, + \, 
\frac{1}{576} \, R^3 , 
\\S_{ij} A_{ij} & = & 2 \, \sigma_2 (A) .
\end{eqnarray}
\noindent
To prove (5.39) and (5.40) we simply use the fact that $A_{ij} \, = \,
E_{ij} \, + \,
\frac{1}{12} \, R g_{ij} , \,
S_{ij} \, = \, - \,
E_{ij} \, + \, \frac{1}{4} \, R g_{ij} .$ 
\enddemo

\proclaim{Lemma}
\begin{eqnarray}
{\rm II}_2 & \equiv & - \, \frac{1}{2} \, \nabla_k w \, \nabla_k A_{ij}
\, S_{ij}  \\
& = & - \, \frac{1}{2} \,
\langle \nabla w , \nabla \, \sigma_2 ( A )  \rangle . \nonumber
\end{eqnarray}
\endproclaim

\demo{Proof} We have
\begin{eqnarray*}
\nabla_k A_{ij} \, S_{ij} & = & 
\left(
\nabla_k E_{ij} \, + \, \frac{1}{12} \, 
\nabla_k R \, g_{ij} \right) \,
\left( - E_{ij} \, + \, \frac{1}{4} \, R g_{ij} \right) \\
& = & - \, E_{ij} \, \nabla_k E_{ij} \, + \, \frac{1}{12} \, R \nabla_k
R \end{eqnarray*}
\begin{eqnarray*}
& = & \nabla_k \left( - \, \frac{1}{2} | E|^2 \, + \, 
\frac{1}{24} \, R^2 \right) \, = \,
\nabla_k \, \sigma_2 ( A )  \end{eqnarray*} 
\hfill ${\displaystyle
\Longrightarrow   -   \frac{1}{2} \, \nabla_k w \nabla_k A_{ij} \, 
S_{ij}  =  - \, \frac{1}{2} \,
\langle \nabla w , \nabla \, \sigma_2 (A) \rangle . }$
\enddemo
 
\proclaim{Lemma}
\begin{eqnarray}
{\rm II}_5 & \equiv & \frac{1}{2} \, R \, \langle \nabla w , \nabla | 
\nabla w|^2 \rangle \\
& = & - \, \frac{1}{2} \, R A_{ij} \nabla_i w \nabla_j w \, + \,
\frac{1}{2} R A^0_{ij} \nabla_i w \nabla_j w \, - \,
\frac{1}{2} \, R | \nabla w |^4 . \nonumber
\end{eqnarray} 
\endproclaim

\demo{Proof} From (5.34),
\begin{eqnarray*}
{\rm II}_5 &= &
\frac{1}{2} \, R \langle \nabla w , \nabla | \nabla w|^2 \rangle   = 
R \nabla_i \nabla_j w \, \nabla_i w \nabla_j w \\
\\
& = & R \left( - \, \frac{1}{2} \, A_{ij} \, + \, 
\frac{1}{2} \, A^0_{i_j} \, - \,
\nabla_i w \nabla_j w \, + \, \frac{1}{2}
| \nabla w|^2 g_{ij} \right) \, \nabla_i w \nabla_j w \\
& = & -\, \frac{1}{2} \, R A_{ij} \nabla_i w \nabla_j w \, + \,
\frac{1}{2} \, R \, A^0_{ij} \nabla_i w \, \nabla_j w  - \, \frac{1}{2} \, R \, | \nabla w |^4 .\\
\noalign{\vskip-36pt}
\end{eqnarray*}
\enddemo
\vglue8pt

\proclaim{Lemma}
\begin{eqnarray}
{\rm II}_6 & \equiv & R_{ikjm} \, S_{ij} \, \nabla_m w \nabla_k w
\\
& = & W_{ikjm} \, S_{ij} \, \nabla_m w \, \nabla_k w \, - \, S_{ik} \,
A_{jk} \, \nabla_i w \nabla_j w \nonumber \\
&& + \, \frac{1}{2} \,
R \, A_{ij} \, \nabla_i w \, \nabla_j w \, + \,
 \sigma_2 (A) \, | \nabla w|^2 .\nonumber 
\end{eqnarray}
\endproclaim
\demo{Proof} This follows directly from (1.2):
\begin{eqnarray*}
 R_{ikjm} \, S_{ij} \, \nabla_m w \, \nabla_k w 
& = & \left(
W_{ikjm} \, + \, \frac{1}{2} \, g_{ij} \, A_{km} \, - \,
\frac{1}{2} \, g_{im} \, A_{jk} \, - \, \frac{1}{2} \, g_{jk} \, A_{im}
\right. \\
&& \left.\quad + \, \frac{1}{2} \, g_{km} \, A_{ij} \right) \, S_{ij} \,
\nabla_m w \nabla_k w \\
& = & W_{ikjm} \, S_{ij} \, \nabla_m w\, \nabla_k w \, - \,
S_{ik} \, A_{jk} \, \nabla_i w \nabla_j w \\
\\
&& + \, \frac{1}{2} \, R \, A_{ij} \, \nabla_i w \, \nabla_j w \, + \,
\frac{1}{2} \, A_{ij} \, S_{ij} \, | \nabla w |^2 , \\
\end{eqnarray*}
and appealing to (5.40) for the last term we get (5.43). \enddemo 

Combining the results of Lemmas 5.10--5.14 we have: \pagebreak

\proclaim{{C}orollary}
\begin{eqnarray}
S_{ij} \nabla_i \nabla_j V & = & 
- \, \frac{1}{4} \, {\rm tr}\, E^3 \, + \,
\frac{1}{48} \, R | E |^2 \, + \, \frac{1}{576} \, R^3 \\
&& - \, \frac{1}{2} \langle \nabla w, \nabla \sigma_2 (A ) \rangle \, 
 - \, \frac{1}{4} \, R | \nabla w |^4
\nonumber \\
&& - \, S_{ij} \nabla_i | \nabla w|^2 \, \nabla_j w \, + \, W_{ikjm} \,
S_{ij} \, \nabla_m w \nabla_k w \nonumber \\
&& - \, \frac{1}{2} \, S_{ij} A_{ik} \, A^0_{jk} \, + \, \frac{1}{4} \,
S_{ij} \, A^0_{ik} \, A^0_{jk} \nonumber \\
&& - \, S_{ij} \, A^0_{ik} \, \nabla_i w \nabla_k w \, + \, \frac{1}{2}
\, S_{ij} \, A^0_{ij} | \nabla w|^2 \nonumber \\
&& + \, \frac{1}{2} \, \nabla_k w \, \nabla_k \, A^0_{ik} \, S_{ij} \, +
\, \frac{1}{2} \, R \, A^0_{ij} \, \nabla_i w \nabla_j w . \nonumber
\end{eqnarray} 
\endproclaim

\proclaim{Proposition}
\begin{eqnarray}
\qquad S_{ij} \, \nabla_i \nabla_j V & \geq & -   \frac{1}{4} \, {\rm tr}\, E^3   + 
\frac{1}{48} \, R | E |^2   +  \frac{1}{576} R^3  \\
&& - \, \frac{1}{2}  \langle \nabla w , \nabla \sigma_2 ( A ) \rangle
 -  \frac{1}{4}  R | \nabla w|^4
\nonumber\\
&& - \, S_{ij}  \nabla_i | \nabla w|^2  \nabla_j w  -  
C  | {\rm Ric} |^2  - 
C  | {\rm Ric}|  | \nabla w |^2   - C .\nonumber 
\end{eqnarray}
\endproclaim 

\demo{Proof}   We begin with a claim:

\vglue4pt {\it Claim} 5.17.
\begin{eqnarray}
|W| & \leq & C , \\
|A^0| & \leq & C , \\
| \nabla A^0 | & \leq & C \, + \, C | \nabla w | . 
\end{eqnarray} 
\vglue4pt
\advance\theoremcount by 1

The proof of (5.46) and (5.47) is a straightforward application of (3.1):
\begin{eqnarray*}
| W | & = & | W |_g \ = \ | W_0 |_0 \,e^{-2w} \  \leq \ C , \\
| A^0 | & = & | A^0|_g \ = \ | A^0 |_0 \, e^{-2w} \ \leq \ C .
\end{eqnarray*}
\noindent
To prove (5.48), note that in any local coordinate system,
\begin{equation}
\nabla_k A^0_{ij} \, = \, \partial_k A^0_{ij} \, - \, \Gamma^m_{ik} \,
A^0_{mj} \, - \, \Gamma^m_{jk} \, A^0_{im}
\end{equation}
\noindent
where the $\Gamma^m_{ik}$ denote the Christoffel symbols relative to the
metric $g$.  By the transformation law for the Christoffel symbols under
a conformal change of metric (see [Ei]),
\begin{equation}
\Gamma^m_{jk} \, = \, \gamma^m_{jk} \, + \, \partial_j w \delta_{km} \,
\, + \, \partial_k w \delta_{jm} \,
- \,
(g_0)^{ms} \, (g_0)_{jk} \partial_s w
\end{equation} 
\noindent
where the $\gamma^m_{jk}$ denote the Christoffel symbols relative to $g_0$.
From (5.49) and (5.50) we conclude (5.48).

Using (5.46)--(5.48) we can estimate the last seven terms of
(5.44) as
follows:
\begin{eqnarray*}
W_{jkim} \, S_{ij} \nabla_m w \nabla_k w & \geq & - \, C \, | {\rm Ric} | \, |
\nabla w |^2,  
\\
- \, \frac{1}{2} \, S_{ij} \, A_{ik} \, A^0_{jk} & \geq & - \, C \, |
{\rm Ric} |^2 ,  
\\
\frac{1}{4} \, S_{ij} \, A^0_{ik} \, A^0_{jk} & \geq & - \, C \, | 
{\rm Ric} | ,  
\\
- \, S_{ij} \, A^0_{ik} \, \nabla_j w \, \nabla_k w & \geq & - \, C \, |
{\rm Ric} | \, | \nabla w|^2 ,  
\\
\frac{1}{2} \, S_{ij} \, A^0_{ij} \, | \nabla w|^2 & \geq & - \, C \, |
{\rm Ric} | \, | \nabla w |^2 ,  
\\
\frac{1}{2} \, \nabla_k w \, \nabla_k \, A^0_{ij} \, S_{ij} & \geq & -
\, C \, | {\rm Ric} |\, | \nabla w | \, -\, C | {\rm Ric} | \, | \nabla w |^2 , 
\\
\frac{1}{2} \, R \, A^0_{ij} \, \nabla_i w \nabla_j w & \geq & - \, C \,
| {\rm Ric} | \, | \nabla w |^2 ,
\end{eqnarray*}
\noindent
and obtain (5.45). \enddemo

\proclaim{Proposition}
For all $\delta > 0$ sufficiently small{\rm ,}
\begin{eqnarray}
0 \, = \, {\rm II} &  \geq & \int - \frac{1}{4} \, {\rm tr}\, E^3 \, + \,
\frac{1}{288} \, R^3  - \, \frac{1}{4} \, R | \nabla w|^4 \,  \\
& & -  \, C \, \delta R^3 \, - \, C
\delta | \nabla w |^6 - \, C \, R^2 \, - \, C . \nonumber
\end{eqnarray}
\endproclaim 

\demo{Proof}  Integrating by parts and using
(1.14) we have 
\begin{eqnarray}\qquad\quad
\int - \, \frac{1}{2} \langle \nabla w , \nabla \sigma_2 (A) \rangle
& = &
\int \frac{1}{2} \, \Delta w \, \sigma_2 ( A )  
\\
& = & \int \frac{1}{2} \left(
- \, \frac{R}{6} \, + \, \frac{1}{6} \, R_0 \, e^{-2w} \, + \, | \nabla
w |^2 \right) \, \sigma_2 ( A ) 
\nonumber \\
& = & \int - \frac{1}{12} \, R \, \sigma_2 ( A ) \, + \, \frac{1}{12}
\, R_0 e^{- 2w} \, \sigma_2 ( A )  
\nonumber \\
&& + \, \frac{1}{2} | \nabla w|^2 \sigma_2 ( A ) .\nonumber 
\end{eqnarray}
\noindent
Similarly, by (5.34) and (5.40),
\begin{eqnarray} 
\noalign{\vskip-3pt}
&&\\
\noalign{\vskip-3pt}
\int - \, S_{ij} \, \nabla_i | \nabla w|^2 \nabla_j w & = &
\int | \nabla w |^2 \, \nabla_i \, S_{ij} \, \nabla_j w \, + \,
| \nabla w |^2 \, S_{ij} \, \nabla_i \nabla_j w 
\nonumber \\
& = & \int | \nabla w |^2 \, S_{ij} \left\{ - \, \frac{1}{2} \, A_{ij}
\, + \, \frac{1}{2} A^0_{ij} \, - \, \nabla_i w \nabla_j w \right. 
\nonumber \\
&& \left. + \, \frac{1}{2} \, | \nabla w |^2 \, g_{ij} \right\} \nonumber
\\
& = & \int - \, \frac{1}{2} \, | \nabla w |^2 \, S_{ij} A_{ij} \, +
\, \frac{1}{2} \, S_{ij} \, A^0_{ij} \, | \nabla w |^2  
\nonumber \\
&& - \, | \nabla w |^2 \, S_{ij} \, \nabla_i w \, \nabla_j w \, + \,
\frac{1}{2} \, R | \nabla w |^4  
\nonumber \\
& = & \int - \, | \nabla w |^2 \, \sigma_2 ( A ) \, + \, \frac{1}{2}
\, S_{ij} \, A^0_{ij} \, | \nabla w |^2  
\nonumber \\
&& + \, | \nabla w |^2 \, R_{ij} \nabla_i w \, \nabla_j w . \nonumber
\end{eqnarray}
Substituting (5.52) and (5.53) into (5.45) we get 
\begin{eqnarray}
\qquad 0 =  {\rm II} & \geq & \int -  \frac{1}{4} 
{\rm tr} E^3  +  \frac{1}{48}  R | E|^2  +  \frac{1}{576}  R^3
\\
&& - \, \frac{1}{12}  R \sigma_2 ( A )  +  \frac{1}{12}  R_0
e^{-2w}  \sigma_2 ( A )  - 
\frac{1}{2}  | \nabla w|^2  \sigma_2 (A) 
\nonumber \\
& & +\, | \nabla w|^2  R_{ij} \nabla_i w \nabla_j w  -  \frac{1}{4}
R | \nabla w |^4  +  \frac{1}{2} 
S_{ij}  A^0_{ij} | \nabla w|^2 
\nonumber \\
& & -\, C \, | {\rm Ric} |^2 \, - \, C \, | {\rm Ric} | \, | \nabla w |^2 \, - \, C .\nonumber 
\end{eqnarray}
By $(*)_{\delta}$,
\begin{eqnarray*}
\int - \, \frac{1}{12} \, R \, \sigma_2 (A) & = &
\int - \, \frac{1}{12} \, R \left( \frac{\delta}{4} \, \Delta R \, - \, 2
\gamma_1 | \eta |^2 \right)  \\
& = & \int - \, \frac{\delta}{48} \, R \Delta R \, + \,
\frac{\gamma_1}{6} \, R | \eta |^2  \\
& = & \int \frac{\delta}{48} \, | \nabla R |^2 \, + \,
\frac{\gamma_1}{6} \, R | \eta|^2  \\
& \geq & \int \frac{\delta}{48} \, | \nabla R |^2 \, - \, C \, | {\rm Ric} |^2 \, - \, C.
\end{eqnarray*} 
\noindent
Similarly,
\begin{eqnarray*}
\int \frac{1}{12} \, R_0 e^{-2w} \, \sigma_2 (A) & = &
\int \frac{1}{12} \, R_0 e^{-2w} \left( \frac{\delta}{4} \, \Delta R
\, - \, 2 \gamma_1 | \eta|^2 \right)  \\
& =& \int \frac{\delta}{48} \, R_0 e^{-2w} \Delta R \, - \,
\frac{\gamma_1}{6} \, R_0 e^{-2w} | \eta |^2 \\
& = & \int - \, \frac{\delta}{48} \, \langle \nabla R_0 , \nabla R
\rangle \, e^{-2w} \, - \,
\frac{\delta}{48} \, \langle \nabla ( e^{-2w}) , \nabla R \rangle R_0
\\ && - \, \frac{\gamma_1}{6} \, R_0 e^{-2w} | \eta |^2 \\
& \geq & \int - \, C \delta | \nabla R | \, - \, C \delta | \nabla w
| \, | \nabla R | \, - C \\
& \geq & \int - \, \frac{\delta}{96} \, | \nabla R |^2 \, - C |
\nabla w |^2 \, - \, C \\
& \geq & \int - \, \frac{\delta}{96} \, | \nabla R |^2 \, - \, C ,
\end{eqnarray*} 
which when combined with the estimate above gives 
\begin{equation}
 \int - \frac{1}{12} \, R \sigma_2 (A) \, + \, \frac{1}{12} \, R_0
e^{-2w} \, \sigma_2 (A)   \geq \, \int \frac{\delta}{96} \, | \nabla R |^2 \, - \,
C \, | {\rm Ric} |^2 \, - \, C .
\end{equation}
\pagebreak

\noindent
Also,
\begin{eqnarray}
 \int - \, \frac{1}{2} \, | \nabla w |^2 \, \sigma_2 (A)& = &
\int - \,  \frac{1}{2} \, | \nabla w |^2 \, \left(
\frac{\delta}{4} \, \Delta R \, - \, 2 \gamma_1 | \eta |^2 \right)
 \\
& \geq &   \int - \, \frac{\delta}{8} \, | \nabla w |^2 \, \Delta R
\,- \,C |\nabla w|^2 
\nonumber \\
& \geq & \, \int  \frac{\delta}{8} \, 
\langle \nabla | \nabla w |^2 ,
\nabla R \rangle \, - \, C \nonumber \\
& \geq& \int - 
\frac{\delta}{4} \, | \nabla^2 w | \, | \nabla w | \,
| \nabla R  \, - \, C  \nonumber \\
& \geq & \int - \, \frac{\delta}{96} \, 
| \nabla R |^2 \, - \, \frac{3}{2} \, 
\delta | \nabla^2 w |^2 \, | \nabla w |^2 \, - \, C .\nonumber
\end{eqnarray} 

\proclaim{Lemma}
\begin{equation}
\int - \delta \, | \nabla^2 w |^2 \, | \nabla w |^2 \, \gtrsim \,
\int - \delta R^3 \, - \,
\delta | \nabla w |^6 \, - \,
R^2 \, - \, 1.
\end{equation}
\endproclaim

\demo{Proof} From (5.34) we have
\begin{eqnarray*}
| \nabla^2 w |^2 & = & \bigg| - \, \frac{1}{2} \, A_{ij} \, + \,
\frac{1}{2} \, A^0_{ij} \, - \,
\nabla_i w \nabla_j w \, + \, \frac{1}{2} \, | \nabla w |^2 \, g_{ij}
\bigg|^2 \\
& \lesssim & | A |^2 \, + \, | \nabla w |^4 \, + \, 1 .
\end{eqnarray*}
Thus,
\begin{eqnarray*}
\int - \delta \, | \nabla^2 w |^2 \, | \nabla w|^2 \, \geq \,
\int - C \delta \, | A |^2 \, | \nabla w|^2 \, - \, C \delta \, |
\nabla w |^6 \, - \, C \delta .
\end{eqnarray*}  
Since $\sigma_2 (A) \, = \, - \, \frac{1}{2} \, | A |^2 \, + \,
\frac{1}{18} \, R^2$,
\begin{eqnarray*} &&\hskip-56pt
\int - \, \delta \, | \nabla^2 w |^2 \, | \nabla w|^2\\
 & \geq &
\int - \, C \delta \, | \nabla w |^2 \, 
\left( - 2 \, \sigma_2 ( A ) \, + \, 
\frac{1}{9} \, R^2 \right) \, - \, C \delta \, 
| \nabla w |^6 \, - C \delta \\
& = & \int - C \delta \, | \nabla w|^2 \, \sigma_2 (A) \, - \, C
\delta R^2 | \nabla w |^2 \\ 
&& - \, C \delta | \nabla w |^6 \, - \, C \delta \\
& = & \int - \, C \delta \, | \nabla w |^2 \,
\left( \frac{\delta}{4} \, \Delta R \, - \, 2 \gamma_1 | \eta |^2
\right) \, - \,
C \delta R^2 | \nabla w |^2 \\ 
&& - \, C \delta \, | \nabla w |^6 \, - \, C \delta 
\\
& = & \int - \, C \delta^2 \, | \nabla w|^2 \, \Delta R \, - \,
C \delta R^2 \, | \nabla w |^2 \, - \,
C \delta \, | \eta |^2 \, | \nabla w |^2 \\
&& - \, C \delta \, | \nabla w |^6 \, - \, C \delta  
\\
& \geq & \int C \delta^2 \langle \nabla | \nabla w|^2 , \nabla R
\rangle \, - \,
C \delta R^2 \, | \nabla w |^2 \, - \, C \delta \, | \nabla w |^6 \, -
\, C \delta 
\\
& \geq & \int - \,
C \delta^2 \, | \nabla^2 w | \, | \nabla w | \, | \nabla R | \, - \,
C \delta R^2 \, | \nabla w|^2 \, - \, C \delta | \nabla w|^6 \, - \,
C \delta  
\\
& \geq & \int - \, C \delta^2 \, | \nabla R|^2 \, - \,
C \delta^2 \, | \nabla^2 w |^2 \, | \nabla w|^2 \, - \, C \delta R^2 \,
| \nabla w|^2 \\
&& - \, C \delta \, | \nabla w |^6 \, - C \delta.
\end{eqnarray*}  
\noindent
Therefore, when $\delta$ is sufficiently small,
$$
\int - \, \delta \, | \nabla^2 w |^2 \, | \nabla w|^2 \, \geq \,
\int - \,
C \delta^2 \, | \nabla R |^2 \, - \,
C \delta R^2 \, | \nabla w|^2 \, - \, C \delta \, | \nabla w |^6 \, - C
\delta .
$$
By Lemma 5.9,
\begin{eqnarray*}
\int - \, C \delta^2 \, | \nabla R |^2 \, \geq \,
\int - \,
C \delta^2 R^3 \, - \,
C \delta R^2 \, - \,
C \delta ,
\end{eqnarray*}
\noindent
and since $R^2 \, | \nabla w|^2 \, \lesssim \, R^3 \, + \,
| \nabla w |^6$ we conclude
\vglue12pt
\hfill ${\displaystyle
\int - \, \delta \, | \nabla^2 w |^2 \, | \nabla w |^2 \, \gtrsim \,
\int - \,
\delta R^3 \, - \, \delta R^2 \, - \,
  \delta | \nabla w |^6 \, - \, \delta . 
}$\enddemo  
 
Substituting (5.57) into (5.56), then (5.56) and (5.55) into
(5.54) we
get 
\begin{eqnarray}
0 \, = \, {\rm II}  & \geq & \int - \, \frac{1}{4} \, {\rm tr}\, E^3 \, + \,
\frac{1}{48} \, R | E |^2 \, + \, \frac{1}{576} \, R^3 \\
& &+\,  | \nabla w|^2 \, R_{ij} \nabla_i w \nabla_j w \, - \, \frac{1}{4}
\, R | \nabla w |^4 \nonumber \\
&& - \, C \delta R^3 \, - \, C \delta \, | \nabla w |^6 \, - \, C
\nonumber \\
&& +\, \frac{1}{2} \, | \nabla w|^2 \, S_{ij} \, A^0_{ij} \, - \,
C |{\rm Ric} |^2 \, - \, C | {\rm Ric}| \, | \nabla w |^2 .\nonumber 
\end{eqnarray}   

We estimate the last three terms in (5.58) as follows:
\begin{eqnarray*}
\int \frac{1}{2} \, | \nabla w |^2 \, S_{ij} \, A^0_{ij} & \geq &
\int - \,
C \, | {\rm Ric} | \, | \nabla w|^2 ,  
\\
\int - \, C \, | {\rm Ric} |^2 & = & \int C \left( 2 \sigma_2 (A) \, - \,
\frac{1}{3} R^2 \right) 
\\
& \geq & \int C \, - \,CR^2,  
\\
\int - \, | {\rm Ric} | \, | \nabla w |^2 & \geq & \int - \, | {\rm Ric} |^2
\, - \, | \nabla w |^4  
\\
& \geq & \int - \, C R^2 \, - \, C , 
\end{eqnarray*} 
\noindent
the last line  following from (3.2).  Therefore, 
\begin{eqnarray}
0 \, = \, {\rm II} \, \geq \, \int &-&\frac{1}{4} \, {\rm tr}\, E^3 \, + \,
\frac{1}{48} \, R | E |^2 \, + \, \frac{1}{576}\, R^3\\
&+& | \nabla w |^2 \, R_{ij} \, \nabla_i w \nabla_j w \, - \,
\frac{1}{4} \, R | \nabla w |^4 \nonumber \\
&-& C \delta R^3 \, - \, C \delta | \nabla w |^6 \, - \, CR^2 \, - \,
C . \nonumber 
\end{eqnarray}

Using (1.7) and integrating by parts, we can estimate the term involving
the Ricci curvature in (5.59) as follows:
\begin{eqnarray*}
\int | \nabla w|^2 \, R_{ij} \, \nabla_i w \nabla_j w & \geq &
\int \frac{3}{R} \, \sigma_2 (A) \, | \nabla w |^4 \\
\\
& = & \int \frac{3}{R} \, \left(
\frac{\delta}{4} \, \Delta R \, - \,
2 \gamma_1 | \eta |^2 \right) \, | \nabla w |^4 \\
\\
& \geq & \int \frac{3}{4} \, \delta \, \frac{\Delta R}{R} \, | \nabla
w |^4 \\
\\
& = & \int - \, \frac{3}{4} \, \delta \, \nabla R \, \nabla (R^{-1})
\, | \nabla w |^4 \, - \,
\frac{3}{4} \, \delta \, \frac{\nabla R}{R} \, \nabla | \nabla w |^4 \\
\\
& = & \int \frac{3}{4} \, \delta \, 
\frac{| \nabla R|^2}{R^2} \, |
\nabla w |^4 \, - \,
3 \delta \, | \nabla w |^2 \, \nabla^2 w \left(
\frac{\nabla R}{R} , \nabla w \right) \\
\\
& \geq & \int \frac{3}{4} \, \delta \frac{|\nabla R|^2}{R^2} \, |
\nabla w |^4 \, - \,
\frac{3}{4} \, \delta \, \frac{| \nabla R|^2}{R^2} \, | \nabla w |^4 \\
\\
&& - \, 3 \delta | \nabla^2 w |^2 \, | \nabla w|^2 \\
\\
& = & \int - \, 3 \delta \, | \nabla^2 w |^2 \, | \nabla w |^2 .
\end{eqnarray*} 
\noindent
Therefore, by Lemma 5.19,
\begin{equation}
\int \,| \nabla w |^2 \, R_{ij} \, \nabla_i w \nabla_j w \, \gtrsim \,
\int -
\, \delta R^3 \, - \, \delta | \nabla w |^6 \, - \, R^2 \, - \, 1 .\hskip.5in
\end{equation}
\noindent
Substituting (5.60) into (5.59) and using Lemma 5.5(i), we arrive at
(5.51). \hfill\qed\vglue12pt

Through (5.51) we can now complete the proof of Proposition 5.2.  We begin
by adding the results of Proposition 5.18 and Corollary 5.7:
\begin{eqnarray}
0   \geq   {\rm I} \, + \, 24 \; {\rm II}  
& = & \int S_{ij} \, \nabla_i \nabla_j R \, + \, 24 \, S_{ij}
\nabla_i \nabla_j V  \\
& \geq & \int \delta \, \frac{( \Delta R)^2}{R} \, + \, 6 {\rm tr}\, E^3 \, +
\, \frac{1}{12} \, R^3 \nonumber \\
&& - \, 6 {\rm tr}\, E^3 \, + \, \frac{1}{12} \, R^3 \, - \, 6 R | \nabla w |^4
\nonumber \\
&& - \, C \delta R^3 \, - \, C \delta | \nabla w |^6 \, - \, CR^2 \, -
\, C \nonumber \\
& = & \int \delta \, 
\frac{(\Delta R)^2}{R} \, + \, \frac{1}{6} \, R^3 \, - \, 6 R | \nabla w
|^4  \nonumber \\
&& - \, C \delta R^3 \, - \, 
C \delta | \nabla w |^6 \, - \, CR^2 \, - C
.\nonumber
\end{eqnarray}
Dividing through by 36 and using H\"{o}lder's inequality we find 
\begin{eqnarray*}
\int \left( \frac{R}{6} \right)^3 & \leq & \int \frac{R}{6} \, |
\nabla w|^4 \, + \, C \delta R^3 \, + \, C \delta | \nabla w |^6 \, + \,
CR^2 \, + \, C \\
\\
& \leq & \left( \int \left( \frac{R}{6} \right)^3 \right)^{\frac{1}{3}}
\,
\left( \int | \nabla w |^6 \right)^{\frac{2}{3}} \, + \,
\int C \delta R^3 \, + \, C \delta | \nabla w |^6 \\
\\
&& + \, C R^2 \, + \, C .
\end{eqnarray*}
\noindent
Since $x y \, \leq \,
\frac{x^3}{3} \, + \,
\frac{2}{3} \, y^{3/2}$ for $x , y \, \geq \, 0$, 
$$
\int \left( \frac{R}{6} \right)^3 \, \leq \,
\frac{1}{3} \, \int \left(
\frac{R}{6} \right)^3 \, + \,
\frac{2}{3} \, \int | \nabla w |^6 \, + \, C \, \int \left[ \delta R^3 \, + \, \delta | \nabla w |^6 \, + \, R^2 \, + \, 1
\right]
$$
which implies (5.4).

It remains to show how (5.4) leads to the $W^{2,3}$-estimate in Theorem
5.1.  This requires several intermediate estimates, beginning with 
\proclaim{Proposition}
\begin{equation}
\int | \nabla^2 w |^2 \, | \nabla w|^2 \, \lesssim \,
\int \delta | \nabla w|^6 \, + \, R^2 \, + \, 1 .
\end{equation}
\endproclaim 

{\it Proof}.  We will need some preparatory
inequalities

\proclaim{Lemma}
\begin{eqnarray}
&  &{ \rm (i)}\quad
\int \left( \frac{R}{6} \right)^2 \, | \nabla w |^2 \, \leq \,
( 1 + C \delta ) \, 
\int | \nabla w |^6 \, + \,  C \int \,R^2 \, + \, C , \\ 
\nonumber \\
&  &{\rm (ii)} \quad
\int ( \Delta w )^2 \, | \nabla w|^2\leq \, \int 2
\Delta w | \nabla w|^4 \\
&&\hskip1.5in\,  + \, C \delta R^3 
\, + \,  C \delta | \nabla w|^6 \, + \, CR^2 \, + \, C ,\nonumber \\ 
\nonumber \\
&  &{\rm (iii)}\quad
\int | \nabla w |^6 \,\leq \,
\int \frac{1}{6} \, R | \nabla w|^4 \,  + \,  C \delta R^3 \, + \, C
\delta | \nabla w |^6 \,  + \,  CR^2 \, + \, C , \\ 
\nonumber \\
&  &{\rm (iv)} \quad
\int \Delta w | \nabla w |^4 \, \lesssim \,
\int \delta R^3 \, + \, \delta | \nabla w |^6 \, + \, R^2 \, + \, 1 .
\end{eqnarray}
\endproclaim

\demo{Proof}  (i) Using $x y \, \leq \,
\frac{2}{3} \, x^{3/2} \, + \, \frac{1}{3} \, y^3$ along with
(5.4), we have
\begin{eqnarray*}
\int \left( \frac{R}{6} \right)^2 \, | \nabla w|^2 & \leq &
\frac{2}{3} \, \int \left( \frac{R}{6} \right)^3 \, + \, \frac{1}{3}
\, \int | \nabla w |^6 \\
\\
& \leq & ( 1 + C \delta ) \, \int | \nabla w |^6 \, + \, \int C
R^2 \, + \, C .
\end{eqnarray*}
 
(ii)  By (1.14),
\begin{eqnarray*}
&&\hskip-48pt
\int \left( \frac{R}{6} \right)^2 \, | \nabla w |^2  -  \int |
\nabla w |^6\\
 &= &
\int \left( - \Delta w \, + \, | \nabla w|^2 \, + \, \frac{1}{6} \,
R_0 e^{-2w} \right)^2 \, | \nabla w |^2 \, - \, \int | \nabla w |^6
\\
& = & \int ( \Delta w )^2 \, | \nabla w |^2 \, + \,
| \nabla w |^6 \, + \, \frac{1}{36} \, R^2_0 \, e^{-4w} \, | \nabla w |^2
\\
&& - 2 \Delta w | \nabla w |^4 \, - \, \frac{1}{3} \, \Delta w \, R_0
e^{-2w} | \nabla w |^2 \, + \, \frac{1}{3} \, R_0 e^{-2w} | \nabla w |^4
\, - \, | \nabla w |^6 
\end{eqnarray*} 
leads to
\begin{eqnarray}
\int ( \Delta w )^2 \, | \nabla w |^2 & = & \int \left[
\left( \frac{R}{6} \right)^2 \, | \nabla w |^2 \, - \, | \nabla w |^6
\right] 
 \\
& &  +\,\int 2 \Delta w | \nabla w |^4 \, - \, \frac{1}{36} \, R^2_0 \,
e^{-4w} | \nabla w |^2 \nonumber \\
&&+ \, \frac{1}{3} \, \Delta w \, R_0 e^{-2w} | \nabla w |^2 \, - \,
\frac{1}{3} \, R_0 e^{-2w} | \nabla w |^4 .\nonumber
\end{eqnarray}
\noindent
The last two terms in (5.67) can be estimated using (1.14) and (3.2):
\begin{eqnarray*}
\int  \frac{1}{3} \, \Delta w R_0 e^{-2w} | \nabla w |^2 &
\lesssim & \int | \nabla w |^4 \, + \, ( \Delta w )^2 \\
\\
& \lesssim & \int 1 \, + \, R^2 ,\\
\int- \frac{1}{3} \, R_0 e^{-2w} \, | \nabla w |^4& \lesssim&\int | \nabla w |^4 \, \leq \, C .
\end{eqnarray*} 
Finally, appealing to (5.63) we get (5.64).

(iii) By (1.16) and (1.14)
\begin{eqnarray}
 &&\int 2 | \nabla w|^2 \, R_{ij} \nabla_i w \nabla_j w \, = \,
\int 2 | \nabla w|^2 \, A_{ij} \, \nabla_i w \nabla_j w \, + \,
\frac{1}{3} \, R \,  | \nabla w |^4  \\
& = & \int 2 | \nabla w|^2 \, \left\{
A^0_{ij} \, - \, 2 \nabla_i \nabla_j w \, - \, 2 \nabla_i w \nabla_j w
\, + \, | \nabla w|^2 g_{ij} \right\} \, \nabla_i w \nabla_j w
\nonumber \\
&& + \, \frac{1}{3} \, R | \nabla w |^4 \nonumber \\
& = & \int 2 | \nabla w |^2 \, A^0_{ij} \nabla_i w \nabla_j w \, - \,
4 \, | \nabla w |^2 \, \nabla_i \nabla_j w \nabla_i w \nabla_j w \nonumber \\
& &  - \, 2 \, | \nabla w|^6 \, + \, \frac{1}{3} \, R | \nabla w |^4
\nonumber \\
& = & \int 2 \, | \nabla w |^2 \, A^0_{ij} \, \nabla_i w \nabla_j
w \, - \,
\nabla_i \, | \nabla w|^4 \, \nabla_i w \, - \, 2 \, | \nabla w |^6
\nonumber \\
&& + \, \frac{1}{3} \, R \, | \nabla w |^4 \nonumber \\
& = & \int 2 | \nabla w|^2 \, A_{ij}^0 \, \nabla_i w \nabla_j w \, +
\, \Delta w | \nabla w |^4 \, - \, 2 | \nabla w |^6 \, + \, \frac{1}{3}
\, R | \nabla w |^4 \nonumber \\
& = & \int 2 | \nabla w |^2 \, A^0_{ij} \, \nabla_i w \nabla_j w \, +
\,
\left( | \nabla w|^2 \, - \, 
\frac{1}{6} R \, + \, \frac{1}{6}   R_0 e^{-2w} \right) \, | \nabla w
|^4 \nonumber \\
&& - \, 2 | \nabla w |^6 \, + \, \frac{1}{3} \, R | \nabla w |^4
\nonumber \\
& = & \int 2 | \nabla w |^2 \, A^0_{ij} \nabla_i w \nabla_j w \, + \,
\frac{1}{6} \, R_0 e^{-2w} | \nabla w |^4 \nonumber \\
&& - \, | \nabla w |^6 \, + \, \frac{1}{6} \, R | \nabla w |^4 \nonumber \\
& \leq & \int C \, - \, | \nabla w |^6 \, + \, \frac{1}{6} \, R |
\nabla w |^4. \nonumber
\end{eqnarray}    
\noindent
Combining (5.68) with (5.60) we get (5.65).

(iv) Substituting (1.14) into (5.65),  we obtain
\begin{eqnarray*}
\int | \nabla w |^6 & \leq & \int \frac{1}{6} \, R | \nabla w |^4
\, + \, C \delta R^3 \, +\, C \delta | \nabla w |^6  \\
&& + \, C R^2 \, + \, C \nonumber \\
& = & \int \left(
\frac{1}{6} \, R_0 e^{-2w} \, - \, \Delta w \, + \, | \nabla w |^2
\right) 
\, | \nabla w |^4 \nonumber \\
&& + \, C \delta R^3 \, + \, C \delta | \nabla w |^6 \, + \, C R^2 \, +
\, C \nonumber \\
& \leq & \int - \Delta w | \nabla w |^4 \, + \, | \nabla w |^6 \, +
\, C \delta R^3 \, + \, C \delta | \nabla w |^6 \nonumber \\
&& + \, C R^2 \, + \, C \nonumber 
\end{eqnarray*} 
\noindent \hglue.15in${\displaystyle\Longrightarrow
\int \Delta w | \nabla w |^4 \, \leq \, \int C \delta R^3 \, + \,
C \delta | \nabla w |^6 \, + \, CR^2 \, + \, C .}$
\enddemo

Now, combining (5.64) and (5.66) we find
\begin{equation}
\int ( \Delta w )^2 \, | \nabla w |^2 \, \lesssim \,
\int \delta R^3 \, + \,
\delta | \nabla w|^6 \, + \, R^2 \, + \, 1.
\end{equation}
\noindent
By the Bochner formula,
$$
\frac{1}{2} \, \Delta | \nabla w |^2 \, = \, | \nabla^2 w |^2 \, + \,
R_{ij} \nabla_i w \nabla_j w \, + \, \langle \nabla w, \nabla ( \Delta w
) \rangle . $$
\noindent
Multiplying both sides by $| \nabla w |^2$ and integrating by parts
give 
\begin{eqnarray}\qquad\qquad
\int | \nabla^2 w |^2 \, | \nabla w |^2 &\hskip-4pt = \hskip-4pt& \int \frac{1}{2} |
\nabla w |^2 \, \Delta | \nabla w|^2 \, - \, | \nabla w|^2 \, R_{ij}
\nabla_i w \nabla_j w  \\
&\hskip-4pt\hskip-4pt& - \, | \nabla w|^2 \, \langle \nabla w , \nabla ( \Delta w ) \rangle
\nonumber \\
&\hskip-4pt =\hskip-4pt & \int - \, \frac{1}{2} | \nabla | \nabla w |^2 |^2 \, - \, |
\nabla w|^2 \, R_{ij} \nabla_i w \nabla_j w \nonumber \\
&\hskip-4pt\hskip-4pt& + \, | \nabla w |^2 ( \Delta w)^2 \, + \, \Delta w \langle \nabla w ,
\nabla | \nabla w |^2 \rangle 
\nonumber \\
&\hskip-4pt \leq \hskip-4pt& \int - \, \frac{1}{2} \, | \nabla | \nabla w |^2 |^2 \, - \,
| \nabla w |^2 \, R_{ij} \nabla_i w \nabla_j w
\nonumber \\
&\hskip-4pt\hskip-4pt& + \, | \nabla w|^2 ( \Delta w )^2 \, + \, \frac{1}{2} | \nabla w |^2
( \Delta w )^2 \, + \, \frac{1}{2} \, | \nabla | \nabla w |^2 |^2 
\nonumber \\
&\hskip-4pt =\hskip-4pt & \int \frac{3}{2} \, | \nabla w |^2 \, ( \Delta w )^2 \, - \,
| \nabla w |^2 \, R_{ij} \nabla_i w \nabla_j w .\nonumber
\end{eqnarray} 

\noindent
Substituting (5.60) and (5.69) into (5.70) and appealing to
(5.4) give 
(5.62), which concludes the proof of Proposition 5.20. \hfill
\qed\vglue12pt

From Lemma 5.21 we can deduce that the inequality in (5.4) can be
reversed:

\advance\theoremcount by -1
\proclaim{{C}orollary}
\begin{equation}
\int | \nabla w |^6 \, \leq \,
( 1 + C \delta ) \, \int \left( \frac{R}{6} \right)^3 \, + \, \int C R^2 \, + \, C .
\end{equation}
\endproclaim

\demo{Proof} By (5.65),
$$
\int | \nabla w |^6 \, \leq \, \int \frac{R}{6} \, | \nabla w |^4 \,
+ \,
C \delta R^3 \, + \, C \delta | \nabla w |^6 \, + \, C R^2 \, + \, C. $$
\noindent
Once again using the inequality $x y \, \leq \, \frac{2}{3} \, x^{3/2}
\, + \, \frac{1}{3} \, y^3 $,  we have
$$
\int | \nabla w |^6 \, \leq \, \int \frac{1}{3} \,
\left( \frac{R}{6} \right)^3 \, + \,
\frac{2}{3} \, | \nabla w |^6 \, + \, C \delta R^3 \, +\, C \delta |
\nabla w |^6 \, + \, CR^2 \, + \, 1 , $$
\noindent
which implies (5.71). \enddemo

The following reverse-H\"{o}lder inequality will be the penultimate
estimate in the proof of Theorem 5.1. 

\proclaim{Proposition}
\begin{equation}
\left( \int | \nabla w |^{12} \right)^{\frac{1}{4}} \, \lesssim \,
\int | \nabla w|^6 \, + \, 1.
\end{equation}
\endproclaim 

\demo{Proof}  In the following, some of our
calculations are done in the background metric $g_0$.  For this reason,
we will carefully distinguish between quantities that are given with
respect to $g$ versus $g_0$.

To begin, recall that the Sobolev embedding theorem implies that
$W^{1,3} \hookrightarrow L^{12}$.  Thus, for any $f \in W^{1,3}$, 
$$
\left( \int | f |^{12} \, dv_0 \right)^{1/4} \, \lesssim \, \int | \nabla_0 f|^3 \, dv_0 \, + \, \int | f |^3 \, dv_0 .
$$
\noindent
If we take $f = | \nabla_0 w | e^{- \frac{2}{3} w}$, then    
\begin{eqnarray*}
\int f^{12} \, dv_0  & = & 
\int | \nabla_0 w|^{12} \, e^{- 8 w} \, dv_0
\\
\\
& = & \int | \nabla w |^{12} \, dv  .
\end{eqnarray*}
\noindent
Therefore,
\begin{eqnarray}
&&\\
\left( \int | \nabla w |^{12} \, dv \right)^{\frac{1}{4}}  & \lesssim
& 
\int | \nabla_0 \left( | \nabla_0 w | e^{- \frac{2}{3}w} \right) |^3
\, dv_0\nonumber \\
&& + \, \int | \nabla_0 w |^3 \, e^{-2w} \, d v _0 \nonumber \\
& = & \int | e^{- \frac{2}{3}w} \, \nabla_0 | \nabla_0 w | \, + \, |
\nabla_0 w | \nabla_0 \left( e^{- \frac{2}{3} w} \right)|^3 \, dv_0
\nonumber \\
&& + \, \int | \nabla_0 w |^3 e^{-2w} \, dv_0 \nonumber \\
& \lesssim & \int \left[ | \nabla^2_0 w |^3 \, e^{-2w} \, + \,
| \nabla_0 w |^6 \, e^{-2w} \, + \,
| \nabla_0 w|^3 \, e^{-2w} \right] dv_0 \nonumber \\
& \lesssim & \int \left[
| \nabla^2_0 w |^3 \, e^{-2w} \, + \,
| \nabla_0w|^6 \, e^{-2w} \right] \, dv_0 \, +\, C . \nonumber
\end{eqnarray}

Comparing the Hessian $\nabla^2$ relative to the metric $g$ with the
Hessian $\nabla^2_0$ relative to $g_0$ we have

$$
| \nabla^2_0 w|^2 \, \lesssim \, 
e^{4w} \, | \nabla^2 w|^2 \, + \, e^{4w} | \nabla w|^4 . $$
\noindent
Thus,
\begin{equation}
| \nabla^2_0 w |^3 \, e^{-2w} \, \lesssim \, e^{4w} | \nabla^2 w |^3 \,
+ \,
e^{4w} | \nabla w|^6 .
\end{equation}
\noindent
Substituting (5.74) into (5.73) we get 
\begin{equation}
\left(
\int | \nabla w|^{12} dv  \right)^{1/4} \, \lesssim \,
\int | \nabla^2 w |^3 \, dv \, + \, \int
| \nabla w |^6 \, dv \, + \,
1 .
\end{equation}
\noindent
By (5.34),
\begin{equation}
| \nabla^2 w |^3 \, \lesssim \, | A|^3 \, + \, | \nabla w |^6 \, + \, 1
, 
\end{equation}
\noindent
so that
\begin{equation}
\left( \int | \nabla w |^{12} \, dv \right)^{1/4} \, \lesssim \, \int \left( | A|^3 \, + \, | \nabla w|^6 \right) \, dv \, + \, 1 . 
\end{equation}
\proclaim{Lemma}
\begin{equation}
\int | A|^3 \, dv \, \lesssim \, \int R^3 \, dv \, + \, 1 .
\end{equation}
\endproclaim

\demo{Proof}  First, notice $|A|^2 = |E|^2 +
\frac{1}{36} R^2$ implies that
\begin{equation}
\int | A|^3 \, dv \, \lesssim \, \int \left( | E|^3 \, + \, R^3 \right) \, dv .
\end{equation}
\noindent
By $(*)_{\delta}$,
$$
| E|^2 \, = \,
\frac{1}{12} \, R^2 \, + \,
4 \gamma_1 | \eta|^2 \, - \, \frac{\delta}{2} \, \Delta R . $$
\noindent
Multiplying by $|E|$ and integrating by parts gives 
\begin{eqnarray}
\int | E|^3 \, dv & = & \int \left[ \frac{1}{12} \, R^2 | E | \, +
\, 4 \gamma_1 | \eta|^2 |E| \, - \,
\frac{\delta}{2} \, | E | \Delta R \right] \, dv \\
& \leq & \int \left[ \frac{1}{12} \, R^2 | E | \, + \,
\frac{\delta}{2} \, \nabla | E | \nabla R \right] \, dv \nonumber \\
& \leq & \int \left[
\frac{1}{12} \,
R^2 |E| \, + \,
\frac{\delta}{2} \, |\nabla E|^2 \, + \, \frac{\delta}{2} \, | \nabla
R|^2 \right] \, dv .\nonumber 
\end{eqnarray}
Using the inequality $R^2 |E| \, \leq \,
\frac{2}{3} \, R^3 \, + \,
\frac{1}{3} \, | E |^3$, we conclude from (5.80) that
\begin{equation}
\int | E |^3 \, dv \, \lesssim \,
\int \left[
\delta | \nabla E|^2 \, + \,
\delta | \nabla R|^2 \, + \,
R^3 \right] \, dv .
\end{equation}
\noindent
Integrating (5.10) over $M^4$ we obtain the identity
\begin{eqnarray*}
\int | \nabla E |^2 \, dv & = & \int \left[
\frac{1}{12} \, | \nabla R|^2 \, - \, 2 {\rm tr}\, E^3 \, - \, \frac{1}{3} \, R |
E |^2 \right. \\
\\
&& + \, 2 W_{ijk \ell} \, E_{ik} \, E_{j \ell} \, + \, 2 E_{ij} \, B_{ij}
\bigg] \, dv , 
\end{eqnarray*}
so that
$$
\int | \nabla E|^2 \, dv \, \leq \,
\int \left[
\frac{1}{12} \, | \nabla R|^2 \, + \,
C |E|^3 \, + \, C \right] \, dv.
$$
\noindent
Substituting this into (5.81) gives 
\begin{eqnarray*}
\int | E|^3 \, dv & \lesssim & \int \left[ R^3 \, + \,
\delta | \nabla R|^2 \, + \,
\delta |E|^3 \, + \, C \right] \, dv
\\
\Longrightarrow
\int | E |^3 \, dv & \lesssim & \int \left[
R^3 \, + \, \delta | \nabla R|^2 \, + \, C \right] \, dv .
\end{eqnarray*} 
\noindent
Finally, appealing to Lemma 5.9 along with inequality (5.79) we get
(5.78). \hfill\qed

Substituting (5.78) into (5.77),  we obtain
$$
\left( \int | \nabla w|^{12} dv \right)^{\frac{1}{4}} \, \lesssim \,
\int \left[
R^3 \, + \, | \nabla w|^6 \, + \, 1 \right] \, dv .
$$
\noindent
Then (5.72) follows from (5.4). \enddemo

{\it Remark}. In the remainder of this section we
return to our convention of computing in the metric $g$ and suppressing
the volume form.

\proclaim{Lemma}
\endproclaim
\begin{itemize}
\item[{\rm (i)}] $ \int |\nabla w|^6 \, \leq \, C$,
\item[{\rm (ii)}] $ \int | \nabla w|^{12} \, \leq \, C$, 
\item[{\rm (iii)}] $ \int R^3 \, \leq \, C $, 
\item[{\rm (iv)}] $ \int | A |^3 \, \leq \, C$,
\item[{\rm (v)}] $ \int | \nabla^2 w |^3 \, \leq \, C$,
\item[{\rm (vi)}] $\parallel w \parallel_{C^\alpha} \, \leq \, C$, for $\alpha
< \frac{1}{3}$ .
\end{itemize}

\demo{Proof} (i) \ 
We begin with integration by parts:
\begin{eqnarray*}
\int | \nabla w|^6 & = & \int \nabla w \nabla w | \nabla w|^4 \, = \,
\int - w \Delta w | \nabla w|^4 \, - \,
w \nabla w \nabla | \nabla w |^4 \\
& \lesssim & \int | w | \, | \nabla^2 w | \, | \nabla w |^4 .
\end{eqnarray*}
\noindent
By H\"{o}lder's inequality,
\begin{eqnarray*}
\int | \nabla w |^6 & \lesssim & \left( \int | \nabla^2 w |^2 \, |
\nabla w|^2 \right)^{\frac{1}{2}} \, \left( 
\int | \nabla w |^6 w^2 \right)^{\frac{1}{2}} \\
& \leq & \left( \int | \nabla^2 w |^2 \, | \nabla w |^2
\right)^{\frac{1}{2}} \,
\left( \int | \nabla w|^{12} \right)^{\frac{1}{8}} \left(
\int | \nabla w |^4 \, |w|^{\frac{8}{3}} \right)^{\frac{3}{8}} .
\end{eqnarray*} 
Appealing to (5.62), (5.72) and (3.4), we obtain
\begin{eqnarray*}
\int | \nabla w |^6 & \lesssim & \left( \int \delta
|\nabla w|^6 \, + \,  R^2 \, + \, 1 \right)^{\frac{1}{2}} \,
\left( \int |\nabla w|^6 + \, 1 \right)^{\frac{1}{2}}
  \\
& \lesssim & \delta^{\frac{1}{2}} \left( \int |\nabla w|^6
\right) \,+ \, \left( \int \,R^2 \, \right)^{\frac{1}{2}}
\left( \int |\nabla w|^6 \, \right)^{\frac{1}{2}} \\
\\
& \, & + \, \left( \int \,R^2 \, \right)^{\frac{1}{2}} +
 \left( \int |\nabla w|^6 \right)^{\frac{1}{2}} + \, 1 ,\nonumber
\end{eqnarray*}
which implies
\begin{equation}
\int | \nabla w |^6  \lesssim \left( \int R^2 \, + 1 \,
\right).
\end{equation}
By (5.4)
\begin{eqnarray*}
 \left( \int R^2 \, \right) \lesssim \left( \int R^3 \,
\right)^{\frac{2}{3}} \lesssim \left( \int |\nabla w|^6 \,
\right)^{\frac{2}{3}} \, + 1 .
\end{eqnarray*}
And combining this with (5.82) gives (i).

\begin{itemize}
\item[{\rm (ii)}] This is immediate from (i) and (5.72).
\item[{\rm (iii)}] follows from (i) and (5.4).
\item[{\rm (iv)}] follows from (iii) and (5.78).
\item[{\rm (v)}] follows from (i), (iv), and (5.76).
\item[{\rm (vi)}] Notice that by (3.3),
\end{itemize}
\begin{eqnarray*}
\int | \nabla_0 w|^6 \, dv_0 & = & \int | \nabla w|^6 \, e^{2w} \,
dv \\
& \leq & \left( \int | \nabla w|^{12} dv \right)^{\frac{1}{2}} \, \left(
\int e^{4w} dv \right)^{\frac{1}{2}} \\
& \lesssim & \left( \int e^{8 w} \, dv_0 \right)^{\frac{1}{2}} \,
\leq \, C.
\end{eqnarray*}
\noindent
The result then follows from the Sobolev embedding theorem. \enddemo

\demo{Proof of Theorem {\rm 5.1}}   Simply apply Lemma
5.24, using (5.74).
\enddemo

In conclusion, we note an important corollary.
\proclaim{{C}orollary}
\endproclaim
\begin{equation}
\int \delta \left( \frac{\Delta R}{R} \right)^2 \, dv \, \leq \, C .
\end{equation}

\demo{Proof}  From (5.61) and Lemma 5.24 we
conclude that 
$$
\int \delta \frac{(\Delta R)^2}{R} \, \leq \, C .
$$
\noindent
Therefore,
\begin{equation}
\int \delta \left( \frac{\Delta R}{R} \right)^2 \, \leq \, \left(
\frac{1}{\min R} \right) \, \int \delta \frac{(\Delta R)^2}{R} .
\end{equation}
\noindent
Then (5.83) follows from (5.84) and the following result:
\proclaim{Lemma}
$$
\min R \, \geq \, C_0 > 0 . 
$$
\endproclaim

\demo{Proof}  By $(*)_\delta$,
$$
\delta \Delta R \, = \,
8 \gamma_1 | \eta|^2 \, - \, 2 | E |^2 \, + \, \frac{1}{6} \, R^2 \,
\leq \, 8 \gamma_1 | \eta|^2 \, + \, \frac{1}{6} R^2 .
$$
Evaluating at the minimum point of $R$ and appealing to Lemma 5.24 (vi),
we have
\begin{eqnarray*}
( \min R )^2 & \geq & - 48 \gamma_1 \, \min | \eta|^2 \\ 
& = & - \, 48 \gamma_1 \, \min \, e^{-4w} | \eta|^2_0 \\
& \geq & 48 ( - \gamma_1 ) \, ( \max e^{4w} )^{-1} \, ( \min | \eta
|^2_0 ) \\
& \geq & C_0 \, > \, 0 . \\
\noalign{\vskip-36pt}
\end{eqnarray*}
\enddemo
\pagebreak

 \section{{\it A priori} $W^{2,s}$ estimates for $s < 5$}

This section will be an extension of the {\it a priori} estimates of
 Section 5.
Our goal is to modify the argument to establish the following:

\proclaim{Theorem}
Let $g = e^{2w} g_0$ be a solution of $(*)_\delta$ with positive scalar
curvature{\rm ,} normalized so that $\int w dv_0 = 0${\rm .}  Assume 
\begin{equation}
\int \sigma_2 ( A_0) \, dv_0 \, = \,
\int \sigma_2 (A) \, dv \, > \, 0.
\end{equation}
\noindent
Then there are constants $C_s = C ( g_0 , s )$ and $\delta_0 < 1$ such
that 
\begin{equation}
\int | \nabla^2_0 w |^s \, d v_0 \, \leq \, C_s
\end{equation}
for any $0 < s < 5$ and $0 < \delta \leq \delta_0${\rm .}
\endproclaim  

As a direct corollary of the Sobolev embedding theorem,
we have the following ${{\cal{C}}^{ 1 , \alpha}}$ {\it a priori} bound for the
solution $w$ of $(*)_\delta$.

\proclaim{{C}orollary}
Under the assumptions of Theorem {\rm 6.1,} there is a constant 
$C_\alpha = C ( g_0 , \alpha)${\rm ,} so that
\begin{equation}
\parallel w \parallel_{{\cal{C}}^{1, \alpha}} \, \leq \,
C_\alpha \ {\rm for \ all} \ \alpha \, < \, \frac{1}{5} \ {\rm and } \
\delta \, \leq \, \delta_0 .
\end{equation}
\endproclaim

The proof of (6.2) follows the same pattern as the proof of (5.2) in
Section~5, with the exception that the terms contributed by $\delta
\Delta R$ are more complicated than before and    need to be handled   with more care.

To start the proof, in analogy with (5.8) and (5.31), for each $0\leq p \leq 2$, 
define 
\begin{eqnarray}
{\rm I}^p & = & 
\int S_{ij} \, \nabla_i \nabla_j \, R^{p+1} , \\
\nonumber \\
{\rm II}^p & = & \int S_{ij} \, \nabla_i ( R^p \nabla_j V ) ,
\end{eqnarray}
\noindent
where $V = \frac{1}{2} | \nabla w  |^2$.  Since $S$ is divergence-free
(see (5.5)), both ${\rm I}^p$ and ${\rm II}^p = 0$.  Thus our strategy is to
show that some combination of the terms ${\rm I}^p$ and ${\rm II}^p$ is bounded below by a
multiple of $\int R^{p+3}$ plus some lower order terms. 

We begin by splitting the terms ${\rm I}^p$, ${\rm II}^p$ as follows:
\begin{eqnarray}
{\rm I}^p & = & \int S_{ij} \, \nabla_i \, \nabla_j \, R^{p+1} \, =
\,
\int S_{ij} \, \nabla_i \, ( (p + 1 ) \, \nabla_j R \, R^p )
\\
& = & ( p + 1 ) \, \int R^p \, S_{ij} \, \nabla_i \nabla_j R \, + \,
p ( p+1 ) \, \int R^{p-1} \, S_{ij} \nabla_i R \, \nabla_j R \nonumber \\
& = & {\rm I}^p_1 \, + \, {\rm I}^p_2 ;\nonumber 
\end{eqnarray}  
\noindent
and 
\begin{eqnarray}
{\rm II}^p & = & \int S_{ij} \, \nabla_i ( R^p \, \nabla_j V )
\\
& = & p \, \int R^{p-1} \, S_{ij} \, \nabla_i R \, \nabla_j V \, + \,
\int R^p \, S_{ij} \, \nabla_i \nabla_j V \nonumber \\
& = & {\rm II}^p_1  \, + \, {\rm II}^p_2 .\nonumber 
\end{eqnarray} 

\demo{Estimating ${\rm I}^p_1$} 
We now apply the identity (5.10), the estimate (5.15), and argue as in
the proof of Lemma 5.5 to obtain
\begin{eqnarray*}
{\rm I}^p_1 & = & ( p + 1) \, \int R^p \, S_{ij} \, \nabla_i \nabla_j
R  \\
& = & ( p + 1 ) \, \int R^p \, \left[ 3 \Delta \sigma_2 ( A ) \, +  \, 3 ( | \nabla E |^2 -
\frac{1}{12} (\nabla R|^2 ) \right. \nonumber \\
&&  + \, \left.6 {\rm tr}\, E^3 \, + \, R |E|^2 \, - \, 6 W_{ikj\ell} \, E_{ij} \, E_{k
\ell} \, - \, 6 B_{ij} \, E_{ij} \right] \nonumber \\
& \geq & 3 (p + 1 ) \, \int \Delta (R^p) \, \sigma_2 (A) \nonumber
\\
&& + \, (p + 1 ) \, \int \left[ \frac{3}{2} \, \delta \, R^{p-1} \, (\Delta R)^2 \, + \,
\frac{3}{2} \, \delta p \, R^{p-2} | \nabla R|^2 \, \Delta R \right. \nonumber
\\
&&+\,  \left. 12 \gamma_1 \, R^{p-1} \, < \, \nabla R, \, \nabla
| \eta |^2 \, > \, - 12 \gamma_1 \, R^{p-2} \, | \eta |^2 \,
| \nabla R|^2 \right] \nonumber\\
&& + \, (p + 1) \, \int ( 6 R^p \, {\rm tr}\, E^3 \, + \, R^{p+1} \, | E |^2 )
 - \, C \, \int R^{p+2} \, - C .\nonumber
\end{eqnarray*}
\noindent
Using $(*)_\delta$ we obtain   
\begin{eqnarray}\qquad
{\rm I}^p_1 & \geq & 3 ( p + 1 ) \, \int \Delta (R^p) \,
\left[\frac{\delta}{4} \, \Delta \, R -2\gamma_1 |\eta|^2 \right] \\
&& + \ ( p + 1 ) \, \int \left[
\frac{3}{2}  \delta \, R^{p-1} \, ( \Delta R )^2 \, + \,
\frac{3}{2} \, \delta p \, R^{p-2} \, | \nabla R|^2 \, \Delta R \right.
 \nonumber  \\
&& \left.  +\ 12 \gamma_1 \, R^{p-1} \, \langle \nabla R, \, \nabla
| \eta |^2 \rangle  - 12 \gamma_1 \, R^{p-2} \, | \eta |^2 \,
| \nabla R|^2 \right]
\nonumber \\
&& +  ( p + 1 ) \, \int \left( 6 R^p \, {\rm tr}\, E^3 \, + \, R^{p+1} \, | E
|^2 \right)
 - \, C \, \int R^{p + 2} \, - C \nonumber \\
& \geq & {\rm I}^p_{1 , \delta}  +
3(p+1)\int -2\gamma_1 \Delta (R^p)|\eta|^2
  \nonumber \\
 \,\,&& + 3(p+1)\int \left[ 4\gamma_1 R^{p-1}<\nabla R, \nabla |\eta|^2>
-4\gamma_1 R^{p-2}|\eta|^2 |\nabla R|^2 \right]\nonumber \\  
&& \ \ + ( p + 1 ) \, \int \left(
6 R^p \, {\rm tr}\, E^3 \, + \, R^{p+1} \, | E |^2 \right)
 - \, C \int R^{p+2} \, - C ,\nonumber 
\end{eqnarray}
where
\begin{eqnarray}
{\rm I}^p_{1,\delta} & = & \frac{3}{4} \,
\delta ( p + 1 ) \, \int \left[
\Delta ( R^p ) \, ( \Delta R ) \, + \, 2 R^{p-1} \, ( \Delta R )^2
\right. 
 \\
&&\left.\hskip1in + \, 2p \, R^{p-2} \, | \nabla R|^2 \, \Delta R
 \right] .\nonumber
\end{eqnarray}
We can estimate the terms involving $\eta$ in (6.8) by integrating
by parts, using the Schwartz inequality, and the fact that $\gamma_1 <0$,
as follows:
\begin{eqnarray*}
&&\hskip-18pt 3(p+1)\int \left[
-2\gamma_1 \Delta(R^p)|\eta|^2+4\gamma_1R^{p-1}<\nabla R, \nabla |\eta|^2>
-4\gamma_1R^{p-2}|\eta|^2|\nabla R|^2\right]\\
&&=\ 3(p+1)\int \left[
2\gamma_1(p+2)R^{p-1}<\nabla R, \nabla |\eta|^2>
- 4\gamma_1 R^{p-2}|\eta|^2|\nabla R|^2\right] \\
&&\geq\ -C \int R^p |\nabla|\eta||^2.
\end{eqnarray*}
 
Since $|\nabla|\eta||^2 = |\nabla (e^{-2w}|\eta|_0)|^2$ and $p \leq 2$,
we can use the results of Lemma 5.24 (i.e., $||\nabla w||_{12} \leq C,\,\,
w \geq -C$) to conclude
$-C \int R^p |\nabla |\eta||^2 \geq -C \int R^{p+2} -C.
$
Substituting this into (6.8) we have
\begin{equation}
{\rm I}_1^p \geq {\rm I}_{1,\delta}^p
+ (p+1)\int (6R^p {\rm tr}\, E^3 + R^{p+1} |E|^2) \, -C\int R^{p+2} -\,C.
\end{equation} 
We now introduce the notation:
\begin{eqnarray}
A_p & =& \int R^{p-1} ( \Delta R )^2,  
 \\
B_p & = & \int R^{p-2} \, | \nabla R|^2 \, \Delta R.\nonumber
\end{eqnarray} 
\noindent
With this notation, we may rewrite (6.9) as
\begin{equation}
{\rm I}^p_{1 , \delta}  =
\frac{3}{4} \,
\delta \, ( p + 1 ) \, \left[
( p + 2 ) \, A_p \, + \, p ( p + 1 ) \, B_p \right].
\end{equation}
\noindent
The following material is fairly technical, but the overall goal is
to establish (6.31) below:
$$
{\rm I}^p_{1, \delta} \, + \,
{\rm I}^p_2 \, + \, 24 ( p + 1 ) \, {\rm II}^p_1 \, \geqslant C \delta
(A_p + C_p ) \, - C \left( \int R^{p+3} \right)^{\frac{p+2}{p+3}} \,
- C,
$$
for any $\delta > 0, p \leq 2$. This will require additional notation
as well. We begin with:

\proclaim{Lemma}
%Denote $\r{\nabla}^2$ 
Denote $\ring{\nabla}^2
R = \nabla^2 R - \frac{1}{4} \Delta R \, g_{ij}${\rm ,} and
\begin{eqnarray*}
C_p & = & \int R^{p-3} \, | \nabla R |^4 , \\
\\
\r{A}_p & = & \int R^{p-1} \, | \r{\nabla}{}^2 R |^2 , \\
\\
D_p & = & \int R^{p-2} \, \nabla_i \nabla_j R \, \nabla_i R \,
\nabla_j R , \\
\\
\r{D}_p & = & \int R^{p-2} \, {\r{\nabla}{}^2}_{ij} R \, 
 \nabla_i R \, \nabla_j R .
\end{eqnarray*} 
\noindent
Then
\begin{eqnarray}
3 p B_p & = & 4 \r{A}_p \, - 3 A_p \, - 2 ( p - 2 ) \, C_p \, + \, 4 ( p
- 2 ) \, \r{D}_p.  \\
&& + \, 4 \, \int {\rm Ric} \, ( \nabla R, \nabla R) \, R^{p-1}.\nonumber 
\end{eqnarray}
\endproclaim

\demo{Proof}   Recall the Bochner identity:
\begin{equation}
\frac{1}{2} \, \Delta \, | \nabla R |^2 \, = \,
| \nabla^2 R|^2 \, + \, {\rm Ric} ( \nabla R , \nabla R ) \, + \,
\langle \nabla R, \nabla \Delta R \rangle .
\end{equation} 
\noindent
Then integration by parts along with (6.14) give 
\begin{eqnarray}
\qquad B_p & = & \int R^{p-2} | \nabla R|^2 \, \Delta R 
\\
& = & \int \Delta \left( R^{p-2} \, | \nabla R|^2 \right) R \nonumber
\\
& = & \int \Delta \left( R^{p-2} \right) \, | \nabla R|^2 \, R + 2 \,
\int \langle \nabla R^{p-2}, \, \nabla | \nabla R|^2 \rangle R
\nonumber \\
&& + \, \int \Delta | \nabla R |^2 \, R^{p-1} \nonumber \\
& = & ( p - 2 ) \, \int R^{p-2} \, |\nabla R|^2 \, \Delta R \, + \, (
p-2) (p-3) \, \int R^{p-3} \, | \nabla R|^4 \nonumber \\
&& + \, 4 ( p-2) \, \int R^{p-2} \, \nabla_i \nabla_j R \, \nabla_i R
\, \nabla_j R 
 + \, 2 \, \int R^{p-1} \, | \nabla^2 R |^2 \nonumber \\
&& + \, 2 \, \int R^{p-1} \, {\rm Ric} ( \nabla R , \nabla R ) \, + \, 2 \,
\int R^{p-1} \langle \nabla R , \, \nabla ( \Delta R ) \rangle .\nonumber 
\end{eqnarray}  
Rewriting the last term in (6.15) and integrating by parts once again,
we obtain
\begin{eqnarray}
\qquad\int R^{p-1} \, \langle \nabla R , \, 
\nabla ( \Delta R ) \rangle & = &
 \frac{1}{p} \, \int \nabla R^p \, 
\nabla ( \Delta R ) 
 = - \frac{1}{p} \, \int \Delta \, R^p \, 
\Delta R \\
& = & - \frac{1}{p} \, ( p \, A_p \, + \, p ( p - 1 ) B_p ) \nonumber\\
& =& - ( A_p \, + \, ( p - 1 ) B_p ) .\nonumber 
\end{eqnarray}
Substituting (6.16) into (6.15), we obtain
\begin{eqnarray}
\qquad B_p &\hskip-4pt  =\hskip-4pt & - 2 A_p - p B_p \, + \, 
(p - 2) \, ( p-3) \, C_p  \\
&\hskip-4pt\hskip-4pt& + \,
4 ( p-2) \, D_p + 2 \, \int R^{p-1} | \nabla^2 R|^2  
 +   2 \, \int R^{p-1} \, {\rm Ric} ( \nabla R , \nabla R) .\nonumber
\end{eqnarray}
\noindent
There are two ways to express the term $D_p$.
First, we can write
\begin{eqnarray}
D_p & = & \int R^{p-2} \, ( \nabla_i \nabla_j R - \, \frac{1}{4} \,
\Delta R g_{ij} ) \, \nabla_i R \, \nabla_j R  \\
&& + \, \frac{1}{4} \,
\int R^{p-2} \, | \nabla R|^2 \, \Delta R  
\nonumber \\
& = & \r{D}_p \, + \, \frac{1}{4} \, B_p ,\nonumber
\end{eqnarray}
\noindent
and substituting (6.18) into (6.17), we get 
\begin{eqnarray}
3 B_p & = & - \frac{3}{2} \, A_p \, + \, 2 \, \r{A}_p \, + \, 4 ( p - 2
) \,
\r{D}_p \, + \, ( p - 2 ) ( p - 3 ) \, C_p  \\
&& + \, 2 \, \int R^{p-1} \, {\rm Ric} \, ( \nabla R , \nabla R) .\nonumber
\end{eqnarray}
Alternatively, we can integrate by parts and express $D_p$ as 
\begin{eqnarray}
D_p  = \frac{1}{2} \, \int R^{p-2} \, \nabla_i R \, \nabla_i |
\nabla R|^2 
 = - \, \frac{1}{2} \, ( B_p + ( p - 2 ) \, C_p ) .
\end{eqnarray} 
Substituting (6.20) back into (6.17), we obtain
\begin{equation}
3 ( p-1) \, B_p = \, - \frac{3}{2} \, A_p + 2 \, \r{A}_p \, -
( p - 1 ) ( p-2) \, C_p + 2 \, \int R^{p-1} \, {\rm Ric} ( \nabla R,\nabla
R).
\end{equation}
Summing (6.19) and (6.21), we obtain the identity (6.13) in the
lemma. \enddemo

\proclaim{{C}orollary} For $p < 2${\rm ,}  
\begin{equation}
3 ( p - \delta ) \, B_p \, \geq \, - 3 \, A_p \, + \, C_p ( 2 - p ) \,
\left(
\frac{1}{2} \, + \, \frac{3}{4} \, p \right) .
\end{equation}
\endproclaim

\demo{Proof} Applying the sharp inequality of
[SW, p.~234], we have
\begin{equation}
| {\r{\nabla}}{^2} R \, ( \nabla R , \nabla R ) | \, \leq \,
\frac{\sqrt{3}}{2} \, | {\r{\nabla}}{^2} R | \, | \nabla R|^2 .
\end{equation}
\noindent
Therefore,
\begin{eqnarray}
4 ( 2 - p) \, | \r{D}_p | & \leq &
2 \sqrt{3} \,
( 2 - p ) \, {\r{A}_p}{^{1/2}} \, C_p^{1/2} 
\\
& \leq & 4 \, \r{A}_p \, + \,
\frac{3}{4} \, ( 2 - p)^2 \, C_p .\nonumber 
\end{eqnarray}

Substituting (6.24) into (6.13), then applying the inequality 
${\rm Ric} ( \nabla R , \nabla R )$ $\geq  \frac{3\sigma_2 (A) }{R}   
 | \nabla R |^2$,  we obtain (6.22). \enddemo

We will now begin to estimate ${\rm II}^p_1$.  Our strategy is to
establish that 
$$
{\rm I}^p_{1, \delta} + {\rm I}^p_2 + 24 ( p + 1 )
{\rm II}^p_1 \, \gtrsim \, \delta ( A_p + C_p ) + \hbox{ lower order terms}
$$
as in (6.31).  

\proclaim{Lemma}
There is a constant $C = C (g_0 )$ such that for any $\varepsilon > 0${\rm ,}
$\eta > 0${\rm ,}  
\begin{eqnarray}\qquad\quad
{\rm II}^p_1  &\geqslant  & -   \frac{1}{2} \, p \varepsilon^2 \, \int R^{p-1} \, S_{ij} \, \nabla_i R \, \nabla_j R 
 - C \, \delta \varepsilon^2 \, \eta \, A_p \, - \, C \delta \varepsilon^2
\eta^{-1} C_p  
\\
&& -\, C \, \varepsilon^{-6} \, \eta^{-1} \left(
\int R^{p+3} \right)^
{\frac{p+1}{p+3}} \, - C p \varepsilon^{-2} \left( \int R^{p+3} \right)^{
\frac{p+2}{p+3}}.\nonumber 
\end{eqnarray}
\endproclaim

\demo{Proof} 
For any $\varepsilon > 0$, we can write ${\rm II}^p_1$ as
\begin{eqnarray}&&\\
{\rm II}^p_1 & = & p \, \int R^{p-1} \, S_{ij} \, \nabla_i R \,
\nabla_j V  
\nonumber\\
& = & \frac{1}{2} \, p \, \int R^{p-1} \, S_{ij} \, \nabla_i \left(
\varepsilon R \, + \,
\frac{1}{\varepsilon } V \right) \, \nabla_j \, 
\left( \varepsilon R \, + \, \frac{1}{\varepsilon} V \right)  
\nonumber \\
&& - \, \frac{1}{2} \, p \varepsilon^2 \, \int R^{p-1} \, S_{ij} \,
\nabla_i R \nabla_j R  - \frac{1}{2} \, p \varepsilon^{-2} \, \int R^{p-1} \, S_{ij} \,
\nabla_i V \nabla_j V .\nonumber 
\end{eqnarray}
We notice that for each $\eta > 0$,  
\begin{eqnarray}
&  &  \int R^{p-1} \, S_{ij} \, 
\nabla_i \left( \varepsilon R \, + \,
\frac{1}{\varepsilon} V \right) \, \nabla_j \, 
\left( \varepsilon R \, + \, \frac{1}{\varepsilon} \, V \right)   \\
& &\qquad\geqslant \frac{3}{4} \, \delta \, \int R^{p-1} \, 
\frac{\Delta R}{R} \, | \nabla 
\left( \varepsilon R \, + \, \frac{1}{\varepsilon} \, V \right) 
\bigg|^2 \nonumber \\
& & \qquad \geqslant - \, \frac{3}{2} \, \delta \, 
\left( 
\int R^{p-1} \, ( \Delta R )^2 
\right)^{1/2} \,
\left( 
\int R^{p-3} \, | \nabla ( \varepsilon R \, + \,
\frac{1}{\varepsilon} V |^4 \right)^{1/2} \nonumber \\
& & \qquad \geqslant  - \, C \, \delta \, A_p^{1/2} \,
\left[
\varepsilon^2 \, C_p^{1/2} \, + \,
\varepsilon^{-2} \, 
\left(
\int R^{p-3} \, | \nabla V |^4 
\right )^{1/2} 
\right] \nonumber \\
& & \qquad \geqslant   - \, C \, \delta \varepsilon^2 \, \eta \, A_p -  C \, \delta
\varepsilon^2 \, \eta^{-1} \, C_p 
 - \, C \, \varepsilon^{-6} \, \eta^{-1} \, \int R^{p-3} \, 
| \nabla V |^4 \, - C .\nonumber
\end{eqnarray}

We now estimate the term $\int R^{p-3} | \nabla V|^4$.  Since $V =
\frac{1}{2} | \nabla w|^2$, we have $\nabla_i V = \nabla_i \nabla_j w
\nabla_j w$, and $| \nabla V | \lesssim | \nabla^2 w | \, | \nabla w|$.
Thus
\begin{eqnarray}
\int R^{p-3} | \nabla V|^4 & \lesssim &
\int R^{p-3} \, | \nabla^2 w |^4 \, | \nabla w |^4
\\
\noalign{\noindent 
and}
\int R^{p-1} \, S_{ij} \nabla_i V \nabla_j V & \lesssim & 
\int R^{p-1} \, | \nabla^2 w |^3 \, | \nabla w |^2   \\
 && + \, \int R^{p-1} \, | \nabla^2 w |^2 \, | \nabla w |^4.\nonumber
\end{eqnarray}

Substituting (6.28) into (6.27), then substituting (6.27), (6.29) into
(6.26), we see that inequality (6.25) in Lemma 6.5 is a direct 
consequence of the following technical lemma. 

\proclaim{Lemma}
Suppose $w$ satisfies {\rm (5.2)} and {\rm (5.3).} Then there is a constant $C = C ( g_0)$
such that
for any nonnegative positive number $a, b , c $ with $s = a + b + \,
\frac{c}{2} \leq 6${\rm ,}  
\begin{equation}
\int R^a \, | \nabla^2 w |^b \, | \nabla w|^c \, \lesssim \, \left(
\int R^s \right)^{ \frac{a+b}{s}} \, + \, C .
\end{equation}
\endproclaim
\pagebreak

\demo{Proof}  Since $|R| \, \lesssim \,
| R_0 | \, + \, | \Delta w | \, + \, | \nabla w |^2$,  by
Holder's inequality
\begin{eqnarray*}
\lefteqn{\int R^a \, | \nabla^2 w |^b | \nabla w|^c} \\
& \lesssim & \int | \nabla^2 w |^{a+b} | \nabla w |^c \, + \, \int | \nabla^2 w |^b \, | \nabla w |^{2 a + c} \\
&& + \, \int | \nabla^2 w |^b \, | \nabla w |^c \\
& \lesssim & \left(
\int | \nabla^2 w |^s \right)^{\frac{a+b}{s}} \,
\left( \int | \nabla w |^{2s} \right )^{\frac{c}{2s}} 
 + \, \left( \int | \nabla^2 w |^s \right)^{\frac{b}{s}} \,
\left( \int | \nabla w |^{2s} \right)^{\frac{2a + c}{2s}} \\
&& + \, \left( \int | \nabla^2 w |^s \right)^{\frac{b}{s}} \,
\left( \int | \nabla w |^{\frac{2cs}{2a+c}} \right)^{\frac{2a+c}{2s}}   .
\end{eqnarray*}
By (5.3), for $s \leq 6$, we have $2s \leq 12$ and $\nabla w \in
L^{12}$; hence
$$
\int R^a \, | \nabla^2 w |^b | \nabla w |^c \, \lesssim \,
\left( \int | \nabla^2 w |^s \right)^{\frac{a+b}{s}} \, + \, C.
$$
To finish the proof of (6.30), we simply observe that since $w$ satisfies
(5.2) (i.e.\ $w$ is bounded) and $s \leq 6$, by elliptic regularity,
\begin{eqnarray*}
\int | \nabla^2 w|^s & \lesssim & \int \left(
| \nabla^2_0 w |^s \, + \, | \nabla_0 w|^{2s} \right) \, dv_0  
 \lesssim \int | \nabla^2_0 w|^s \, dv_0 \, + \, C \\
& \lesssim & \int ( \Delta_0 w)^s \, dv_0 \, + \, C 
 \lesssim \int R^s \, dv_0 \, + \, C . \\
\noalign{\vskip-36pt}
\end{eqnarray*}  
\enddemo

\vglue8pt
\proclaim{Lemma}
There is a constant $C$ such that for each $\delta > 0${\it ,} $p < 2${\rm ,} 
\begin{equation}
{\rm I}^p_{1, \delta} \, + \,
{\rm I}^p_2 \, + \, 24 ( p + 1 ) \, {\rm II}^p_1 \, \geqslant C \delta
(A_p + C_p ) \, - C \left( \int R^{p+3} \right)^{\frac{p+2}{p+3}} \,
- C.
\end{equation}
\endproclaim

\demo{Proof}  Combining (6.6), (6.12) and (6.25), we
have for each $\varepsilon > 0$ small enough so $12 \varepsilon^2 < p+1$, $\eta >  0$,
\begin{eqnarray*}
&&\hskip-24pt {\rm I}^p_{1, \delta} \, + \,
{\rm I}^p_{2} \, + \, 24 ( p+1) \, {\rm II}^p_1 
\nonumber \\
& & \geqslant\, \frac{3}{4} \, \delta ( p + 1 ) \,
[ ( p + 2 ) \, A_p \, + \, p ( p+1) B_p] \nonumber \\
&&\quad + \, \left( p ( p+1) - 12 p ( p + 1 ) \varepsilon^2 \right) \,
\int R^{p-1} \, S_{ij} \, \nabla_i R \, \nabla_j R
 - \, C \delta \varepsilon^2 \eta \, A_p \, \nonumber \\
&&\quad -\, C \delta \varepsilon^2 \,
\eta^{-1} \, C_p
 - \, C \varepsilon^{-6} \, \eta^{-1} \left(
\int R^{p+3} \right)^{\frac{p+1}{p+3}} \, - C p \varepsilon^{-2} \left(
\int R^{p+3} \right)^{\frac{p+2}{p+3}}. 
\end{eqnarray*}
 By the fact that
\begin{eqnarray*}
\int R^{p-1}S_{ij}\nabla_i R\nabla_j R &\geq&
\int 3R^{p-2} \sigma _2(A)|\nabla R|^2\\
&=& \int 3R^{p-2}|\nabla R|^2( \frac{\delta}{4}\Delta R - 2\gamma_1 |\eta|^2)\\
&\geq& \frac {3}{4} \delta B_p,
\end{eqnarray*}
the preceding estimate becomes:
\begin{eqnarray}
&&\\
 && \hskip-16pt {\rm I}^p_{1, \delta} \, + \,
{\rm I}^p_{2} \, + \, 24 ( p+1) \, {\rm II}^p_1 \nonumber \\
& &\quad\geqslant   \frac{3}{4} \, \delta ( p + 1 ) \,
\left[ ( p + 2 ) \, A_p \, + \, p ( p + 1 ) \, B_p \, + \, ( p - 12 p
\varepsilon^2 ) \, B_p \right]
 - \, C \delta \varepsilon^2 \, \eta \, A_p \, \nonumber \\
&&\qquad -\, C \delta \varepsilon^2 \,
\eta^{-1} \, C_p
 - \, C \varepsilon^{-6} \, \eta^{-1}
\left(
\int R^{p+3} \right)^{\frac{p+1}{p+3}} \, - \,
C p \varepsilon^{-2} \left( \int R^{p+3} \right)^{\frac{p+2}{p+3}}.\nonumber 
\end{eqnarray}

Thus if $\delta < 1 \leq p < 2$, we may apply (6.22) and (6.32) to
obtain that for all $\eta > 0$,
\begin{eqnarray}&&\\
&&\hskip-24pt {\rm I}^p_{1, \delta} \, + \, {\rm I}^p_{2} \, + \,
24 ( p + 1 ) \, {\rm II}^p_1    \nonumber\\
& &\geqslant \, 9 \delta \varepsilon^2 \, ( p+1) \, A_p \, + \, a_p \, \delta
C_p \, - \, a_p \delta^2 A_p 
 - \, C \delta \varepsilon^2 \, \eta \, A_p \, \nonumber \\
&&\quad -\, C \delta \varepsilon^2 \,
\eta^{-1} C_p 
 - \, C \varepsilon^{-6} \, \eta^{-1} \, \left(
\int R^{p+3} \right)^{\frac{p+1}{p+3}} \, - \, C p \varepsilon^{-2} \,
\left( \int R^{p+3} \right)^{\frac{p+2}{p+3}} ,\nonumber
\end{eqnarray}
where $a_p$ is a positive constant depending only on $p$.

Thus if we first choose $\eta$ small enough so that $C \eta < 8 ( p +
1)$, and then choose $\varepsilon$ sufficiently small so that $a_p > C
\varepsilon^2 \eta^{-1}$, then for $\delta$ sufficiently small we conclude
from (6.33) that (6.31) holds.  \enddemo

We will now estimate the term ${\rm II}^p_2 = \int R^p S_{ij} \nabla_i \nabla_j
V$. 
\proclaim{Proposition}
There is a constant $C$ such that for $p < 2${\rm ,}   for each $\gamma > 0${\rm ,} 
\begin{eqnarray}
\qquad{\rm II}^p_2 & \geqslant & \int R^p \, \left(
- \frac{1}{4} \, {\rm tr}\, E^3 \, + \, \frac{1}{288} \, R^3 \right)
\\
&& - \, C \gamma \delta A_p \, - \, C \gamma \delta C_p \, - \, 
C \gamma^{-1} \delta \, \int R^{p+3} \, -  C \, \int R^{p+2} \,
- C . \nonumber
\end{eqnarray} 
\endproclaim 

\demo{Proof}  The proof of this proposition
follows the   pattern of the proof of Proposition 5.18.  However, the
estimates are less delicate because we already know $w \in
L^\infty$ and $\nabla w \in L^{12}$ in view of Theorem 5.1.  We
will outline the proof but skip some of the details.

To begin with, we have from (5.45) in Proposition 5.16,
\begin{eqnarray}
R^p \, S_{ij} \, \nabla_i \nabla_j V & \geqslant &
R^p \left( - \, \frac{1}{4} \, {\rm tr}\, E^3 \, + \,
\frac{1}{48} \, R |E|^2 \, + \,
\frac{1}{576} \, R^3 \right)  
 \\
&& - \, \frac{1}{2} \, R^p \, \langle
\, \nabla w , \nabla \sigma_2 ( A ) \rangle
\nonumber \\
&& - \, \frac{1}{4} \, R^{p+1} \, | \nabla w|^4 \, - \,
R^p \, S_{ij} \, \nabla_i | \nabla w|^2 \, \nabla_j w  
\nonumber \\
& & -\, C \, R^p | {\rm Ric} |^2 \, - C \, R^p | {\rm Ric} | \, | \nabla w |^2 \, - C
R^p .\nonumber
\end{eqnarray}
By $(*)_\delta$,  
\begin{eqnarray}
\int R^{p+1} \, |E|^2 & = & \int R^{p+1} \, \left(
\frac{1}{12} \, R^2 \, - \, \frac{\delta}{2} \, \Delta R \, + \, 4
\gamma_1 \, | \eta |^2 \right)  
 \\
& \geqslant & \frac{1}{12} \, \int R^{p+3} \, - C \, \int R^{p+1} ,\nonumber\\
\int R^p \, \langle \nabla w , \nabla \sigma_2 (A) \rangle&= & -
\, \int R^p \, \Delta w \sigma_2 (A) \, - \, \int \nabla \, (R^p) \, \nabla w \, \sigma_2 (A),\nonumber
\end{eqnarray}
so that
\begin{eqnarray}
\qquad&&\bigg| \, \int R^p
\langle \nabla w \, \nabla \sigma_2 (A)\rangle \bigg|  \\ 
&&\lesssim \,
 \bigg| \, \int R^p \, \Delta w \delta \Delta R \, \bigg|
+ \bigg| \, \int R^p \, \Delta w \, | \eta|^2 \bigg| \nonumber \\
& &\quad + \, \int R^{p-1} \, | \nabla R| \, \delta | \Delta R | \, | \nabla
w | \, + \,
2\gamma_1 \int \nabla (R^p)   \nabla w   |\eta|^2  
\nonumber \\
&&  \lesssim \,\delta \, A_p^{1/2} \, \left(
\int R^{p+1} \, ( \Delta w )^2 \right)^{1/2} \, + \, \int R^p | \Delta w | \nonumber \\
&&\quad +\, \delta \, A_p^{1/2} \, C_p^{1/4} \, \left(
\int R^{p+1} | \nabla w |^4 \right)^{1/4} \, -
2\gamma_1 \int R^p\left[
\Delta w + \langle \nabla w, \nabla |\eta|^2\rangle \right] .\nonumber 
\end{eqnarray}

Now,   for $p < 2$,
\begin{eqnarray}
&&\\
\int R^{p+1} ( \Delta w )^2 & \lesssim & \int R^{p+3} \, + \, \int R^{p+1} \, | \nabla w |^4 \, + \, \int R^{p+1} 
\nonumber\\
& \lesssim & \int R^{p+3} \, + \,
\left(
\int R^{p+3} \right)^{\frac{p+1}{p+3}} \, 
\left( \int | \nabla w |^{2 (p+3)} \right)^{\frac{2}{p+3}} 
 + \, \int R^{p+1} \nonumber \\
& \lesssim & \int R^{p+3} \, + \, C .\nonumber 
\end{eqnarray}

Applying a similar argument as (6.38) to each of the terms in (6.37), we
obtain
\begin{equation}
\bigg| \int R^p \, \nabla w \, \nabla \sigma_2 (A) \bigg|  \leq 
\gamma \delta \, A_p \, + \, \gamma^{-1} \, \delta \, \int R^{p+3}
 + \, \gamma \, \delta \, C_p \, + \, C\hskip.25in
\end{equation}
\noindent
for any $\gamma > 0$.

We now observe that   integration of the rest of the terms on the
right-hand side of (6.35) can be estimated similarly to the
corresponding terms in Proposition 5.18.

Combining (6.36), (6.37), (6.40) and our observation above, we obtain the
desired estimate (6.34) in Proposition 6.8. \enddemo

\demo{Proof of Theorem {\rm 6.1}}   As explained before,
our strategy of proof is\break the same as the strategy of proof of Theorem
5.1.  That is, we add up\break ${\rm I}^p + 24 (p+1) {\rm II}^p$, so that the
coefficient of the term $\int R^p {\rm tr}\, E^3$ in the sum becomes zero and the
rest of the terms in the sum are dominated from below by $\int R^{p+3}$.
To be more precise, first we combine (6.6), (6.7), (6.8) and (6.34) to
obtain
\begin{eqnarray}
0 & = & {\rm I}^p \, + \, 24 (p+1) \, {\rm II}^p  
 \\
& = & {\rm I}^p_1 \, + \, {\rm I}^p_2 \, + \, 24 ( p + 1 ) \,
{\rm II}^p_1 \, + \, 24 ( p+1) \, {\rm II}^p_2  
\nonumber \\
& \geqslant & {\rm I}^p_{1,\delta} \, + \,
{\rm I}^p_2 \, + \, 24 ( p + 1 ) \, {\rm II}^p_1  
\nonumber \\
&& + \, ( p + 1 ) \, \left(
\int 6 R^p \, {\rm tr}\, E^3 \, + \, \frac{1}{12} \, R^{p+3} \right) 
\nonumber \\
&& + \, 24 (p+1) \, \left( \int - \, \frac{1}{4} \, R^p \, {\rm tr}\, E^3 \, +
\, \frac{1}{288} \, R^{p+3} \right) 
 - \, C \, \gamma \delta \, A_p \,  
\nonumber \\
&&  - \, C \gamma \delta \, C_p \, - \,
C \gamma^{-1} \delta \, \int R^{p+3} 
 - \, C \, \int R^{p+2} \, - C .\nonumber
\end{eqnarray} 

We then apply (6.31), Lemma 6.7 to   estimate the term ${\rm I}^p_{1,
\delta} + {\rm I}^p_2 + 24 ( p + 1 ) {\rm II}^p_1$ in (6.40) above.  We
now choose $\gamma$ small enough so that in the combined
expression the coefficients of the $\delta A_p$ and $\delta C_p$ terms are
positive.  Thus we conclude that there is a constant $C = C ( g_0, p)$,
so that for all $p < 2$,  
\begin{eqnarray}\qquad
0   =   {\rm I}^p  +  24 ( p + 1 )  {\rm II}^p  
& \geqslant & \frac{1}{6}  ( p + 1 )  \int R^{p+3} -  C
\delta  \int R^{p+3}  \\
&& -\, C \left( \int R^{p+3}
\right)^{\frac{p+2}{p+3}}  -  C  \int R^{p+2}  - C.\nonumber
\end{eqnarray}
It follows from (6.41) that for $\delta$ sufficiently small and $p < 2$
there is some constant $C = C ( g_0 , p )$ so that $\int R^{p+3} \, \leq
\, C$.  Thus $\int | \Delta w |^{p+3} \leq C$; from this and the fact
that $w \in L^\infty$, $\nabla w \in L^{12}$, we conclude that
(6.2) holds. \enddemo

\section{Smoothing via the Yamabe flow}

\hskip.2in In light of our estimates in Sections 4--6, we now have {\it a priori}
$C^{1 , \alpha}$ bounds for solutions of $(*)_\delta$ with positive
scalar
curvature.  However, for technical reasons we seem to be unable to
improve on this.  For  example, the integral estimates of Section 6
break down when we attempt to establish an $L^p$-bound for the scalar
curvature as soon as $p \geqslant 5$.  In this section we show that once
$p > 4$, we can use the Yamabe flow to smooth solutions of $(*)_\delta$
and obtain metrics with $\sigma_2 (A) > 0$.

\proclaim{Theorem}
Let $g = e^{2w} g_0$ be a solution of $(*)_\delta$ with positive scalar
curvature{\rm ,} normalized so that $\int w dv_0 = 0${\rm .}  Assume $\int \sigma_2
(A) dv > 0${\rm .}  If $\delta$ is sufficiently small{\rm ,} then there is a
smooth conformal metric $h = e^{2v} g$ such that $\sigma_2 (A_h) > 0${\rm .}
\endproclaim 

The proof of Theorem 7.1 is based on estimates for solutions of
the Yamabe flow using parabolic Moser iteration.  The nonlinear nature of
the flow obviously complicates matters, but we will see that our cause
is aided by the fact that the evolution of the quantity $f \equiv
\sigma_2 (A) + 2 \gamma_1 | \eta|^2$ is fairly simple to analyze (more
precisely, we will study the quantity $\frac{f}{R}$; see (7.24)).

A curious feature of the analysis in this section is the necessity of {\it a
priori} $L^p$ bounds for the curvature of the initial data with $p > 4$.
Typically, the smoothing effects of semi-linear heat flows, like the
Yamabe or Ricci flows, only require $p > \frac{n}{2} = 2$.  But to obtain
in addition a positive lower bound for $\sigma_2 (A)$ we actually need
$p > 4$; see the proof of Theorem 7.1 at the end of the section.

We begin with a basic short-time existence result, based on the work of
[Ha], [Ye].
\proclaim{Proposition}
Let $g$ satisfy the hypotheses of Theorem {\rm 7.1.}  Consider
\begin{equation}
\left\{
\begin{array}{l} 
\frac{\partial h}{\partial t} 
\, = \, - \, \frac{1}{3} \, Rh ,  
\\[5pt]
h ( 0, \cdot ) \, = \, g \, = \, e^{2w} g_0 .
\end{array}
\right.
\end{equation}
\noindent
Then there exists a $T_0 = T_0 (g_0)$ such that {\rm (7.1)} has a unique
smooth solution for $t \in [ 0 , T_0)${\rm .}
\endproclaim

\demo{Proof}  On a compact $n$-dimensional
Riemannian manifold, consider the normalized Yamabe flow
\begin{equation}
 \left\{
\begin{array}{l}
\frac{\partial h}{\partial t} \,  =  \, - \, \frac{1}{(n-1)} \, ( R - r ) h
, \\[5pt]
r ( t ) \,  = \, \int R \, dv \, / \int dv , \\[5pt]
h ( 0 , \cdot ) \, = \, h_0 .
\end{array}
\right. \speqnu{7.1$'$} 
\end{equation}
Then (7.1)$^\prime$ is known to admit a unique smooth solution for all
time (see [Ha], [Ye]).  When $n = 4$, (7.1) and (7.1)$^\prime$ differ
only by a rescaling in time and space in order to normalize the volume.
Therefore, our result will follow if we can produce a time interval
(depending on $g_0$ alone) on which the volume of $h$ is controlled.

To this end, let us record some basic consequences of (7.1).

\proclaimtitle{See, for example, [Ch]}
\proclaim{Lemma}
  Under {\rm (7.1),}
\begin{eqnarray}
\frac{\partial}{\partial t} \, dv & =&- \, \frac{2}{3} \, R dv,
\\
\frac{\partial R}{\partial t}& =& \Delta R \, + \, \frac{1}{3} \,
R^2 .
\end{eqnarray}
\endproclaim

{\it Remark}. Since the initial metric $h ( 0 , \cdot )
= g$ has positive scalar curvature, it follows from applying the minimum
principle to (7.3) that the scalar curvature of $h$ remains positive for
as long as the solution exists.  Indeed, by Lemma 5.26 the scalar
curvature must satisfy 
\begin{equation}
R \geq \, C ( g_0 ) > 0 .
\end{equation}
From (7.2) we see that the volume is decreasing:
$$
\frac{d}{dt} \, \int dv \, = \, \int - \, \frac{2}{3} \, R dv \,
< \,  0 .
$$
\noindent
Also,   
\begin{eqnarray}
 \frac{d}{dt} \, \int d v & \geq&
- \, \frac{2}{3} \, \left(
\int R^2 d v \right)^{ \frac{1}{2}} \, \left(
\int d v \right)^{ \frac{1}{2} }
\\
\Longrightarrow
\frac{d}{dt} \, \left( \int dv \right)^{\frac{1}{2}}&\geq& -
\, \frac{1}{3} \, \left( \int R^2 dv \right)^{\frac{1}{2}}.\nonumber
\end{eqnarray} 
By (7.3),
\begin{eqnarray}
\frac{d}{dt} \, \int R^2 \, dv & = &
\int 2R \left( \Delta R \, + \, \frac{1}{3} \, R^2 \right) \, dv \,
+ \, R^2 \left( - \, \frac{2}{3} \, R \, dv \right) \\
& = & \int - \, 2 | \nabla R|^2 \, dv \, \leq \, 0 . \nonumber
\end{eqnarray} 
\noindent
From (7.5) and (7.6) we conclude that 
$$
\left[ {\rm vol}\, ( h ( 0 , \cdot ) )^\frac{1}{2} \, - \,
C_0 t \right]^2 \, \leq \,
{\rm vol}\, h(t) \, \leq \,
{\rm vol}\, ( h ( 0 , \cdot ))
$$
\noindent
where $C_0 = C_0 ( \parallel R_g \parallel_{L^2} )$.  By Lemma 5.24,
$\parallel R_g \parallel_{L^2} \leq C ( g_0 )$, and this completes
the proof. \enddemo

\proclaim{Proposition}
Let $g$ satisfy the hypotheses of Theorem {\rm 7.1.}  Fix $s \in ( 4 , 5 )${\rm .}
Then there is a $T_1 = T_1 ( g_0 ) < T_0$ such that for $ t \leq T_1${\rm ,}
the solution $h = e^{2v} g$ of {\rm (7.1)} satisfies
\begin{itemize}
\item[{\rm (i)}] $\parallel {\rm Ric}_h \parallel_{L^s} \, \leq \, 2 \parallel {\rm Ric}_g
\parallel_{L^s}$;
\item[{\rm (ii)}] 
$\parallel {\rm Ric}_h \parallel_\infty \, \leq \, C_2t^{- \frac{2}{s}}$, 
where $C_2 = C_2 (g_0 )$;
\item[{\rm (iii)}] 
$\parallel v \parallel_\infty \, \leq \, C ( g_0 )$.
\end{itemize} 
\endproclaim 

\demo{Proof}  The proof of Proposition 7.4 relies
on estimates for the Ricci flow derived in [Ya], summarized in the following 
lemma:
\proclaimtitle{See [Ya]}
\proclaim{Lemma}
  Assume that with respect to the metric 
$h = h ( t )${\rm ,}
$0 \leq t \leq T${\rm ,} 
the following Sobolev inequality holds\/{\rm :}\/
\begin{equation}
\left( \int \left|
\varphi \right|^{\frac{2n}{n-2}} \, dv 
\right)^{\frac{n-2}{n}} \, \leq \,
C_S \left[
\int \left| \nabla \varphi 
\right|^2 \, dv \, + \,
\int \varphi^2 \, dv 
\right] , \,
\varphi \in W^{1,2} (M^n) .\quad
\end{equation}
\noindent
Also{\rm ,} let $b$ be a nonnegative function on $M^n \times [ 0 , T ]$ such
that 
$$
\frac{\partial}{\partial t} \, dv \, \leq \, b \, dv .
$$
\noindent
Let $ q > n${\rm ,} and suppose $u \geq 0$ is a function on $M^n \times [
0 , T ]$ satisfying
$$
\frac{\partial u}{\partial t} \, \leq \, \Delta u \, + \, b u ,
$$
\noindent
and that
$$
\displaystyle{\sup_{0 \leq t \leq T}} \, \parallel b \parallel_{L^{q/2}}
\, \leq \, \beta .
$$
\noindent
Given $p_0 > 1${\rm ,} there exists a constant $C = C ( n , q , p_0 , C_S, \beta )$
such that for $0 \leq t \leq T${\rm ,}
\begin{equation}
\parallel u(t, \cdot ) \parallel_\infty \, \leq \,
C e^{C t} t^{- \frac{n}{2p_0}} \,
\parallel u ( 0 , \cdot ) \parallel_{{p_0}} .
\end{equation}
Moreover{\rm ,} given $p \geq p_0 > 1$, the following inequality holds
for $0 \leq t \leq T${\rm :}
\begin{equation}
\frac{d}{dt} \, \int u^p dv \, + \,
\int \left|
\nabla 
\left( 
u^{p/2} 
\right) 
\right|^2 \, dv \, \leq \,
C \, p^{\frac{2n}{q - n}} \,
\int u^p \, dv
\end{equation}
\noindent
where $C = C ( n , q , p_0 , C_S )${\rm .}
\endproclaim

To apply Lemma 7.5, we need  first to  gain control of the Sobolev
constant $C_S$ defined in (7.7).  To this end, let $\widehat{T}_2 \, \leq \,
T_0$ denote the first time at which
\begin{equation} 
\int R^s_h \, d v_h \, = \,
2 \, \int R^s_g \, dv_g ,
\end{equation}
\noindent
(if (7.10) never occurs then take $\widehat{T}_2 = + \infty$)
and define $T_2 = \min \left( \widehat{T}_2 , 1 \right)$. By the definition
of the Yamabe invariant, for $t \leq T_0$ and any $\varphi \in W^{1,2}
(M^4)$,
\begin{equation}
Y(g_0) \, 
\left(
\int \varphi^4 \, dv_h 
\right)^{\frac{1}{2}} \, \leq \,
\int 6 | \nabla \varphi |^2 \, dv_h \, + \,
\int R_h \varphi^2 d v_h .
\end{equation}
\noindent
Using (7.10) and the fact that $s > 2$ we can estimate the second term on the right
in (7.11) as follows:
\begin{eqnarray*}
\int R_h \, \varphi^2 \, dv_h & \leq &
\left( \int R^s_h \, dv_h \right)^{\frac{1}{s}} \,
\left( \int | \varphi|^{\frac{2s}{s-1}} dv_h \right)^{\frac{s - 1}{s}}  
\nonumber \\
& \leq & \left( 2 \, \int R^s_g \, dv_g \right)^{\frac{1}{s}} \,
\left( \int \varphi^4 \, dv_h \right)^{\frac{1}{s}} \,
\left( \int \varphi^2 \, dv_h \right)^{ \frac{s-2}{s}} .
\end{eqnarray*}
\noindent
By (7.10) we therefore have
\begin{eqnarray*}
\int R_h \, \varphi^2 \, dv_h & \leq &
C ( g_0 ) \, \left(
\int \varphi^4 \, dv_h 
\right)^{\frac{1}{s}} \,
\left( \int \varphi^2 \, dv_h 
\right)^{\frac{s-2}{s}} 
\\
& \leq & \frac{1}{2} \, Y (g_0 ) \, 
\left( \int \varphi^4 \, dv_h 
\right)^{\frac{1}{2}} \, + \, C \, \int \varphi^2 \, dv_h .
\end{eqnarray*}
\noindent Substituting this into (7.11), we see that (7.7) holds for $0
\leq t \leq T_2$, with $C_S = C_S (g_0)$.

We now invoke Lemma 7.5 with $u = R, b = R$, $q = 2s > 8$, and $p_0 =
s$.  By (7.9),
\begin{equation}
\frac{d}{dt} \, \int R^s \, dv \, \leq \, C \, \int R^s
\end{equation}
\noindent
where $C = C ( g_0 , s )$.  Integrating (7.12) we obtain 
\begin{equation}
\int R^s_h \, dv_h \, \leq \, e^{Ct} \, \int R^s_g \, dv_g .
\end{equation}
\noindent
If $T_2 < 1$, then taking $t = T_2$ in (7.13) we conclude that
\begin{equation}
T_2  \, = \, \frac{1}{C} \, \log 2 \, \geq \, C ( g_0 ) .
\end{equation}
Also, from (7.8) we see that 
\begin{equation}
R \, \leq \, C t^{- \frac{2}{s}}
\end{equation}
\noindent
for $0 \leq t \leq T_2$, with $C = C ( g_0)$.  Writing $h = e^{2v} g$
and differentiating, we see that $v$ satisfies
\begin{equation}
\left\{
\begin{array}{l}
\frac{\partial v}{\partial t} \, = \, - \, \frac{1}{6} \, R ,
  \\[5pt]
v ( 0 , \cdot ) \, = \, 0 .
\end{array}
\right.
\end{equation} 
\noindent
Integrating (7.16) using (7.15) we find 
\begin{equation}
\parallel v \parallel_\infty \, \leq \, C \, T_2^{1 - \frac{2}{s}} \,
\leq \, C .
\end{equation}
Now let $\widehat{T}_1 \leq T_0$ denote the first time at which 
\begin{equation}
\int | {\rm Ric}_h |^s \, dv_h \, = \, 2 \, \int | {\rm Ric}_g |^s \, dv_g ,
\end{equation}
\noindent
(if (7.18) never occurs take $\widehat{T}_1 = + \infty$)
and define 
$T_1 \, = \, \min 
\left\{ \widehat{T}_1 , T_2 
\right\}$.  In
order to apply Lemma 7.5 to the evolution for the Ricci tensor, we first
derive the evolution equation for $|{\rm Ric}|^2$.

\proclaim{Lemma} Under {\rm (7.1),} 
\begin{eqnarray}
\frac{\partial}{\partial t} \, | {\rm Ric} |^2 & = & \Delta | {\rm Ric} |^2 - \, 2|
\nabla \, {\rm Ric} |^2 \, - \, 4 tr \, Ric^3 \, + \, 3 R | {\rm Ric} |^2
\\
&& - \, \frac{1}{3} \, R^3 \, + \, 4 W_{ikj\ell} \, R_{ij} \, R_{k\ell}\, + 
\, 4 B_{ij} \, R_{ij}\nonumber 
\end{eqnarray}
\noindent
where $B_{ij}$ is the Bach tensor{\rm .}
\endproclaim

\demo{Proof}  A simple calculation gives 
$$
\frac{\partial R_{ij}}{\partial t} \, = \, \frac{1}{3} \, \nabla_i
\nabla_j R \, + \, \frac{1}{6} \, \Delta R \, g_{ij}.
$$
\noindent
Therefore{\rm ,} by {\rm (1.18),}
\begin{eqnarray*}
\frac{\partial}{\partial t} \, R_{ij} & = & \Delta R_{ij} \, - \,
2 R_{ik} \, R_{jk} \, + \, \frac{1}{2} \, | {\rm Ric} |^2 \, g_{ij} \, + \,
\frac{2}{3} \, R \, R_{ij} \\
&& - \, \frac{1}{6} \, R^2 g_{ij} \, + \,
2 W_{ikj \ell} \, R_{k \ell} \,+ \, 2 B_{ij}
\end{eqnarray*}
\noindent
and (7.19) follows. \enddemo

\proclaim{{C}orollary}
For $t \leq T_1${\rm ,}
\begin{equation}
\frac{\partial}{\partial t} \, | {\rm Ric} | \, \leq \, \Delta | {\rm Ric} | \, + \,
C \, | {\rm Ric} |^2 .
\end{equation}
\endproclaim
\demo{Proof}  By (7.19),
\begin{eqnarray*}
\frac{\partial}{\partial t} \, | {\rm Ric} |^2 & \leq & \Delta | {\rm Ric} |^2 \, -
\, 2 | \nabla {\rm Ric} |^2 \, + \, C \, | {\rm Ric} |^3  
\\
&& + \, 4 | W | \, | {\rm Ric} |^2 \, + \, 4 \, | B | \, | {\rm Ric} | .
\end{eqnarray*}
\noindent
Since $| W | \, = \,
| W_h | \, = \, e^{-2v} | W_g| , \, | B | \, = \,
| B_h| \, = \, e^{-4v} | B_g |$, by (7.17) it follows that for $ t \leq
T_1$,
\begin{eqnarray}\qquad
\frac{\partial}{\partial t} \, | {\rm Ric} |^2 & \leq & \Delta | {\rm Ric} |^2 \, - \,
2 | \nabla {\rm Ric} |^2 \, + \, C \, | {\rm Ric} |^3 \, + \, C \, | {\rm Ric} |^2 \, + \,
C \, | {\rm Ric} | 
 \\
& \leq & \Delta | {\rm Ric} |^2 \, - \,
2 | \nabla {\rm Ric} |^2 \, + \, C \, |{\rm Ric} |^3 \, + \,
C \, |{\rm Ric} | .\nonumber
\end{eqnarray} 
By (7.4), $C_0 < R \, \lesssim \, | {\rm Ric} |$, so $|{\rm Ric}| \, \lesssim \,
|{\rm Ric}|^2$.  Applying this inequality to (7.21) we obtain (7.20). \enddemo

As a consequence of (7.20), we may apply Lemma 7.5 with $u = |{\rm Ric}|$, $b
= C | {\rm Ric} |$, $q = 2s > 8$, and $p_0 = s$.  By (7.9), 
\begin{equation}
\frac{d}{dt} \, \int | {\rm Ric} |^s dv \, \leq \, C \, \int | {\rm Ric} |^s
\, dv.
\end{equation}
\noindent
Integrating (7.22) in time gives 
\begin{equation}
\int | {\rm Ric} |^s \, dv \, \leq \, e^{C t} \, \int | {\rm Ric}_g |^s \,
dv_g .
\end{equation}
Now, if  $T_1 < \min \{ 1 , T_2 \}$, then taking $t = T_1$ in (7.23) we
see that $T_1 \, \geq \, C ( \int | {\rm Ric}_g |^s \, dv_g ) = C (
g_0 )$.  On the other hand, if $T_1 \, \geq \min \{ 1 , T_2 \}$,
then by (7.14) we still conclude that $T_1 \geq C ( g_0 )$. 

Finally, note that (7.8) implies part (ii) of Proposition 7.4, thus
completing the proof. 
 \hfill\qed

\proclaim{Proposition}\hskip-9pt
Define $f = \sigma_2 (A) + 2 \gamma_1 | \eta|^2${\rm .}  If $t \leq T_1${\rm ,} then
under~{\rm (7.1)}
\begin{eqnarray}
\frac{\partial}{\partial t} \, \left(
\frac{f}{R} \right) & \geq & \Delta \, \left(
\frac{f}{R} \right) \, + \,
\frac{2}{R} \, {\rm tr}\, E^3 \, + \,
\frac{1}{3} \, | E |^2 \, - \, \frac{1}{3} \, f   
 \\
&& - \, 2 R^{-1} \, W_{ikj \ell} \, E_{ij} \, E_{k \ell} \, - \,
2 R^{-1} \, B_{ij} \, E_{ij} \, - C.\nonumber
\end{eqnarray}
\endproclaim 

{\it Proof}. The proof of (7.24) requires
several intermediate lemmas, beginning with

\proclaim{Lemma}
Under {\rm (7.1),}
\begin{eqnarray} \qquad\quad
\frac{\partial f}{\partial t} & = & \Delta f \, + \, \left(
| \nabla E|^2 \, - \frac{1}{12} \, | \nabla R|^2 \right) \, + \,
2 {\rm tr}\, E^3 \, + \, \frac{1}{3} \, R | E |^2  
 \\
&& - \, 2 W_{ikj \ell} \, E_{k \ell} \, E_{ij} \, - \, 2 B_{ij} \,
E_{ij}  + \, \left(
\frac{4}{3} \, \gamma_1 \, R | \eta |^2 \, - \,
2 \gamma_1 \, \Delta \, | \eta |^2 \right).\nonumber
\end{eqnarray}
\endproclaim

{\it Proof}.   Since
$$
\frac{\partial}{\partial t} \, | \eta |^2 \, = \,
\frac{2}{3} \, R \, | \eta |^2 ,
$$
by combining (7.3) and (7.19) we get (7.25).
\proclaim{Lemma}
\begin{equation}
| \nabla E |^2 \, - \,
\frac{1}{12} \, | \nabla R|^2 \, \geq \, - 2 \bigg\langle \nabla f, \,
\frac{\nabla R}{R} \bigg\rangle \, + \,
2 f \, \frac{| \nabla  R |^2}{R^2} 
\, + \, 4 \, \gamma_1 \, | \nabla | \eta | \, |^2 .\qquad
\end{equation} 
\endproclaim

\demo{Proof}  We argue as we did in Lemma 5.6.  Namely, 
\begin{eqnarray*}
\nabla f & = & \nabla ( \sigma_2 (A) \, + \, 2 \gamma_1 | \eta |^2 ) \\
& = & \nabla \left( - \frac{1}{2} \, | E|^2 \, + \, \frac{1}{24} \, R^2 \, +
\, 2 \gamma_1 | \eta |^2 \right) \\
& = & - | E | \nabla | E | \, + \, \frac{1}{12} \, R \nabla R \, + \, 4
\gamma_1 | \eta | \nabla | \eta | .
\end{eqnarray*}
Therefore,
\begin{eqnarray*}
- \, \bigg\langle \nabla f, \frac{\nabla R}{R} \bigg\rangle & = & 
\frac{|E|}{R} \, \bigg\langle \nabla|E| , \nabla R \bigg\rangle \, - \,
\frac{1}{12} \, | \nabla R|^2 \, - \, 4 \gamma_1 \, \frac{| \eta |}{R}
\, \langle \nabla | \eta | , \nabla R \rangle \\
& \leq & \frac{1}{2} | \nabla E|^2 \, + \,
\frac{1}{2} \, \frac{|E|^2}{R^2} \, | \nabla R|^2 \, - \,
\frac{1}{12} | \nabla R |^2 \\
&& - \, 4 \gamma_1 \, \frac{| \eta |}{R} \, \langle \nabla | \eta | ,
\nabla R \rangle \\
& = & \frac{1}{2} \, | \nabla E|^2 \, + \, \frac{|\nabla R|^2}{R^2} \,
\left[
- f\, + \, \frac{1}{24} \, R^2 \, + \, 2 \gamma_1 | \eta |^2 \right] \\
&& - \, \frac{1}{12} \, | \nabla R|^2 \, - \, 4 \gamma_1 \,
\frac{|\eta|}{R} \, \langle \nabla | \eta | , \nabla R \rangle \\
& = & \frac{1}{2} \, \left(
| \nabla E|^2 \, - \, \frac{1}{12} | \nabla R|^2 \right) \, - \, f \,
\frac{|\nabla R|^2}{R^2} \\
&& + \, 2 \gamma_1 \, | \eta |^2 \, \frac{| \nabla R|^2}{R^2} \, - \, 4
\gamma_1 \, \frac{|\eta|}{R} \, \langle \nabla | \eta | , \nabla R
\rangle \end{eqnarray*}
\begin{eqnarray*}
\phantom{- \, \bigg\langle \nabla f, \frac{\nabla R}{R} \bigg\rangle}& \leq & \frac{1}{2} \left( | \nabla E|^2 \, - \,
\frac{1}{12}
\, |
\nabla R|^2 \right) \, - \, f \, \frac{|\nabla R|^2}{R^2} \\
&& + \, 2 \gamma_1\, | \eta |^2 \, \frac{| \nabla R|^2}{R^2} \, - \,
2 \gamma_1 \, | \nabla | \eta | \, |^2 \\
& = & \frac{1}{2} \left( | \nabla E |^2 \, - \, \frac{1}{12} \, | \nabla
R|^2 \right) \, - \,
f \, \frac{|\nabla R|^2}{R^2} \, - \, 2 \gamma_1 | \nabla | \eta| \, |^2
.\\
\noalign{\vskip-36pt}
\end{eqnarray*} 
\enddemo

\vglue6pt

\proclaim{Lemma}
For $t \leq T_1${\rm ,}
\begin{equation}
\frac{4}{3}\, \gamma_1 R | \eta |^2 \, - \,
2 \gamma_1 \, \Delta | \eta |^2 \, + \,
4 \gamma_1 | \nabla | \eta | \,|^2 \, \geq \, - CR - C.
\end{equation}
\endproclaim 

\demo{Proof} Let us write $h = e^{2v} g = e^{2 ( v + w)} g_0
\equiv e^{2z} g_0 .$  Let $L = \Delta - \frac{1}{6} R$ denote the
conformal Laplacian; then $L_h \varphi = e^{-3z} L_{g_0} ( e^z \varphi )
= e^{-3z} L_0 (e^z \varphi ).$  Therefore, 
\begin{eqnarray}
\Delta | \eta |^2 & = & L_h | \eta |^2 \, + \, \frac{1}{6} \, R | \eta |^2
 \\
& \geq & e^{-3z} \, L_0 ( e^{z} | \eta |^2 ) \nonumber \\
& = & e^{-3z} \, L_0 
\left( 
e^{-3z} | \eta |^2_0 \right) \nonumber \\
& = & e^{-3z} 
\left(
\Delta_0 
\left(
e^{-3z} | \eta |^2_0 
\right) \, - \,
\frac{1}{6} \, R_0 e^{-3z} \, | \eta |^2_0 
\right] \nonumber \\
& = & e^{-3z} \, 
\left[
| \eta|^2_0 \, \Delta_0 
\left( e^{-3z} 
\right) \, + \, e^{-3z} \,
\Delta_0 | \eta|^2_0 
\right. \nonumber \\
\nonumber \\
&& \left.
+ \, 2 \bigg\langle \nabla_0 
\left( e^{-3z} 
\right), \,
\nabla_0 
\left( | \eta |^2_0 
\right) \bigg\rangle 
\, - \,
\frac{1}{6} \, R_0 \, e^{-3z} | \eta |^2_0 
\right] 
\nonumber \\
& = & e^{- 6z} \, | \eta |^2_0 \, \left[
- 3 \Delta_0 z \, + \, 9 | \nabla_0z |^2 \right] \nonumber \\  
&& +\, e^{-6z} \Delta_0 | \eta|^2_0 \, - \,
12 e^{-6z} \,
\langle \nabla_0 z , \nabla_0 | \eta|_0 \rangle \, | \eta |_0 \nonumber
\\
&& \hspace{.5in} - \, \frac{1}{6} \, R_0 e^{-3z} \, | \eta |^2_0 .\nonumber
\end{eqnarray}
By (1.10),
$$
\Delta_0z \, + \,
| \nabla_0 z|^2 \, + \, \frac{1}{6} \, R e^{2z} \, = \, \frac{1}{6} \,
R_0 .
$$
\noindent
Thus,
$$
e^{-6z} \, | \eta|^2_0 \, \left[
- 3 \Delta_0 z \, + \, 9 | \nabla_0 z|^2 \right]  
= \, e^{-6z} | \eta|^2_0 \, \left[
12 | \nabla_0 z|^2 \, + \,
\frac{1}{2} \, R e^{2z} - \, \frac{1}{2} \, R_0 \right]. $$
Substituting this into (7.28), we get
\begin{eqnarray*}
\Delta | \eta|^2 & \geq & 12 \, e^{-6z} \, | \eta |^2_0 \, |
\nabla_0 z|^2 \, - \,
12 \, e^{- 6z} \, \langle \nabla_0 z , \nabla_0 | \eta|_0 \rangle \, |
\eta|_0 \\
&& + \, e^{-6z} \, \Delta_0 | \eta|^2_0 \, - \,
\frac{2}{3} \, R_0 e^{-6z} \, | \eta |^2_0 .
\end{eqnarray*}
\noindent
Since $\gamma_1 < 0$, this implies 
\begin{eqnarray*}
- \, 2 \gamma_1 \Delta | \eta |^2 & \geq & - \, 24 \gamma_1 \,
e^{-6z} \, | \eta |^2_0 \, | \nabla_0 z |^2 \, + \, 24 \gamma_1 \,
e^{-6z} \,
\langle \nabla_0 z , \nabla_0 | \eta |_0 \rangle \, | \eta |_0 \\
&& - \, 2 \gamma_1 \, e^{-6z} \, \Delta_0 | \eta|^2_0 \, + \,
\frac{4}{3} \, \gamma_1 \, R_0 e^{-6z} | \eta|^2_0 .
\end{eqnarray*} 
Now, 
\begin{eqnarray*}
4 \gamma_1 | \nabla |\eta| \, |^2 & = & 4 \gamma_1 e^{- 6z} | \nabla_0 (
e^{-2z} | \eta |_0 ) |^2 \\
& = & 4 \gamma_1 e^{-6z} | \, - 2e^{-2z} \, \nabla_0 z \, | \eta |_0 \,
+ \,
e^{-2z} \, \nabla_0 | \eta|_0 |^2 \\
& = & 16 \gamma_1 \, e^{-6z} | \nabla_0 z|^2 \, | \eta|^2_0 \, - \,
16 \gamma_1 \, e^{-6z} \,
\langle \nabla_0 z, \nabla_0 | \eta|_0 \rangle \, | \eta |_0 \\
 &&
+ \, 4 \gamma_1 e^{-6z} | \nabla_0 | \eta|_0|^2 .
\end{eqnarray*}
Therefore,
\begin{eqnarray*}
&&\hskip-.75in - \, 2 \gamma_1 \Delta | \eta|^2 \, + \, 4 \gamma_1 | \nabla | \eta | \,
|^2\\
& \geq & - \, 8 \gamma_1 \, e^{-6z} | \eta|^2_0 \, | \nabla_0 z|^2
\\
&&+ \, 8 \gamma_1 e^{-6z} \, 
\langle \nabla_0 z , \nabla_0 | \eta|_0 \rangle | \eta|_0  
\\
&& + \, 4 \gamma_1 e^{-6z} | \nabla_0 |\eta |_0 |^2 \, - \, 
2 \gamma_1 \, e^{-6z} \Delta_0 | \eta|^2_0 \, + \,
\frac{4}{3} \, \gamma_1 R_0 e^{-6 z} | \eta |^2_0 
\\
& \geq & 6 \gamma_1 e^{-6z} | \nabla_0 | \eta|_0 |^2 \, - \,
2 \gamma_1 \, e^{-6z} \Delta_0 | \eta|^2_0  + \, \frac{4}{3} \,
 \gamma_1 \, R_0 \, e^{-6z} | \eta|^2_0 .
\end{eqnarray*} 
By Proposition 7.4 and (5.3), $ \parallel z \parallel_\infty \, \leq
\, \parallel w \parallel_\infty \, + \, \parallel v \parallel_\infty \,
\leq C ( g_0 )$.  Hence,
$$
- \, 2 \gamma_1 \, \Delta | \eta |^2 \, + \, 4 \gamma_1 | \nabla | \eta
| \, |^2 \, \geq \, - C . $$
Also,
$$
\frac{4}{3} \, \gamma_1 \, R | \eta|^2  =  \frac{4}{3} \, \gamma_1 \,
R e^{-4z} \, | \eta |^2_0 
 \geq  - CR,
$$
so that (7.27) follows. \enddemo

To complete the proof of (7.24), we compute $\frac{\partial}{\partial t}
\, (\frac{f}{R})$ using the results of Lemmas 7.9 and 7.10:
\begin{eqnarray*}
\frac{\partial}{\partial t} \, \left(
\frac{f}{R} \right) 
& = & R^{-1} \, 
\frac{\partial f}{\partial t} \, + \,
f \, \frac{\partial}{\partial t} \, 
\left( R^{-1} \right) \\
\\
& = & R^{-1} \, 
\frac{\partial f}{\partial t} \, - \, R^{-2} f \, 
\left(
\Delta R \, + \, \frac{1}{3} \, R^2 \right) \\
\\
& \geq & R^{-1} \, \Delta f \, - \,
2 R^{-2} \, \langle \nabla f, \nabla R \rangle \, + \, 2 f \, R^{-3} |
\nabla R|^2 \\
\\
& - & f \, R^{-2} \Delta R \, + \,
2R^{-1} \, {\rm tr}\, E^3 \, + \, 
\frac{1}{3} \, | E |^2 \, - \, \frac{1}{3} f \\
\\
& - & 2 R^{-1} W_{ikj \ell} \, E_{ij} \, E_{k \ell} \, - \, 2 R^{-1} \,
B_{ij} \, E_{ij} \, - \, CR^{-1} \, - \, C .
\end{eqnarray*}
\noindent
Note that
\begin{eqnarray*}
\Delta \, \left( \frac{f}{R} \right) & = &
R^{-1} \Delta f \, +\,
f \Delta ( R^{-1} ) \, + \, 2 \langle \nabla f , \nabla ( R^{-1})
\rangle \\
\\
& = & R^{-1} \Delta f \, - \,
f R^{-2} \, \Delta R \, + \, 2 f R^{-3} | \nabla R |^2 \, - \, 2R^{-2} \,
\langle \nabla f ,   \nabla R ), 
\end{eqnarray*}
\noindent
so that (7.24) follows. \hfill\qed

\proclaim{Proposition}
Define $\varphi = \max \left\{
- \, \frac{f}{R} , 0 \right\}$.  Then for $t \leq T_1${\rm ,}
\begin{equation}
\frac{\partial \varphi}{\partial t}\, \leq \, \Delta \varphi \, + \, C_1
\, |{\rm Ric}|\, \varphi \, + \, C_1 | {\rm Ric} |
\end{equation}
\noindent
where $C_1 = C_1 ( g_0)${\rm .}
\endproclaim 

\demo{Proof}   We begin by analyzing the
curvature terms in (7.24).  As in the proof of Corollary 5.8, we have
the sharp inequality
\begin{eqnarray*}
2 R^{-1} {\rm tr}\, E^3 \, + \, \frac{1}{3} \, | E|^2 & = & \frac{1}{3R} \,
\left(
6 {\rm tr}\, E^3 \, + \, R | E|^2 \right) 
\\
& \geq & \frac{|E|^2}{3R} \, \left( - 2 \sqrt{3} \, | E | \, + \, R
\right).
\end{eqnarray*}
\noindent
Therefore,
\begin{eqnarray} 
&&\\
2R^{-1} {\rm tr}\, E^3 \, + \, \frac{1}{3} \, | E |^2 &
\geq & \frac{|E|^2}{3R (2 \sqrt{3} |E| + R )} \, ( 2 \sqrt{3} |E|
\, + \, R ) \, ( - 2 \sqrt{3} | E| \, + \, R )  
\nonumber\\
& = & \frac{|E|^2}{3R(2 \sqrt{3} |E| + R)} \, ( -12 |E|^2 + R^2 )
\nonumber \\
& = & \frac{8 |E|^2}{R ( 2 \sqrt{3} | E| + R )} \, \sigma_2 (A)
\nonumber \\
& = & \frac{8 |E|^2}{R(2 \sqrt{3} |E| + R)} \, ( f - 2 \gamma_1 |
\eta|^2 ) \nonumber \\
& \geq & \frac{8 | E|^2}{(2 \sqrt{3} | E | + R )} \, \left(
\frac{f}{R} \right) .\nonumber 
\end{eqnarray}
Also, by Proposition 7.4 (iii),
\begin{eqnarray}
&& \hskip-.75in - \, 2 R^{-1} \, W_{ikj \ell} \, E_{ij} \, E_{k \ell} \, - \, 2
R^{-1} \, B_{ij} \, E_{ij} 
 \\
& \geq & - \, C \, \frac{|E|^2}{R} \, - \, C \, \frac{|E|}{R}
\nonumber \\
& \geq & - \, C \, \frac{|E|^2}{R} \, - \,
C \, \frac{1}{R} \nonumber \\
& = & \frac{C}{R} \, \left(
- \, \frac{1}{2} \, | E |^2 \, + \, \frac{1}{24} \, R^2 \right) \, - \,
\frac{C}{24} \, R \, - \, \frac{C}{R} \nonumber \end{eqnarray}
\begin{eqnarray}
& = & \frac{C}{R} \, ( f - 2 \gamma_1 | \eta|^2 )\, - \,
\frac{C}{24} \, R \, - \, \frac{C}{R} \nonumber \\
& \geq & C \, \frac{f}{R} \, - \,
CR \, - \, \frac{C}{R} .\nonumber
\end{eqnarray} 
Finally, combining (7.24), (7.30), (7.31), and using (7.4) we conclude
\begin{eqnarray*}
\frac{\partial}{\partial t} \, \left( \frac{f}{R} \right) & \geq &
\Delta \, \left( \frac{f}{R} \right) \, + \,
\frac{8 | E|^2}{(2 \sqrt{3} |E|  + R )} \, \left(
\frac{f}{R} \right) \, - \, \frac{1}{3} \, R \left( \frac{f}{R} \right)
\\
&& - \, C \, \left( \frac{f}{R} \right) \, - \,
CR \, - \, C.\nonumber 
\end{eqnarray*}
\noindent
If we let $\varphi = \max \{ - \, \frac{f}{R} , 0 \}$, then (in the
$W^{1,2}$-sense)
$$
\frac{\partial \varphi}{\partial t} \, \leq \,
\Delta \varphi \, + \, C ( | E | \, + \, R \, + \,1 ) \varphi \, + \,
CR \, + \, C.
$$
Since 
$| {\rm Ric} | \, \geq \, \frac{1}{2} \, R \, \geq \, C >
0, \, | E | \, + \, R \, + \, 1 \, \lesssim \, | {\rm Ric}|$,
and (7.29) follows. \enddemo 

\demo{Proof of Theorem {\rm 7.1}} 
\hspace{.5em} Let us begin by summarizing
(7.29) and Proposition 7.4(ii): for $t \leq T_1$, 
\begin{equation}
\frac{\partial \varphi}{\partial t} \, \leq \, \Delta \varphi \, + \,
C_1 | {\rm Ric} | \varphi \, + \, C_1 | {\rm Ric} |,
\end{equation}
\begin{equation}
\parallel {\rm Ric} \parallel_\infty \, \leq \, C_2 t^{- \frac{2}{s}} .
\end{equation}
\noindent
Define 
$\varphi_1 = \varphi_1 (t)\, = \, \exp 
\left\{
\frac{s}{s-2} \, C_1 C_2 t^{\frac{s-2}{s}} 
\right\} - 1$. 
Since $s > 4$,
$\varphi_1 (0) = 0$ and it is easily verified that 
$\partial_t \varphi_1 \, = \, C_1 C_2 ( 1 + \varphi_1 ) 
t^{- \frac{2}{s}}$ for $t > 0$. Let $u = \varphi - \varphi_1$. Then
\begin{eqnarray*} 
\frac{\partial u}{\partial t} & = & \frac{\partial \varphi}{\partial t} \, -
\, \frac{\partial \varphi_1}{\partial t}  
\\
& \leq & \Delta \varphi \, + \, C_1 | {\rm Ric} | \varphi \, + \, C_1 | {\rm Ric}|
\, - \, \frac{\partial \varphi_1}{\partial t}  
\\
& = & \Delta u \, + \, C_1 | {\rm Ric} |u \, + \, C_1 |{\rm Ric}| \varphi_1 \, + \,
C_1 |{\rm Ric} | \, - \, \frac{\partial \varphi_1}{\partial t}  
\\
& \leq & \Delta u \, + \, C_1 |{\rm Ric}|u \, + \, C_1 C_2 ( 1 \, + \,
\varphi_1 ) \, t^{- \frac{2}{s}} \, - \, \frac{\partial
\varphi_1}{\partial t} 
\\
& = & \Delta u \, + \, C_1 |{\rm Ric}|u .
\end{eqnarray*} 
Appealing once more to Lemma 7.5 with $b = C_1 |{\rm Ric}|$, $p_0 = 2$, $q =
2s$, we conclude that for $t \leq T_1$,
$$
\parallel \varphi - \varphi_1 \parallel_\infty \, \leq \,
Ct^{- 1} \, \parallel \varphi (0, \cdot ) \, - \,
\varphi_1 ( 0 ) \parallel_{L^2} .
$$
Now, $\varphi_1 ( 0 ) \, = \, 0$, and by $(*)_\delta$
$$
\parallel \varphi (0 , \cdot ) \parallel_{L^2}  \leq  
\parallel \, \sigma_2 ( A ) \, + \,
2 \gamma_1 | \eta |^2 \parallel_{L^2}   =  \parallel \, \frac{\delta}{4} \left(\frac{ \Delta_g R_g}{R_g}\right) \parallel_{L^2}
.
$$
Therefore, by Corollary 5.25,
$$
\parallel \varphi - \varphi_1 \parallel_\infty \, \leq \, C
\delta^{\frac{1}{2}} \, t^{- 1 }
$$
$$ \Longrightarrow
\varphi \, \leq \, \varphi_1 ( t ) \, + \,
C \delta^{\frac{1}{2}} \, t^{- 1} , \,\,\,  t \leq T_1 .
$$
\noindent
By the definition of $\varphi$, this implies
$$
\frac{1}{R} \left(
\sigma_2 (A) \, + \, 2 \gamma_1 | \eta|^2 \right) \, \geq \,
- \, \varphi_1 ( t ) \, - \, C \delta^{\frac{1}{2}} \,
t^{- 1} ,\,\,\,  \ t \leq T_1 . $$
\noindent
By Taylor's Theorem, $\varphi_1 ( t) \, \leq \, C t^{1- \frac{2}{s}}$.
Also, by (7.33), $R \leq C t^{- \frac{2}{s}}$.  Therefore,
\begin{eqnarray*}
\sigma_2 (A) +   2 \gamma_1 | \eta|^2   &\geq &  R \left(
- \varphi_1 ( t ) -   C \delta^{\frac{1}{2}}   t^{- 1}
\right) \geq  C \, t^{- \frac{2}{s}}  \left(
- t^{1- \frac{2}{s}}   -  
\delta^{\frac{1}{2}}  t^{- 1} \right)
\\[6pt]  \Longrightarrow \sigma_2 (A) &\geq&  - 2 \gamma_1 | \eta |^2 \, - \, C_3 \,
t^{1 - \frac{4}{s}} \, - \,
C_3 \, \delta^{ \frac{1}{2}} \,
t^{-(1+ \frac{2}{s})} ,\,\,\, \ t \leq T_1 .
\end{eqnarray*}
Recall that $| \eta|^2 = e^{-4 (v+w)} | \eta|^2_0 \, \geq \,
C ( g_0 ) > 0  $.
Thus, there is a constant $C_4 = C_4 ( g_0 ) > 0 $ such that
\begin{equation}
\sigma_2 (A ) \geq C_4 - \, C_3 \, t^{1- \frac{4}{s}} \, - \,
C_3 \, \delta^{\frac{1}{2}} \, t^{-(1 + \frac{2}{s})} , \,\,\,\ t \, \leq \, T_1
.
\end{equation}
Let $\widehat{t}_0 \, = \, \widehat{t}_0 ( g_0)$ satisfy
$$
C_3 \, \widehat{t}^{1-\frac{4}{s}}_0 \, = \, \frac{1}{4} \, C_4
$$
and define $t_0 \, = \, \min \left\{ T_1 , \widehat{t}_0 \right\}$.  Then
$t_0$ satisfies
\begin{equation}
C_3 \, t^{1- \frac{4}{s}}_0 \, \leq \, \frac{1}{4} \, C_4
\end{equation}
because $s > 4$.

It now follows from (7.34) and (7.35) that the metric $h = h (t_0 ,
\cdot )$ satisfies
$$
\sigma_2 ( A_h) \, \geq \, \frac{3}{4} \, C_4 \, - \, C_3 \,
\delta^{\frac{1}{2}} \, t_0^{ -(1+ \frac{2}{s})} .
$$
\noindent
Therefore, once $\delta < \delta_0 \, = \, \left( C_4 t_0^{1+\frac{2}{s}}
\bigg/ 4 C_3 \right)^2$,
\vglue6pt
\hfill ${\displaystyle
\sigma_2 (A_h) \, \geq \, \frac{1}{2} \, C_4 \, > \, 0 .
}$
\enddemo 

\section{Examples}

In this section we consider the class of 4-manifolds that admits a
positive conformal structure (i.e., the Yamabe constant of the conformal
class is positive) satisfying the condition $\int \sigma_2 (A) dv > 0$.
As a consequence of our result these manifolds carry metrics of positive
Ricci curvature, hence their \pagebreak fundamental group must be finite.  We consider
simply connected manifolds satisfying these two conditions. The
homeomorphism classification of simply connected 4-manifolds proceeds
according to the algebraic classification of the intersection form (see
[DK]). There are two families of possible quadratic forms according to
whether the manifold carries a spin structure.  For the nonspin case, the
quadratic form is of odd type and is of the form
\setcounter{equation}{-0}
\begin{equation}
{\underbrace{1 \, \oplus \, 1 \, \oplus \, \cdots \, \oplus 1}_{k}} \, \oplus
\, {\underbrace{- 1 \, \oplus
\, \cdots \, \oplus \, - 1}_\ell}
\end{equation}
\noindent
where $k \ne \ell$.  By reversing orientation if necessary we may assume
$k > \ell$.  $M$ is then homeomorphic to  
$k ( \Bbb {C P}^2) \, \# \, \ell \, \left(
\overline{ \Bbb {C P}^2} \right)$. 
\noindent
The index formula
$$
12\pi^2\tau ( M ) \, = \, \frac{1}{4} \, \int \left( | W^+|^2 \, - \, | W^-|^2 \right) $$
\noindent
and the Gauss Bonnet formula
$$
8\pi^2\chi(M)=\frac{1}{4} \int \left( |W^+|^2+|W^-|^2\right)+\int \sigma_2 $$
\noindent
combine to give
\begin{equation}
\left\{
\begin{array}{l}
4\pi^2(2 \chi \, + \, 3 \tau )\, = \, \frac{1}{2} \, \int   | W^+ |^2 \, + \, \int\sigma_2  , 
 \\[5pt]
4\pi^2(2 \chi \, - \, 3 \tau )\, = \, \frac{1}{2} \, \int   | W^-|^2 \, + \,\int \sigma_2  .
\end{array} 
\right.
\end{equation}
\noindent
Since $\chi (M) = k + \ell + 2$, $\tau (M) = k - \ell$. In order that
 $\int _M\sigma_2>0$ for $M = k ( \Bbb {C P}^2   ) \, \# \, \ell
{\overline{( \Bbb {C P}^2  )}}$, 
it is necessary that
$4 + 5 \ell > k$.  In particular, when $\ell = 0$, $k < 4$.  In fact
Lebrun et al.\ [LNN] show  that $k (  \Bbb {C P}^2 )$ for $k \leq 3$ admits a
conformal structure which is self-dual and satisfies $\frac{R}{2
\sqrt{3}} > |E|$.  When $\ell = 1$, $k < 9$, and each such topology is
homeomorphic to a K\"ahler metric with positive first Chern class.
 When $\ell > 1$,
the constraint $2 \chi - 3 \tau = 4 - k + 5 \ell > 0$ is satisfied by a
large number of manifolds; however it is not clear whether all such
topologies admit a positive conformal class for which $\int \sigma_2 dv >
0$.

In case $M$ is spin, the quadratic form is of even-type and hence of the
form 
\begin{equation}
k \, {0 \ 1 \choose 1 \ 0} \, \oplus  \, 2 \ell \, E_8
\end{equation} 
\noindent
where $E_8$ is the matrix corresponding to the Dynkin diagram of the
exceptional Lie group $E_8$.  The vanishing theorem of Lichnerowicz
requires the $\hat{A}$ genus to vanish; hence the signature is zero.
This means no $E_8$ component appears in the quadratic form, so that $M$
is homeomorphic to $k ( S^2 \times S^2)$.

To summarize the above discussion: the simply connected 4-manifolds that
admit a positive conformal structure with $\int \sigma_2 dv > 0$ must be
homeomorphic to
$
k (  \Bbb {C P}^2 ) \, \# \, \ell \, {\overline{( \Bbb {C P}^2 )}}$ or 
$k ( S^2 \times S^2 ).$ 

We remark that Sha-Yang [ShYa] have constructed metrics of positive Ricci
curvature on $k (  \Bbb {C P}^2 ) \, \# \, \ell {\overline{(
  \Bbb {C P}^2 )}}$ and $k ( S^2 \times S^2 )$ without constraints on $k
, \ell$.  Thus, the class of 4-manifolds admitting metrics with 
positive $\sigma_2$
 are necessarily a proper subset of those admitting positive
Ricci curvature metrics. 

In the following we give a construction: 
\proclaim{Proposition}
Given two positive conformal structures {\rm (}$M^\prime , g^\prime )$ and
$(M^{\prime\prime}, g^{\prime\prime} )${\rm ,} there exists a positive
conformal structure on $M_1 \# M_2$ so that
$$\int_{M^\prime} 
| W_+|^2 \, + \, \int_{M^{\prime \prime}} | W_+|^2  \, - \,  \int_{M^\prime \# M^{\prime\prime}} |W_+|^2 $$
is arbitrarily small{\rm .}
\endproclaim 

{\it  Remark}. This fact is probably well known in
view of the work of Schoen-Yau [ScY] and Gromov-Lawson [GL].  However,
we cannot find it in the literature and so will provide an argument.

\demo{Proof}  Take a geodesic coordinate system
$B^\prime$ centered at some $p^\prime \in M^\prime$ in which
$g^\prime_{ij} ( x ) = \delta_{ij} + G^\prime_{imjn} x_m x_n + O (
| x|^3 )$ and likewise for 
$B^{\prime\prime}$, $p^{\prime\prime}$ and
$g^{\prime\prime}_{ij} = \delta_{ij} + G^{\prime\prime}_{imjn} x_m
x_n + O ( | x|^3)$.  Subject $B^\prime \setminus \{ p^\prime \} $ to the
conformal change of metric
\begin{equation}
\tilde{g}^\prime \, = \,
|x|^{-2} \, g^\prime_{ij} \, d x^i dx^j \, = \,
dt^2 \, + \,
\left(
\delta_{ij} \, + \,
G^\prime_{ijmn}\sigma_m \sigma_n \, e^{2 t} \, + \,
O (e^{3t} ) \right) \, d \sigma^i d \sigma^j,
\end{equation}
\noindent
where $t = log |x| $, $ \sigma = \frac{x}{|x|}$.  To compute the scalar
curvature of $\tilde{g}^\prime$, we notice that
\begin{eqnarray} \noalign{\vskip-8pt}
&&\\
- 6 \Delta ( | x |^{-1}) \, + \, R | x |^{-2} & = & 
( - 6 \partial^2_r \, - 6H \partial_r ) \, |x|^{-1} \, + \, R |x|^{-1}
\nonumber \\
& = &- 12 | x |^{-3} \, + \, 6 \left(O( |x|)
 \, + \,
\frac{3}{|x|} \right) 
\, | x |^{-2} \, + \, R |x|^{-1} \nonumber \\
& = & ( 6 \, + \, O ( |x|)) \, | x |^{-3},\nonumber
\end{eqnarray} 
where $H$ is the mean curvature.
So the scalar curvature is
\begin{equation}
 \tilde{R} \, = \, 6 \, + \,
O (| x|).
\end{equation}

\noindent
On the annuli $ \{
e^{- r} < | x | < e^{-s} \}$ and 
$\{ e^{-r} < | x^\prime | < e^{ - s} \}$, 
introduce the cylindrical coordinates
$- r < t < -s $ 
and 
$\sigma \in S^3$ for 
$B_1$ and $- r < t^\prime < - s$ 
and $\sigma^\prime \in S^3$ for 
$B^\prime_1$; where 
$t = \log | x |$, $t^\prime = \log |x^\prime |$,
$\sigma = \frac{x}{|x|}$, 
$\sigma^\prime = \frac{x^\prime}{|x^\prime|}.$
We make the identification 
$( t , \sigma ) \sim ( t^\prime,
\sigma^\prime)$ 
if $t + r = -(s + t^\prime )$ and $\sigma = \sigma^\prime$, 
to form the connected sum 
$( B_1 \setminus \{ |x| < e^{- r} \} ) \cup ( B_1^\prime \setminus 
\{ | x^\prime | < e^{- r} \} )$.  
Let $\rho$ be a smooth function on $[ 0 , 1]$,
$\rho ( 0 ) = 0$ and 
$\rho (1) = 1$ ,
$\rho^\prime (0) = \rho^\prime (1) = 0$.  
Then set 
$\varphi (t) = \rho \left(
\frac{t + r}{r - s} \right)$ 
and define the gluing metric\pagebreak

\begin{eqnarray*}
\noalign{\vskip-18pt}
h_{ij} \, dx^i \, dx^j & = &
(\varphi ( t ) \tilde{g}^\prime_{ij} \, ( t , \sigma ) \, + \,
( 1 - \varphi ) (t) \, \tilde{g}^{\prime\prime}_{ij} \, (t^\prime ,
\sigma ))  dx^i \, dx^j 
\nonumber \\
& = & dt^2 \, + \,
( \delta_{ij} + ( \varphi G^\prime_{ijmn}\sigma_m \sigma_n \, + \,
( 1 - \varphi) \, G^{\prime\prime}_{i j m n }\sigma_m \sigma_n ) \, e^{2t} \\
&& + \,
O (e^{3t} ) ) \, d \sigma^i \, d \sigma^j.
\end{eqnarray*}
 
Over the annuli $\{e^{-s} <|x|<1\}$ and $\{e^{-s}<|x'|<1\}$ we introduce the conformal metric $g=e^{2t+v(t)}\tilde{g}'$ to
smoothly join the metric $g'$ to $\tilde{g}'$ (respectively $g''$ to $\tilde{g}'')$, while keeping the scalar curvature positive.

Observe that for functions $f$ of the $t$ variable we have $\tilde{\Delta}'=\partial^2_t +\tilde{H}'\partial_t$, where
$\tilde{H}'$ is the mean curvature of the $t$ slice $\{t\}\times S^3$.  The scalar curvature $R_g$ is then given by
\begin{eqnarray*}
R_g e^{3(t+v)} &=& - 6\Delta_{\tilde{g}'} e^{t+v}+R_{\tilde{g}'}e^{t+v}\\
&=& R_{g'} e^{3t} e^v-12e^{t+v} \partial_t v-6\left(v''+\tilde{H}'v'+|v'|^2\right) e^{t+v}\\
&=& \left\{ R_{g'} e^{2t}-12v' - 6\left( v''+|v'|^2 +o(|v'|)\right)\right\}e^{t+v}\;.\end{eqnarray*}
It is thus clear that if $s$ is taken sufficiently large, we can find a function $v$ defined on $[-s,0]$ which agrees with the
constant zero function to second order at $t=0$ and agrees with the linear function $-t$ up to second order at $t=-s$, while
keeping the scalar curvature positive.  Similarly we can join the metrics $g''$ to $\tilde{g}''$ over the annulus
$\{e^{-s}<|x'|<1\}$.
   For the
gluing region, write
\begin{equation}
h_{ij} \, = \,
(h_0 )_{ij} \, + \,
e^{2t} \left( \varphi G^\prime_{imjn} \, + \,
( 1 - \varphi ) \,
G^{\prime\prime}_{imjn} \right) 
\, \sigma_m \sigma_n \, + \, O ( e^{3 t} )
\end{equation}
where $(h_0)$ = cylinder metric.  In computing $\Gamma^m_{ij}$ we will
write 
$$\Gamma^m_{ij} \, = \,
\frac{1}{2} \, h^{mn} \, ( h_{n j , i} \, + \,
h_{i n , j} \, - h_{i j , n} ) 
$$
\noindent
and will replace indices with dots to write in abbreviation. Thus for example
$G^\prime_{\cdot\cdot}$ will mean $G^\prime_{imjn}\sigma_m \sigma_n$. 
\begin{eqnarray*}
\Gamma^\cdot_{\cdot \cdot} & = & 
\frac{1}{2} \, h^{\cdot\cdot} \, 
( \partial_\cdot \, h_{\cdot\cdot} )
\\
& = & \frac{1}{2} \, h^{\cdot\cdot}  
( 
\partial_\cdot (h_0 )_{\cdot\cdot} \, + \,
( \varphi \partial_\cdot {G^\prime}_{\cdot\cdot} \, + \, (1 - \varphi) 
 \, \partial_\cdot
{G^{\prime\prime}}_{\cdot\cdot} ) e^{2t} 
\, + \, \partial_\cdot \varphi ( G_{\cdot\cdot} ) e^{2t} \,\\
& + & \, \hbox{(lower
 order  terms )}.
\end{eqnarray*} 
Observe that when the differentiation falls on $\varphi$,  $ |
\nabla_\cdot \varphi | = O \left( \frac{1}{r-s} \right)$; hence for $r-s$
large it is of lower order.  Thus we write 
$$
\Gamma^\cdot_{\cdot\cdot} \, = \,
( \Gamma_0 )^{\cdot}_{\cdot\cdot} \, + \,
\frac{1}{2} \, h^{\cdot\cdot} \left[
\varphi \partial_\cdot ({ G^\prime}_{\cdot\cdot} e^{2t} ) \, + \,
( 1 - \varphi) \,
\partial_\cdot ( {G^{\prime\prime}}_{\cdot\cdot} e^{2t} ) \, + \,
O (e^{3t} ) \right].
$$
\noindent
Hence
\begin{equation}
\Gamma^\cdot_{\cdot\cdot} \, = \,
( \Gamma_0 )^\cdot_{\cdot\cdot} \, + \,
O ( e^{2t}) \, + \,
O \left(
\frac{1}{r-s} \right) \, e^{2t} \, + \, O ( e^{3t} ).
\end{equation} 
\noindent
Likewise
\begin{eqnarray}
\qquad \quad\Gamma^\cdot_{\cdot\cdot} \, \Gamma^\cdot_{\cdot\cdot}  
& = & \frac{1}{2} \, h^{\cdot\cdot} \left\{
\partial_\cdot h_{0\cdot\cdot} \, + \, \partial_\cdot (e^{2t} 
G^\prime_{\cdot\cdot} ) \, \varphi
\, + \, \partial_\cdot (e^{2t} G^{\prime\prime}_{\cdot\cdot}) \, ( 1 - \varphi ) \, +
\,( {\rm l.o.t.}) \right\}  \\
&& \frac{1}{2} \, h^{\cdot\cdot} \left\{
\partial_\cdot  h_{0\cdot\cdot} \, + \,
\partial_\cdot ( e^{2t} G^\prime_{\cdot\cdot} ) \, \varphi \, 
+ \, \partial_\cdot (
e^{2t} G^{\prime\prime}_{\cdot\cdot}) \, (1 - \varphi) \, 
+ \,( {\rm l.o.t.}) \right\}
\nonumber \\
& = & ( \Gamma_0)^\cdot_{\cdot\cdot} 
(\Gamma_0)^\cdot_{\cdot\cdot} \, + \,
O ( e^{2t} ) \, + \,
O \left( e^{2t} \, \frac{1}{r-s} \right) \, + \,({\rm l.o.t.}) \ . \nonumber
\end{eqnarray}
Hence, 
\begin{equation}
R^n_{ijm} \, = \,
\r{R}^{\lower6pt\hbox{$\scriptstyle n$}}_{ijm} \, + \,
O (e^{2t}) \, + \,
O (e^{2t} \, \frac{1}{r-s} ) \, + \, ({\rm l.o.t.}) 
\end{equation}
where $\r{R}^{\lower6pt\hbox{$\scriptstyle n$}}_{ijm}$ is the curvature tensor of the cylinder.  This
shows  
\begin{equation}
R \, = \,R_0 \, + \, O (e^{2t})
\end{equation}
and
\begin{equation}
|W|^2 \, = \,
|W_0|^2 \, + \,
O (e^{2t} ) \, = \,
O (e^{2t} ) .
\end{equation}
\noindent
Upon integration, the length of the cylinder being $r-s$, we find
\begin{equation}
\int_{- r < t < - s} \, | W^+|^2 \, dv \, \leq \, O ( |r-s| ) \,
( e^{-4r} ).
\end{equation}
\vglue-20pt
\enddemo
\vglue12pt

We can apply the gluing construction above to show that $M= 2 ( S^2
\times S^2 )$ admits a conformal structure satisfying
$$
\frac{1}{4}\int_{M} | W^+|^2 \, \leq \, \frac{64 \pi^2}{3} \, + \, \varepsilon.
$$
\noindent
Hence
$$
\frac{1}{4 \pi^2} \, {\displaystyle{\int_M }}\, \sigma_2 \, = \,
12 - \, \frac{1}{8 \pi^2} \, \int  |W^+|^2 \, \geq \,
12 \, - \frac{32}{3} \, - \, \varepsilon > 0 .
$$

\vglue8pt
\demo{{R}emark}  This calculation becomes critical
at $3 (S^2 \times S^2)$.  

Similarly $\Bbb {C P}^2 \, \# 2 ( S^2
\times S^2 )$ admits a conformal structure satisfying 
$$
\frac{1}{4}\int |W^+|^2 \,
\leq \, \frac{64 \pi^2}{3} \, + \, 12 \pi^2 \, + \, \varepsilon ,
$$
\noindent
so that
$$
\frac{1}{4 \pi^2} \, \int \sigma_2 \, \geq \,
17 \, - \, \frac{50}{3} \, - \, \varepsilon \, > \, 0 .$$
\enddemo

\demo{Remarks added in proof}  

1.  The regularity of general weak solutions to critical exponential variational equations such as Euler equations of the
functional $F_\delta$ appearing in this paper for $\delta\ne 0$ has been established in the article of Uhlenbeck and
Viaclovsky ([UV]).  Thus in the proof of Proposition 4.3, one may quote this regularity result to verify smoothness of
solutions.

2.  Since the submission of this paper, we were able to use the techniques of this paper to establish a conformally invariant
sphere theorem for 4-manifolds with positive Yamabe invariant satisfying an equality involving the Euler number and the
$L^2$ integral of the Weyl tensor.  We will address this result in another publication.
\AuthorRefNames [CGY-1]

\bye
\bigskip

\end{document}